\documentclass[leqno,12pt]{article}  
\usepackage{amsmath}
\usepackage{amssymb}

\usepackage[german,english]{babel}
\usepackage{amsthm}

\title{Supersingular Hecke modules as Galois representations}

\author{\textsc{Elmar Grosse-Kl\"onne}}
\date{}

\theoremstyle{plain} 
\newtheorem{satz}{Theorem}[section]  
\newtheorem{kor}[satz]{Corollary}  
\newtheorem{lem}[satz]{Lemma}  
\newtheorem{pro}[satz]{Proposition}  
 
\theoremstyle{remark}

\theoremstyle{definition}

\newcommand{\C}{\ensuremath{\mathbb{C}}}

\newcommand{\0}{\ensuremath{\overrightarrow{0}}}

\begin{document}

\maketitle

%


\begin{abstract} Let $F$ be a local field of mixed characteristic $(0,p)$, let $k$ be a finite extension of its residue field, let ${\mathcal H}$ be the pro-$p$-Iwahori Hecke $k$-algebra attached to ${\rm GL}_{d+1}(F)$ for some $d\ge1$. We construct an exact and fully faithful functor from the category of supersingular ${\mathcal H}$-modules to the category of ${\rm Gal}(\overline{F}/F)$-representations over $k$. More generally, for a certain $k$-algebra ${\mathcal H}^{\sharp}$ surjecting onto ${\mathcal H}$ we define the notion of $\sharp$-supersingular modules and construct an exact and fully faithful functor from the category of $\sharp$-supersingular ${\mathcal H}^{\sharp}$-modules to the category of ${\rm Gal}(\overline{F}/F)$-representations over $k$.
\end{abstract}

\tableofcontents

\vspace{0.6cm}

\begin{center}{\Large{\bf Introduction}}\end{center} 

Let $F$ be a local field of mixed characteristic $(0,p)$, let $\pi\in {\mathcal O}_F$ be a uniformizer, let $k$ be a finite extension of the residue field ${\mathbb F}_q$ of $F$. Let $d\in{\mathbb N}$. An important line of current research in number theory is concerned with relating smooth representations of $G={\rm GL}_{d+1}(F)$ over $k$ with finite dimensional representations of ${\rm Gal}(\overline{F}/F)$ over $k$.

At present, the smooth representation theory of $G$ is understood only up to identifying, constructing and describing the still elusive {\it supercuspidal} representations of $G$, or equivalently, the {\it supersingular} representations of $G$. An important role in better understanding this theory is played by the module theory of the pro-$p$-Iwahori Hecke $k$-algebra ${\mathcal H}$ attached to $G$ and a pro-$p$-Iwahori subgroup $I_0$ in $G$. There is a notion of supersingularity for ${\mathcal H}$-modules which, in contrast to that of supersingularity for $G$-representations, is transparent and concrete. The notions are compatible in the following sense: At least after replacing $k$ by an algebraically closed extension field, a smooth admissible irreducible $G$-representation $V$ is supersingular if and only if its space of $I_0$-invariants $V^{I_0}$ (which carries a natural action by ${\mathcal H}$) is supersingular, if and only if $V^{I_0}$ admits a supersingular subquotient; see \cite{ollvig17}. It is true that the functor $V\mapsto V^{I_0}$ from $G$-representations to ${\mathcal H}$-modules often looses information. But the potential of taking into account also its higher derived functors, which again yield (complexes of) ${\mathcal H}$-modules, has been barely explored so far. 

The purpose of the present paper is to explain a method for converting (supersingular) ${\mathcal H}$-modules into ${\rm Gal}(\overline{F}/F)$-representations over $k$.

For $F={\mathbb Q}_p$ we had constructed in \cite{dfun} an exact functor from finite dimensional ${\mathcal H}$-modules to ${\rm Gal}(\overline{\mathbb Q}_p/{\mathbb Q}_p)$-representations over $k$. The construction was inspired by Colmez's functor from ${\rm GL}_2({\mathbb Q}_p)$-representations to ${\rm Gal}(\overline{\mathbb Q}_p/{\mathbb Q}_p)$-representations. It was geometric-combinatorial in that it invoked coefficient systems on the Bruhat Tits building of ${\rm GL}_n({\mathbb Q}_p)$. Unfortunately, we see no way to generalize this geometric-combinatorial method to arbitrary finite extensions of $F$ of ${\mathbb Q}_p$. However, when trying to extract its "algebraic essence", we found that the functor indeed admits a generalization to any $F$, albeit now taking on an entirely algebraic and concrete shape. But in fact, it is this concreteness which allows us to not only investigate its behaviour on irreducible objects, but also to prove that it accurately preserves extension structures. In this way, even for $F={\mathbb Q}_p$ we significantly improve on our previous work \cite{dfun}.  

Let ${\rm Rep}({\rm Gal}(\overline{F}/F))$ denote the category whose
objects are projective limits of finite dimensional ${\rm Gal}(\overline{F}/F)$-representations over $k$. Let ${\rm Mod}_{ss}({\mathcal H})$ denote the category of supersingular ${\mathcal H}$-modules which are inductive limits of their finite dimensional submodules.\\

{\bf Theorem A:} {\it There is an exact and fully faithful functor$${\rm Mod}_{ss}({\mathcal H})\longrightarrow {\rm Rep}({\rm Gal}(\overline{F}/F)),\quad\quad M\mapsto V(M).$$We have ${\rm dim}_k(M)={\rm dim}_k(V(M))$ for any $M\in{\rm Mod}_{ss}({\mathcal H})$.}\\

The radical elimination of the group $G$ (and its building) from our approach allows us to improve Theorem A further as follows. We construct $k$-algebras ${\mathcal H}^{\sharp\sharp}$ and ${\mathcal H}^{\sharp}$ by looking at a certain small set of distinguished generators of ${\mathcal H}$ and by relaxing resp. omitting some of the usual (braid) relations between them. In this way we get a chain of surjective $k$-algebra morphisms ${\mathcal H}^{\sharp\sharp}\to{\mathcal H}^{\sharp}\to {\mathcal H}$. There is a again a notion of supersingularity for ${\mathcal H}^{\sharp\sharp}$-modules and for ${\mathcal H}^{\sharp}$-modules (which are inductive limits of their finite dimensional submodules; we assume this for all ${\mathcal H}^{\sharp\sharp}$-, resp. ${\mathcal H}^{\sharp}$-, resp. ${\mathcal H}$-modules appearing in this paper). The simple supersingular modules are the same for ${\mathcal H}^{\sharp\sharp}$, for ${\mathcal H}^{\sharp}$ and for ${\mathcal H}$, but there are more extensions between them in the category of ${\mathcal H}^{\sharp\sharp}$-modules, resp. of ${\mathcal H}^{\sharp}$-modules, than in the category of ${\mathcal H}$-modules. A particular useful category ${\rm Mod}_{ss}({\mathcal H}^{\sharp})$ is formed by what we call $\sharp$-supersingular ${\mathcal H}^{\sharp}$-modules. It contains the category of supersingular ${\mathcal H}$-modules as a full subcategory (but is larger). Now it turns out that the above functor is actually defined on the category of supersingular ${\mathcal H}^{\sharp\sharp}$-modules, and again with ${\rm dim}_k(M)={\rm dim}_k(V(M))$ for any $M$. When restricting to ${\rm Mod}_{ss}({\mathcal H}^{\sharp})$ we furthermore get:\\

{\bf Theorem A$^{\sharp}$:} {\it There is an exact and fully faithful functor$${\rm Mod}_{ss}({\mathcal H}^{\sharp})\longrightarrow {\rm Rep}({\rm Gal}(\overline{F}/F)),\quad\quad M\mapsto V(M).$$}

We do not know if the $k$-algebra ${\mathcal H}^{\sharp}$ admits a group
theoretic interpretation, as does the double coset algebra ${\mathcal H}\cong
k[I_0\backslash G/I_0]$. However, already from the Galois representation
theoretic point of view we think that the additional effort taken in proving
Theorem A$^{\sharp}$ (rather than just Theorem A) is justified, since in this way we identify an even larger abelian subcategory of ${\rm Rep}({\rm Gal}(\overline{F}/F))$ as a (supersingular) module category of a very concretely given $k$-algebra. In fact, the additional effort is mostly notational.

We define a standard supersingular ${\mathcal H}$-module to be an ${\mathcal H}$-module induced from a supersingular character of a certain subalgebra ${\mathcal H}_{\rm aff}$ of ${\mathcal H}$ with $[{\mathcal H}:{\mathcal H}_{\rm aff}]=d+1$. Each simple supersingular ${\mathcal H}$-module is a subquotient of a standard supersingular ${\mathcal H}$-module. We also define the notion of a $(d+1)$-dimensional standard cyclic ${\rm Gal}(\overline{F}/F)$-representation; in particular, each irreducible ${\rm Gal}(\overline{F}/F)$-representation of dimension $d+1$ is a $(d+1)$-dimensional standard cyclic ${\rm Gal}(\overline{F}/F)$-representation.\\

{\bf Theorem B:} {\it The functor $M\mapsto V(M)$ induces a bijection between standard supersingular ${\mathcal H}$-modules and $(d+1)$-dimensional standard cyclic ${\rm Gal}(\overline{F}/F)$-representations. $M$ is irreducible if and only if $V(M)$ is irreducible.}\footnote{A numerical version (i.e. comparing cardinalities) of Theorem B was known for quite some time, due to work of Ollivier and Vign\'{e}ras \cite{vigneras}.} \\

However, we emphasize that it is rather the much deeper Theorem A (and A$^{\sharp}$) which proves that supersingular modules are of a strong inherent arithmetic nature. 

In subsection \ref{essim} we gather some generic statements which come close to describing the image of the functor $M\mapsto V(M)$.

Let us now indicate the main features of the construction of the functor. We fix once and for all a Lubin-Tate group for $F$. More precisely, as this simplifies many formulae, we work with the Lubin Tate group associated with the Frobenius power series $\Phi(t)=t^q+\pi t$. On the $k$-algebra $k[[t]][\varphi]$ with commutation relation $\varphi\cdot
t=t^q\cdot \varphi$ we let $\Gamma={\mathcal O}_F^{\times}$ act by $\gamma\cdot
\varphi=\gamma'\varphi$ and $\gamma\cdot
t=[\gamma]_{\Phi}(t)$, where $[\gamma]_{\Phi}(t)\in k[[t]]$ describes multiplication with $\gamma$ with respect to $\Phi$ and where $\gamma'\in k^{\times}$ means the image of $\gamma\in\Gamma$ in $k^{\times}$. We view a supersingular ${\mathcal H}^{\sharp\sharp}$-module (or ${\mathcal H}^{\sharp}$-module, or ${\mathcal H}$-module) $M$ as a $k[[t]]$-module by means of $t|_M=0$. In $k[[t]][\varphi]\otimes_{k[[t]]}M$ we then use the ${\mathcal H}$-action on $M$ to define a certain submodule $\nabla(M)$ by giving very explicitly a certain number of generators of it. This is done in such a way that $\Delta(M)=k[[t]][\varphi]\otimes_{k[[t]]}M/\nabla(M)$ naturally receives an action by $\Gamma$ and is a torsion $k[[t]]$-module. A very general construction then allows us to endow $\Delta(M)^*\otimes_{k[[t]]}k((t))$ with the structure of a $(\varphi,\Gamma)$-module over $k((t))$. The notion of a $(\varphi,\Gamma)$-module over $k((t))$ with respect to the chosen Lubin-Tate group $\Phi$ is explained in full detail in the book \cite{peterlec}, where it is also explained that this category is equivalent with the category of representations of ${\rm Gal}(\overline{F}/F)$ over $k$.

It was pointed out by C\'{e}dric P\'{e}pin that the syntax of the functor $M\mapsto V(M)$ bears strong resemblance with that of Fontaine's various functors (using "big rings"). 

One may wonder which of our results remain valid if the coefficient field $k$ is allowed to be a more general field $k$ containing ${\mathbb F}_q$, i.e. not necessarily finite. First, this finiteness is invoked for the equivalence of categories between Galois representations and $(\varphi,\Gamma)$-modules. But it is also invoked in the proofs of Proposition \ref{schutzengel} (our main result in section \ref{may2019} on recovering a supersingular ${\mathcal H}^{\sharp}$-module from subquotients) and of Theorem \ref{shmichaelmas} (on recovering $M$ from $\Delta(M)$).

In subsection \ref{monachstras} we list some automorphisms of ${\mathcal H}$ (and of ${\mathcal H}^{\sharp}$ and ${\mathcal H}^{\sharp\sharp}$). They induce autoequivalences of the category of supersingular ${\mathcal H}$-modules\footnote{but this is not so evident, if true at all, for the category of $\sharp$-supersingular ${\mathcal H}^{\sharp}$-modules}; thus, precomposing them to $M\mapsto V(M)$ we get more functors satisfying Theorems A, A$^{\sharp}$, B.  

We end this paper somewhat speculatively by discussing assignments of ${\rm Gal}(\overline{F}/F)$-representations to supersingular $G$-representations. The functor $M\mapsto V(M)$ invites us to search for meaningful assigments of (complexes of) supersingular ${\mathcal H}$-modules to supersingular $G$-representations $Y$. First we suggest studying the left derived functor of the functor taking $Y$ to the maximal supersingular ${\mathcal H}$-submodule of $Y^{I_0}$. This entails working in derived categories and appears to be the most natural approach. Nevertheless, as a variation of this theme we then suggest an {\it exact} functor from (suitably filtered) $G$-representations to supersingular ${\mathcal H}$-modules. It builds on a general procedure of turning complexes of ${\mathcal H}$-modules into new ${\mathcal H}$-modules, applied here to complexes arising from $E_1$-spectral sequences attached to the said left derived functor. 

Apparantly, the constructions and results of the present paper call for generalizations
into various directions. We mention here just the obvious question of what
happens if the pro-$p$-Iwahori
Hecke algebra ${\mathcal H}$ attached to $G={\rm
  GL}_{d+1}(F)$ is replaced by pro-$p$-Iwahori
Hecke algebras ${\mathcal H}$ attached to other $p$-adic reductive groups
$G$. In extrapolation of what we did here, the general
Langlands philosophy suggests searching for a functor from
${\mathcal H}$-modules to Galois representations such that in some way the algebraic $k$-group with
root datum dual to that of $G$ shows up on the Galois side --- just as it does here in Theorem B. In a
subsequent paper we will propose such a functor. However, in its formal shape it will {\it not} precisely specialize to
the functor discussed here if $G={\rm GL}_{d+1}(F)$\footnote{but of course, it
  will be closely related}, and Theorem A will {\it not} be a special case of what we will then prove for general $G$'s.

{\it Acknowledgments} I thank Laurent Berger, Peter Schneider and Gergely Z\'{a}br\'{a}di for
helpful discussions related to this work. I thank Marie-France Vign\'{e}ras for a very close reading of the text and for detailed suggestions for improvement. I thank the anonymous referees for their careful reading and helpful recommendations. I thank Rachel Ollivier for the
invitation to UBC Vancouver in the spring of 2017; some progress on this work was
obtained during that visit.\\


{\bf Notations:} Let $F/{\mathbb Q}_p$ be a finite field extension. Let ${\mathbb F}_q$ be the
residue field of $F$ (with $q$ elements). Let $\pi$ be a uniformizer in
${\mathcal O}_F$. Let $k$ be a finite field extension of ${\mathbb F}_q$.

As explained in \cite{peterlec} Proposition 1.3.4, attached to the Frobenius (or: Lubin-Tate) formal power series $\Phi(t)=\pi t+t^q$ is associated a commutative formal group law (the associated Lubin Tate (formal) group law) $F_{\Phi}(X,Y)$ over ${\mathcal O}_F$ such that $\Phi(t)\in{\rm End}_{{\mathcal O}_F}(F_{\Phi}(X,Y))$. There is a unique injective homomorphism of rings $${\mathcal O}_F\to {\rm End}_{{\mathcal O}_F}(F_{\Phi}(X,Y)),\quad a\mapsto [a]_{\Phi}(t)$$such that $\Phi(t)=[\pi]_{\Phi}(t)$, see \cite{peterlec} Proposition 1.3.6, where we recall that, by definition,$${\rm End}_{{\mathcal O}_F}(F_{\Phi}(X,Y))=\{h\in {\mathcal O}_F[[t]]\,;\,h(0)=0\mbox{ and }h(F_{\Phi}(X,Y))=F_{\Phi}(h(X),h(Y))\}.$$

\begin{lem}\label{dominikanerjub} Assume that $F\ne{\mathbb Q}_2$. Writing $[a]_{\Phi}(t)=at+\sum_{i\ge2}a_it^i$
  (with $a_i\in{\mathcal O}_F$), we have
$a_i=0$ whenever $i-1\notin(q-1){\mathbb N}$. If $a^{q-1}=1$ we even have $a_i=0$ for all $i\ge2$.
\end{lem}

{\sc Proof:} As $\Phi(t)=\pi t+t^q$, the power series
$[a]_{\Phi}(t)=at+\sum_{i\ge2}a_it^i$ is characterized by the formula$$\pi
[a]_{\Phi}(t)+([a]_{\Phi}(t))^q=[a]_{\Phi}(\pi t+t^q).$$If $a^{q-1}=1$ we see that $[a]_{\Phi}(t)=at$ satisfies this formula. Given a general $a$, consider the equalities $[a]_{\Phi}([b]_{\Phi}(t))=[b]_{\Phi}([a]_{\Phi}(t))$ for all $b\in{\mathcal O}_F$ with $b^{q-1}=1$. Since we know $[b]_{\Phi}(t)=bt$, and since $F\ne{\mathbb Q}_2$ implies the existence of primitive such $b's$ different from $1$, we indeed obtain $a_i=0$ whenever $i-1\notin(q-1){\mathbb N}$.\hfill$\Box$\\

\section{Lubin-Tate $(\varphi,\Gamma)$-modules }

In the first two subsections we transpose some constructions and results from the theory of cyclotomic $(\varphi,\Gamma)$-modules over $k$ (i.e. where $F={\mathbb Q}_p$ and where the underlying Lubin Tate group is ${\mathbb G}_m$) to the context of $(\varphi,\Gamma)$-modules over $k$ with respect to the Lubin Tate group attached to $\Phi(t)=\pi t+t^q$ (with arbitrary $F$). Namely, we define an exact functor from admissible (torsion) $k[[t]]$-modules with commuting semilinear actions by $\Gamma={\mathcal O}_F^{\times}$ and $\varphi$ to \'{e}tale $(\varphi,\Gamma)$-modules over $k$. The former category is closely related to that of $\psi$-stable lattices in \'{e}tale $(\varphi,\Gamma)$-modules ${\bf D}$, and we are lead to transpose some of Colmez's constructions in \cite{col} involving the $\psi$-stable lattices ${\bf D}^{\natural}$ and ${\bf D}^{\sharp}$ to our context. One difference is that in our context the $\psi$-operator on
  $k((t))$ does not satisfy $\psi(1)=1$, but this necessitates only minor modifications. 

We then identify a category of admissible (torsion) $k[[t]]$-modules with actions by $\Gamma$ and $\varphi$ on which the above functor is fully faithful.

\label{muluta}

\subsection{$(\varphi,\Gamma)$-modules and torsion $k[[t]]$-modules}

Put $\Phi(t)=\pi t+t^q$. Put $\Gamma={\mathcal O}_F^{\times}$. The formula $\gamma\cdot
t=[\gamma]_{\Phi}(t)$ with $\gamma\in\Gamma$ defines an action of $\Gamma$ by $k$-algebra automorphisms on
$k[[t]]$ and on $k((t))$. Consider the $k$-algebra$${\mathfrak O}=k[[t]][\varphi,\Gamma]$$ with
commutation rules given by$$\gamma
\varphi=\varphi\gamma,\quad\quad \gamma
t=[\gamma]_{\Phi}(t)\gamma,\quad\quad\varphi
t=t^q\varphi$$for $\gamma\in\Gamma$. (Here we read $[\gamma]_{\Phi}(t)\gamma=([\gamma]_{\Phi}(t))\gamma$.)\footnote{As $t^q=\Phi(t)=[\pi]_{\Phi}(t)$ in $k[[t]]$ one may also think of ${\mathfrak O}$ as ${\mathfrak O}=k[[t]][{\mathcal O}_F-\{0\}]$ with commutation rules $a t=[a]_{\Phi}(t) a$ for all $a\in {\mathcal O}_F-\{0\}$.}\\

{\bf Definition:} A $\psi$-operator on $k[[t]]$ is a $k$-linear map $\psi:k[[t]]\to
k[[t]]$ such that
$\psi(\gamma\cdot t)=\gamma\cdot (\psi(t))$
for all $\gamma\in\Gamma$ and such
that the following holds true\footnote{We
  do not require $\psi(1)=1$.}: If we view $\varphi$ as acting on
$k[[t]]$, then\begin{gather}\psi(\varphi(a)x)=a\psi(x)\quad\quad\mbox{ for }a,x\in k[[t]].\label{psichar}\end{gather}

\begin{lem}\label{erleicht} There is a surjective $\psi$-operator on
  $k[[t]]$ which extends to
a surjective $k$-linear operator $\psi=\psi_{k((t))}$ on $k((t))$ satisfying formula (\ref{psichar}) analogously.

We may choose $\psi_{k((t))}$ on $k((t))$ such that for $m\in{\mathbb Z}$ and
$0\le i\le q-1$ we have\footnote{Notice that $\frac{q}{\pi}=0$ (in $k$) if $F\ne{\mathbb Q}_p$.}\begin{gather}\psi_{k((t))}(t^{mq+i})=\left\{\begin{array}{l@{\quad:\quad}l}\frac{q}{\pi}t^m& i=0\\0
&  1\le i\le q-2\\t^m& i=q-1\end{array}\right.\label{psieqqp}.\end{gather}
\end{lem}

{\sc Proof:} This is explained in \cite{multivar}; it relies on \cite{schven} section
3.\hfill$\Box$\\

In the following we fix the surjective $\psi$-operator $\psi$ on $k[[t]]$ satisfying formula
(\ref{psieqqp}). We extend it to a $k$-linear operator $\psi=\psi_{k((t))}$ on $k((t))$ as in Lemma \ref{erleicht}.\\

{\bf Definition:} An \'{e}tale $(\varphi,\Gamma)$-module over
$k((t))$ is an ${\mathfrak O}\otimes_{k[[t]]}k((t))$-module ${\bf D}$
which is finite dimensional over $k((t))$ such that the
$k((t))$-linearized structure map$${\rm id}\otimes\varphi:k((t))\otimes_{\varphi,k((t))} {\bf
  D}\stackrel{\cong}{\longrightarrow}{\bf
  D}$$ is bijective. We define ${\rm Mod}^{et}(k((t)))$ to be the category of \'{e}tale $(\varphi,\Gamma)$-module over
$k((t))$.

\begin{satz}\label{sosego} (Fontaine, Kisin-Ren, Schneider) There is an equivalence between ${\rm Mod}^{et}(k((t)))$ and the category of continuous representations of ${\rm Gal}(\overline{F}/F)$ on finite dimensional $k$-vector spaces.  
\end{satz}

{\sc Proof:} For $F={\mathbb Q}_p$ and the Frobenius power series $(1+t)^p-1$ (instead of $\Phi(t)=\pi t+t^q$) this is a theorem of Fontaine, see par 1.2. in \cite{fon}. The analog of the theorem (for an arbitrary Frobenius power series) for a coefficient field of characteristic $0$ (hence not $k$) is due to Kisin and Ren, see \cite{kire}. A detailed proof of the theorem stated here can be found in Schneider's book \cite{peterlec}.\hfill$\Box$\\

{\bf Definition:} A torsion $k[[t]]$-module $\Delta$ is called admissible if
$$\Delta[t]=\{x\in\Delta\,;\,tx=0\}$$is a finite dimensional $k$-vector
space.\\

We remark that admissible $k[[t]]$-modules on which $t$ acts
surjectively are precisely the Pontrjagin duals of finitely generated torsion
free, and hence free $k[[t]]$-modules.\\

{\bf Definition:} ${\rm Mod}^{\rm ad}({\mathfrak O})$ is the category of ${\mathfrak O}$-modules which are finitely generated over $k[[t]][\varphi]$ and admissible (in particular: torsion) over
$k[[t]]$.\\

\begin{lem} The categories ${\rm Mod}^{et}(k((t)))$ and ${\rm Mod}^{\rm ad}({\mathfrak O})$ are abelian.
\end{lem}

{\sc Proof:} An ${\mathfrak O}\otimes_{k[[t]]}k((t))$-module subquotient of an
\'{e}tale $(\varphi,\Gamma)$-module is again an \'{e}tale
$(\varphi,\Gamma)$-module: To see that the \'{e}taleness condition (the
bijectivity of ${\rm id}\otimes\varphi$) is preserved under passing to such
subquotients, just notice that it is equivalent with saying that the matrix of
$\varphi$ in an arbitrary $k((t))$-basis is invertible. Thus, ${\rm
  Mod}^{et}(k((t)))$ is abelian. (Of course, one could also point to Theorem
\ref{sosego}.)

 An ${\mathfrak O}$-module subquotient of an object in ${\rm Mod}^{\rm
   ad}({\mathfrak O})$ is again an object in ${\rm Mod}^{\rm
   ad}({\mathfrak O})$: This is shown in \cite{em} Proposition 3.3. Thus, ${\rm Mod}^{\rm ad}({\mathfrak O})$ is abelian.\hfill$\Box$\\

{\bf Definition:} For a $k$-vector space $\Delta$ we write $\Delta^*={\rm
  Hom}_k(\Delta,k)$ (algebraic dual). For a $k[[t]]$-module $\Delta$ we endow $\Delta^*$ with a $k[[t]]$-action by
putting$$(S\cdot f)(\delta)=f(S\delta)$$for $S\in k[[t]]$, $f\in\Delta^*$,
$\delta\in\Delta$. If $\Delta$ even carries a $k[[t]][\Gamma]$-module structure then also $\Delta^*$ receives one, with $\gamma\in\Gamma$ acting as$$(\gamma\cdot
f)(\delta)=f(\gamma^{-1}\delta)$$for $\gamma\in\Gamma$, $f\in\Delta^*$, $\delta\in\Delta$. 

\begin{pro}\label{nopsi} For $\Delta\in {\rm Mod}^{\rm ad}({\mathfrak O})$
  there is a natural structure of an \'{e}tale $(\varphi,\Gamma)$-module on
  $\Delta^*\otimes_{k[[t]]}k((t))$. The contravariant functor\begin{gather}{\rm Mod}^{\rm ad}({\mathfrak O})\longrightarrow {\rm Mod}^{et}(k((t))),\quad \Delta\mapsto
  \Delta^*\otimes_{k[[t]]}k((t))\label{wienschluss}\end{gather}is exact.
\end{pro}

{\sc Proof:} The map ${\rm
    id}\otimes\varphi:k[[t]]\otimes_{\varphi,k[[t]]}\Delta\to\Delta$ gives rise to the $k[[t]]$-linear
map\begin{gather}\Delta^*\stackrel{({\rm
    id}\otimes\varphi)^*}{\longrightarrow}(k[[t]]\otimes_{\varphi,k[[t]]}\Delta)^*.\label{multvarquer1}\end{gather}On
the other hand, we have the $k[[t]]$-linear
map\begin{gather}k[[t]]\otimes_{\varphi,k[[t]]}(\Delta^*)\longrightarrow
(k[[t]]\otimes_{\varphi,k[[t]]}\Delta)^*\label{multvarquer2}\\a\otimes\ell\mapsto[b\otimes
  x\mapsto\ell(\psi(ab)x)].\notag\end{gather}It is shown in \cite{multivar}
that the respective base extended maps
(\ref{multvarquer1})$\otimes_{k[[t]]}k((t))$ and
(\ref{multvarquer2})$\otimes_{k[[t]]}k((t))$ are bijective. Composing (\ref{multvarquer2})$\otimes_{k[[t]]}k((t))$
with the inverse of (\ref{multvarquer1})$\otimes_{k[[t]]}k((t))$ thus yields a $k((t))$-linear
isomorphism$$k((t))\otimes_{\varphi,k((t))}(\Delta^*\otimes_{k[[t]]}k((t)))=k((t))\otimes_{\varphi,k[[t]]}(\Delta^*)\longrightarrow
\Delta^*\otimes_{k[[t]]}k((t))$$and hence the desired $\varphi$-operator on
$\Delta^*\otimes_{k[[t]]}k((t))$. The exactness of
$\Delta\mapsto\Delta^*\otimes_{k[[t]]}k((t))$ follows from the exactness of
taking duals and of applying $(.)\otimes_{k[[t]]}k((t))$.\hfill$\Box$\\
  
\subsection{$\psi$-stable lattices in $(\varphi,\Gamma)$-modules}

\begin{lem}\label{vigiltaufe} Let ${\bf D}\in{\rm Mod}^{et}(k((t)))$. There is a natural additive operator $\psi:{\bf D}\to {\bf D}$ satisfying $$\psi(a\varphi(x))=\psi(a)x\quad\quad\mbox{and}\quad\quad\psi(\varphi(a)x)=a\psi(x)$$for
all $a\in
k((t))$ and all $x\in {\bf
  D}$, and commuting with the action of $\Gamma$.
\end{lem} 

{\sc Proof:} We define the composed map$$\psi:{\bf D}\longrightarrow
k((t))\otimes_{\varphi,k((t))}{\bf D}\longrightarrow{\bf D}$$where the first
arrow is the inverse of the structure isomorphism ${\rm id}\otimes\varphi$,
and where the second arrow is given by $a\otimes x\mapsto \psi(a)x$. By
construction, it
satisfies $\psi(a\varphi(x))=\psi(a)x$. To see $\psi(\varphi(a)x)=a\psi(x)$
observe that by assumption we may write $x=\sum_ia_i\varphi(d_i)$ with $d_i\in
{\bf D}$ and $a_i\in k((t))$. We then
compute$$\psi(\varphi(a)x)=\sum_i\psi(\varphi(a)a_i\varphi(d_i))=\sum_i\psi(\varphi(a)a_i)d_i$$$$=a\sum_i\psi(a_i)d_i=a\sum_i\psi(a_i\varphi(d_i))=a\psi(x).$$Finally, let
$\gamma\in\Gamma$. As $\gamma$ and $\varphi$ commute on $k[[t]]$, and as
$\Gamma$ acts semilinearly on ${\bf D}$, the additive map
$$k((t))\otimes_{\varphi,k((t))}{\bf D}\to k((t))\otimes_{\varphi,k((t))}{\bf
  D},\quad a\otimes d\mapsto \gamma(a)\otimes\gamma(b)$$ is the map corresponding
to $\gamma$ on ${\bf D}$ under the isomorphism ${\rm id}\otimes\varphi$, and under
$a\otimes x\mapsto \psi(a)x$ it commutes with $\gamma$ on ${\bf D}$ since $\gamma$
and $\psi$ commute on $k((t))$.\hfill$\Box$\\

In the following, by a lattice in a $k((t))$-vector space ${\bf D}$ we mean a free $k[[t]]$-submodule containing a $k((t))$-basis of ${\bf D}$.

\begin{lem}\label{hitze3} Let ${\bf D}\in{\rm Mod}^{et}(k((t)))$ and let $D$ be a lattice in (the $k((t))$-vector space underlying) ${\bf D}$. Let $\psi:{\bf D}\to {\bf D}$ be the operator
  constructed in Lemma \ref{vigiltaufe}.

(a) $\psi(D)$ is a $k[[t]]$-module.

(b) If $\varphi(D)\subset D$ then $D\subset\psi(D)$.

(c) If $D\subset k[[t]]\cdot \varphi(D)$ then $\psi(D)\subset D$.

(d) If $\psi(D)\subset D$ then $\psi(t^{-1}D)\subset t^{-1}D$, and for each $x\in {\bf D}$ there is some $n(x)\in{\mathbb N}$ such that for all $n\ge n(x)$ we have $\psi^n(x)\in t^{-1}D$.

\end{lem}

{\sc Proof:} (a) Use $\psi(\varphi(a)x)=a\psi(x)$ for $a\in
k((t))$ and $x\in {\bf
  D}$.

(b) Choose $a\in k[[t]]$ with $\psi(a)=1$. For $d\in D$ we have $d=\psi(a\varphi(d))$ which belongs to $\psi(D)$ since $\varphi(D)\subset D$.

(c) Let $d\in D$. By assumption there are $e_i\in D$ and $a_i\in k[[t]]$ with $d=\sum_ia_i\varphi(e_i)$, hence $\psi(d)=\sum_i\psi(a_i)e_i\in D$.

(d) For $i\ge1$ we
have \begin{gather}\psi(\varphi^i(t^{-1})D)\subset\varphi^{i-1}(t^{-1})\psi(D)\subset\varphi^{i-1}(t^{-1})D\label{colmepeda}\end{gather}where
the second inclusion uses the assumption. From $\varphi(t^{-1})=t^{-q}$ we get $$\psi(t^{-1}D)\subset\psi(\varphi(t^{-1})D)\subset t^{-1}D.$$Moreover,
if $n(x)\in{\mathbb N}$ is such that $x\in\varphi^n(t^{-1})D$ for $n\ge n(x)$,
then iterated application of formula (\ref{colmepeda}) shows $$\psi^n(x)\in
\psi^n(\varphi^n(t^{-1})D)\subset
\psi^{n-1}(\varphi^{n-1}(t^{-1})D)\subset\ldots\subset t^{-1}D$$ for $n\ge n(x)$.\hfill$\Box$\\

\begin{lem} (a) There are lattices $D_0$, $D_1$ in ${\bf D}$ with $$\varphi(D_0)\subset t D_0\subset D_0\subset D_1\subset k[[t]]\cdot \varphi(D_1).$$(b) For $D_0$, $D_1$ as in (a) and for $n\ge0$ we have $\psi^n(D_0)\subset\psi^{n+1}(D_0)\subset D_1$.
\end{lem}

{\sc Proof:} (a) This is a (simplified) subclaim in the proof of Lemma
2.2.10 in \cite{peterlec} (which follows \cite{col} Lemme II 2.3). One
proceeds as follows. Let $d_1,\ldots, d_r$ be a $k((t))$-basis of ${\bf
  D}$. Then also $\varphi(d_1),\ldots, \varphi(d_r)$ is $k((t))$-basis of ${\bf
  D}$. We therefore find $\tilde{f}_{ij}, \tilde{g}_{ij}\in
k((t))$ with $\varphi(d_j)=\sum_{i=1}^r\tilde{f}_{ij}d_i$ and
$d_j=\sum_{i=1}^r\tilde{g}_{ij}\varphi(d_i)$, for any $1\le j\le r$. Choose some $n\ge0$ with $t^{n(q-1)}\tilde{f}_{ij}\in tk[[t]]$ and $t^{n(q-1)}\tilde{g}_{ij}\in
tk[[t]]$ for all $i, j$. Then $D_0=\sum_{i=1}^rt^nk[[t]]d_i$ and
$D_1=\sum_{i=1}^rt^{-n}k[[t]]d_i$ work as desired.

(b) Choose $a\in k[[t]]$ with $\psi(a)=1$. For $x\in D_0$ we have $\psi^n(x)=\psi^{n+1}(a\varphi(x))\in\psi^{n+1}(D_0)$ since $\varphi(D_0)\subset tD_0$ implies $\varphi(x)\in D_0$ and hence $a\varphi(x)\in D_0$. This shows $\psi^n(D_0)\subset\psi^{n+1}(D_0)$. As $D_0\subset D_1\subset k[[t]]\cdot\varphi(D_1)$, an induction using Lemma \ref{hitze3} (c) shows $\psi^{n+1}(D_0)\subset D_1$.\hfill$\Box$\\

\begin{pro}\label{goldabreise} There exists a unique lattice ${\bf D}^{\sharp}$ in ${\bf D}$ with $\psi({\bf D}^{\sharp})={\bf D}^{\sharp}$ and such that for each $x\in {\bf D}$ there is some $n\in{\mathbb N}$ with $\psi^n(x)\in {\bf D}^{\sharp}$.

For any lattice $D$ in ${\bf D}$ we have $\psi^n(D)\subset {\bf D}^{\sharp}$ for all $n>>0$.

For any lattice $D$ in ${\bf D}$ with $\psi(D)=D$ we have $t{\bf D}^{\sharp}\subset D\subset {\bf D}^{\sharp}$.
\end{pro} 

{\sc Proof:} Using the previous Lemmata, the proof is the same as the one given in \cite{col} Proposition II.4.2. \hfill$\Box$\\

\begin{pro}\label{wienabsage} (a) For any lattice $D$ in ${\bf D}$ contained
  in ${\bf D}^{\sharp}$ and stable under
  $\psi$ we have $\psi(D)=D$.

(b) The intersection ${\bf D}^{\natural}$ of all lattices in ${\bf D}$
  contained in ${\bf D}^{\sharp}$ and stable under
  $\psi$ is itself a lattice, and it satisfies $\psi({\bf D}^{\natural})={\bf
    D}^{\natural}$.
\end{pro} 

{\sc Proof:} (See \cite{col} Proposition II.5.11 and Corollaire II.5.12) 

(a) Since ${\bf D}^{\sharp}$ as well as $D$ and $\psi(D)$ are
lattices in ${\bf D}^{\sharp}$, both ${\bf D}^{\sharp}/D$ and ${\bf
  D}^{\sharp}/\psi(D)$ are finite dimensional $k$-vector spaces. $\psi$ induces an isomorphism $\psi({\bf D}^{\sharp})/D={\bf D}^{\sharp}/\psi(D)$ (as $\psi(D)\subset D$), hence $\psi(D)=D$.

(b) For any $D$ as in (a) we have $t{\bf D}^{\sharp}\subset D$ by what we saw
in (a) together with proposition
\ref{goldabreise}. This shows $t{\bf D}^{\sharp}\subset{\bf D}^{\natural}$, hence
${\bf D}^{\natural}$ is indeed a lattice, and $\psi({\bf D}^{\natural})={\bf
    D}^{\natural}$ follows by applying (a) once more.\hfill$\Box$\\

\begin{lem}\label{staga} ${\bf D}^{\natural}$ and ${\bf D}^{\sharp}$ are stable under the action of $\Gamma$.
\end{lem} 

{\sc Proof:} If $D$ is a lattice in ${\bf D}$, then so is $\gamma\cdot D$ for
any $\gamma\in\Gamma$. If in addition $\psi(D)\subset D$, resp. $\psi(D)= D$,
then also $\psi(\gamma\cdot D)\subset \gamma\cdot D$, resp. $\psi(\gamma\cdot
D)= \gamma\cdot D$. From these observations we immediately get
$\gamma\cdot{\bf D}^{\natural}={\bf D}^{\natural}$ and $\gamma\cdot{\bf D}^{\sharp}={\bf D}^{\sharp}$.
\hfill$\Box$\\

\begin{pro}\label{irrnopsi} The functor ${\rm Mod}^{\rm ad}({\mathfrak O})\to {\rm Mod}^{et}(k((t)))$ in Proposition \ref{nopsi} sends simple objects to simple objects.
\end{pro}

{\sc Proof:} Let $\Delta\in {\rm Mod}^{\rm ad}({\mathfrak O})$ be simple. By construction, $\psi$ on $\Delta^*\otimes_{k[[t]]}k((t))$, when
restricted to $\Delta^*$, is the adjoint of $\varphi$ on
$\Delta$. Therefore the simplicity of $\Delta$ as an ${\mathfrak O}$-module means that $\Delta^*$
admits no non-trivial $k[[t]]$-submodule stable under $\Gamma$ and $\psi$. If ${\bf D}$ is a non-zero $(\varphi,\Gamma)$-submodule of
$\Delta^*\otimes_{k[[t]]}k((t))$ then also ${\bf D}^{\natural}$ is non-zero
and stable under $\Gamma$ and $\psi$, see Proposition \ref{wienabsage} and
Lemma \ref{staga}. As ${\bf D}^{\natural}\subset
(\Delta^*\otimes_{k[[t]]}k((t)))^{\natural}\subset \Delta^*$ we get ${\bf
  D}^{\natural}=\Delta^*$ (since $\Delta^*$ is stable under $\psi$), as desired.\hfill$\Box$\\

\begin{lem}\label{goldhoch} Let $f:{\bf D}_1\to {\bf D}_2$ be a morphism in ${\rm Mod}^{et}(k((t)))$. 

(a) $f({\bf D}_1^{\sharp})\subset {\bf D}_2^{\sharp}$ and $f({\bf D}_1^{\natural})\subset {\bf D}_2^{\natural}$.

(b) If $f:{\bf D}_1\to {\bf D}_2$ is injective (resp. surjective), then so is $f:{\bf D}_1^{\sharp}\to {\bf D}_2^{\sharp}$.

(c) If $f:{\bf D}_1\to {\bf D}_2$ is injective (resp. surjective), then so is $f:{\bf D}_1^{\natural}\to {\bf D}_2^{\natural}$.

\end{lem}

{\sc Proof:} (a) $f({\bf D}_1^{\sharp})$ is a free $k[[t]]$-submodule of ${\bf D}_2$ on
which $\psi$ acts surjectively. Thus $f({\bf D}_1^{\sharp})+{\bf
  D}_2^{\sharp}$ is a lattice satisfying the defining condition for ${\bf
  D}_2^{\sharp}$ given in \ref{goldabreise}, hence $f({\bf D}_1^{\sharp})+{\bf
  D}_2^{\sharp}={\bf
  D}_2^{\sharp}$, hence $f({\bf D}_1^{\sharp})\subset {\bf
  D}_2^{\sharp}$. Next, let $D=\{x\in {\bf D}_1^{\natural}\,;\,f(x)\in {\bf
  D}_2^{\natural}\}$. It is a lattice in ${\bf D}_1$ since ${\bf
  D}_1^{\natural}$ is a lattice, $f({\bf D}_1^{\natural})\subset f({\bf
  D}_1^{\sharp})\subset {\bf D}_2^{\sharp}$ and ${\bf D}_2^{\sharp}/{\bf
  D}_2^{\natural}$ is a finite dimensional $k$-vector space. It is also stable
under $\psi$, hence contains ${\bf D}_1^{\natural}$, hence $f({\bf D}_1^{\natural})\subset {\bf D}_2^{\natural}$.

(b) and (c) If $f:{\bf D}_1\to {\bf D}_2$ is injective then obviously so are
$f:{\bf D}_1^{\sharp}\to {\bf D}_2^{\sharp}$ and $f:{\bf D}_1^{\natural}\to
{\bf D}_2^{\natural}$. If $f:{\bf D}_1\to {\bf D}_2$ is surjectice then
$f({\bf D}_1^{\natural})$ is a lattice in ${\bf D}_2$ stable under $\psi$,
hence contains ${\bf D}_2^{\natural}$. To see $f({\bf D}_1^{\sharp})={\bf
  D}_2^{\sharp}$ we proceed as in \cite{col} Proposition II.4.6 (iii). Namely, choose a lattice $D'$ in ${\bf D}_1$ with $f(D')={\bf
  D}_2^{\sharp}$. Put $D=\sum_{n\ge0}\psi^n(D')$. By construction we have $\psi(D)\subset D$
as well as $f(D)={\bf
  D}_2^{\sharp}$ (since $\psi ({\bf
  D}_2^{\sharp})={\bf
  D}_2^{\sharp}$). Proposition \ref{goldabreise} shows that $D$ is
again a lattice. Let $x\in {\bf
  D}_2^{\sharp}$. By Proposition \ref{goldabreise} we find some $n\in{\mathbb N}$ such that $\psi^n(D)\subset {\bf D}_1^{\sharp}$. For such an $n$, choose $x_n\in {\bf
  D}_2^{\sharp}$ and $\tilde{x}_n\in D$ with $\psi^n(x_n)=x$ and
$f(\tilde{x}_n)=x_n$. Put $u_n=\psi^n(\tilde{x}_n)\in {\bf D}_1^{\sharp}$. By their construction in Lemma \ref{vigiltaufe}, the operators $\psi$ on ${\bf D}_1$ and ${\bf D}_2$ commute with $f$, thus we may compute$$f(u_n)=f(\psi^n(\tilde{x}_n))=\psi^n(f(\tilde{x}_n))=\psi^n(x_n)=x.$$\hfill$\Box$\\

\begin{lem}\label{goldhochzeit} Let $0\to {\bf D}_1\to {\bf D}_2\to{\bf
    D}_3\to0$ be an exact sequence in ${\rm Mod}^{et}(k((t)))$. For each $i$
  let $D_i\subset {\bf D}_i$ be a lattice with $\psi(D_i)=D_i$, and suppose
  that the above sequence restricts to an exact
  sequence \begin{gather}0\longrightarrow {D}_1\longrightarrow
    {D}_2\longrightarrow{D}_3\longrightarrow0.\label{sommer}\end{gather}If we have
  $D_1={\bf D}_1^{\natural}$ and $D_3={\bf
    D}_3^{\natural}$, then we also have $D_2={\bf
    D}_2^{\natural}$. If we have $D_1={\bf D}_1^{\sharp}$ and $D_3={\bf
    D}_3^{\sharp}$ then we also have $D_2={\bf D}_2^{\sharp}$.   
\end{lem} 

{\sc Proof:} By Lemma \ref{goldhoch} the sequence $0\to {\bf
  D}_1^{\natural}\to {\bf D}_2^{\natural}\to{\bf D}_3^{\natural}\to0$ is exact
on the left and on the right. Comparing it with the sequence (\ref{sommer})
via ${\bf D}_1^{\natural}=D_1$, ${\bf D}_2^{\natural}\subset D_2$ and ${\bf
  D}_3^{\natural}=D_3$, we immediately get ${\bf D}_2^{\natural}=D_2$. Next,
by Lemma \ref{goldhoch} the sequence $0\to {\bf D}_1^{\sharp}\to {\bf
  D}_2^{\sharp}\to{\bf D}_3^{\sharp}\to0$ is exact on the left and on the
right. We compare it with the sequence (\ref{sommer}) via $D_1={\bf
  D}_1^{\sharp}$, $D_2\subset {\bf D}_2^{\sharp}$ and $D_3={\bf
  D}^{\sharp}_3$. We claim $$\psi({\bf D}_1\cap {\bf D}_2^{\sharp})={\bf
  D}_1\cap {\bf D}_2^{\sharp}.$$Of course, $\psi({\bf D}_1\cap {\bf D}_2^{\sharp})\subset{\bf
  D}_1\cap {\bf D}_2^{\sharp}$ is clear. To see ${\bf
  D}_1\cap {\bf D}_2^{\sharp}\subset\psi({\bf D}_1\cap {\bf D}_2^{\sharp})$
take $x\in {\bf D}_1\cap {\bf D}_2^{\sharp}$. Choose $y\in {\bf D}_2^{\sharp}$
with $\psi(y)=x$. Choose $y'\in D_2$ mapping to the same element in ${\bf
  D}^{\sharp}_3=D_3$ as $y$. We then have $\psi(y')\in D_2\cap {\bf D}_1=D_1$
and $\psi(y-y')-x\in D_1$, hence there is some $z\in D_1$ with
$\psi(z)=\psi(y-y')-x$, hence $x=\psi(y-y'-z)\in \psi({\bf D}_1\cap {\bf
  D}_2^{\sharp})$ since $y-y'\in {\bf D}_1\cap {\bf
  D}_2^{\sharp}$ and $z\in {\bf D}_1\cap {\bf
  D}_2^{\sharp}$. 

The claim is proven. By the definition of ${\bf D}_1^{\sharp}$ it implies ${\bf D}_1\cap
{\bf D}_2^{\sharp}={\bf D}_1^{\sharp}$, hence ${\bf D}_1\cap
{\bf D}_2^{\sharp}=D_1$ since $D_1={\bf D}_1^{\sharp}$. Thus, $D_2={\bf D}_2^{\sharp}$.\hfill$\Box$\\

{\bf Remark:} An \'{e}tale $\varphi$-module over
$k((t))$ is a $k[[t]][\varphi]\otimes_{k[[t]]}k((t))$-module ${\bf D}$
which is finite dimensional over $k((t))$ such that the
$k((t))$-linearized structure map ${\rm id}\otimes\varphi$ is bijective. The above theory of the operator $\psi$ and the lattices ${\bf D}^{\sharp}$ and ${\bf D}^{\natural}$ works analogously for \'{e}tale $\varphi$-modules ${\bf D}$ over
$k((t))$, i.e. the $\Gamma$-action is not really needed.

\subsection{Partial full faithfulness of $\Delta\mapsto
  \Delta^*\otimes_{k[[t]]}k((t))$}

\begin{lem} \label{abstrtor} Let $N$ be a $k$-vector space, and suppose that we are given a $k$-linear automorphism $\tau$ of $N$, a basis
  ${\mathcal N}$ of $N$, integers $0\le k_{\nu}\le q-1$ and units $\alpha_{\nu}\in k^{\times}$ for ${\nu}\in {\mathcal
  N}$. View $N$ as a $k[[t]]$-module with $t N=0$ and let $\Delta$ denote
the quotient of $k[[t]][\varphi]\otimes_{k[[t]]}N$ by the $k[[t]][\varphi]$-submodule $\nabla$
generated by the elements $$1\otimes
{\nu}+\alpha_{\nu}t^{k_{\nu}}\varphi\otimes\tau({\nu})$$with ${\nu}\in {\mathcal
  N}$. We then have:

(a) $k[[t]][\varphi]\otimes_{k[[t]]}N$ is
a torsion $k[[t]]$-module.

(b) The map $N\to \Delta[t]$ sending $n\in N$ to the class of $1\otimes n$ is
an isomorphism.

In particular, $\Delta$ is admissible if $N$ is a finite dimensional $k$-vector space. 

(c) The action of $\varphi$ on $\Delta$ is injective.
\end{lem}

{\sc Proof:} (a) As $\varphi t=t^q\varphi$ in $k[[t]][\varphi]$ we may write any element in $k[[t]][\varphi]\otimes_{k[[t]]}N$ as a
finite sum of elements of the form $a\varphi^n\otimes x$ with $a\in k[[t]]$,
$n\ge 0$ and $x\in N$. It is therefore enough to show \begin{gather}a\varphi^n\otimes x=0\quad\quad\mbox{ for each }a\in t^{q^n}k[[t]]\label{tdmneu}\end{gather}where $n\ge 0$ and $x\in N$. We may write
$a=a_0 t^{q^n}$ with $a_0\in k[[t]]$ and compute$$a\varphi^n\otimes x=a_0 t^{q^n}\varphi^n\otimes
x=a_0\varphi^nt\otimes x=0.$$ 

(b) and (c) We may
write$$k[[t]][\varphi]\otimes_{k[[t]]}N\cong \bigoplus_{\nu\in{\mathcal
    N}}\bigoplus_{i\ge 0}\bigoplus_{0\le\theta\le
  q^{i}-1}k.t^{\theta}\varphi^i\otimes \tau(\nu).$$Indeed, that
$k[[t]][\varphi]\otimes_{k[[t]]}N$ is a quotient of the right hand side
follows from formula (\ref{tdmneu}). It is in fact an isomorphic quotient since all
relations between $\varphi$ and $t$ in $k[[t]][\varphi]$ can be generated from
$\varphi t=t^q\varphi$.

Consider the three $k$-subvector spaces \begin{gather}1\otimes N=\bigoplus_{\nu\in{\mathcal N}}k\otimes
\tau(\nu)=\bigoplus_{\nu\in{\mathcal N}}k\otimes
\nu,\notag\\C=\bigoplus_{\nu\in{\mathcal N}}\bigoplus_{i>0}\bigoplus_{0\le\theta<q^{i-1}k_{\nu}}k.t^{\theta}\varphi^i\otimes \tau(\nu),\label{bohrnerv}\\\nabla=\bigoplus_{\nu\in{\mathcal
    N}}\bigoplus_{i>0}\bigoplus_{\epsilon\ge 0}k.t^{\epsilon}\varphi^{i-1}(1\otimes
{\nu}+\alpha_{\nu}t^{k_{\nu}}\varphi\otimes\tau({\nu})).\label{bohrnerv0}\end{gather}Using the formula $\varphi
t=t^q \varphi$ we see$$t^{\epsilon}\varphi^{i-1}(1\otimes
{\nu}+\alpha_{\nu}t^{k_{\nu}}\varphi\otimes\tau({\nu}))\in
k^{\times}.t^{\epsilon+q^{i-1}k_{\nu}}\varphi^i\otimes
\tau({\nu})+k[[t]]\varphi^{i-1}\otimes\nu.$$We also see that in the sum
(\ref{bohrnerv0}) all summands with
$\epsilon\ge(q-1)q^{i-1}k_{\nu}-1$ vanish. Equivalently, in the sum
(\ref{bohrnerv0}) only those summands are nonzero for which $\theta=\epsilon+q^{i-1}k_{\nu}$ satisfies $q^{i-1}k_{\nu}\le \theta\le
q^{i}-1$. Thus we
find\begin{gather}k[[t]][\varphi]\otimes_{k[[t]]}N\quad\cong\quad 1\otimes
N\quad\bigoplus\quad \nabla\quad\bigoplus\quad
C.\label{bohrnerv1}\end{gather}

Let $C'$, resp. $C''$, denote the $k$-subspace of $C$ spanned by all $t^{\theta}\varphi^i\otimes
\tau(\nu)$ with $\nu\in{\mathcal N}$, $i>1$ and
  $0\le\theta<q^{i-1}k_{\nu}$, resp. by all $t^{\theta}\varphi\otimes
\tau(\nu)$ with $\nu\in{\mathcal N}$ and
  $0\le\theta<k_{\nu}$. Then $\varphi(C)\subset C'$ and $\varphi:C\to C'$
  is injective. On the other hand, $\varphi(1\otimes
  N)\subset C''$ and $\varphi:1\otimes
  N\to C''$ is injective. Since $C'\cap C''=0$ we
  conclude that $\varphi$ acts injectively on
$\Delta$. Now consider the composed
map$$C\quad\longrightarrow\quad
k[[t]][\varphi]\otimes_{k[[t]]}N\quad\stackrel{t(.) }{\longrightarrow}\quad
k[[t]][\varphi]\otimes_{k[[t]]}N\quad\longrightarrow\quad 1\otimes
N\quad\bigoplus\quad C$$where the first arrow is the inclusion, the last arrow
is the projection. This map is bijective, the critical point being the
computation$$t (k.t^{q^{i-1}k_{\nu}-1}\varphi^i\otimes
\tau(\nu))=k.t^{q^{i-1}k_{\nu}}\varphi^i\otimes
\tau(\nu)=k.\varphi^{i-1}t^{k_{\nu}}\varphi\otimes \tau(\nu)\equiv
k.\varphi^{i-1}\otimes \nu$$ modulo $\nabla$ (for $i>0$). It follows that indeed the image of $1\otimes N$ in $\Delta$ is the kernel of $t$ acting on $\Delta$.\hfill$\Box$\\

{\bf Definition:} An object $\Delta\in {\rm Mod}^{\rm ad}({\mathfrak O})$ is called standard
cyclic if it is generated over $k[[t]][\varphi]$ by ${\rm ker}(t|_{\Delta})=\Delta[t]$ and if there is a basis of
$\Delta[t]$ consisting of $\Gamma$-eigenvectors $e_0,\ldots, e_d$ such that $$t^{k_i}\varphi e_{i-1}=\rho_ie_i\quad\quad\mbox{ for all }\quad 0\le i\le d$$(reading
$e_{-1}=e_d$), for certain $0\le k_i\le q-1$ and $\rho_i\in k^{\times}$ such
that $k_i>0$ for at least one $i$, as well as $k_i<q-1$
for at least one $i$.\\

In the following, we extend any indexing by $0,\ldots,d$ to an indexing by
${\mathbb Z}$ in the obvious way (i.e. $k_i=k_{i+d+1}$,
$e_i=e_{i+d+1}$, $\rho_i=\rho_{i+d+1}$, $\eta_i=\eta_{i+d+1}$ for all
$i\in{\mathbb Z}$). Let $\nabla$ denote the $k[[t]][\varphi]$-submodule of $k[[t]][\varphi]\otimes_{k[[t]]}\Delta[t]$ generated by the elements $t^{k_i}\varphi\otimes e_{i-1}-1\otimes \rho_ie_i$. The inclusion $\Delta[t]\to \Delta$ extends to a natural $k[[t]][\varphi]$-linear map \begin{gather}k[[t]][\varphi]\otimes_{k[[t]]}\Delta[t]/\nabla\longrightarrow\Delta.\label{bohrlaerm}\end{gather} 

\begin{pro}\label{altfund} Let $\Delta\in {\rm Mod}^{\rm ad}({\mathfrak O})$ be standard
cyclic, with $e_i$, $k_i$, $\rho_i$, $\rho_i$ as above.

(a) $t$ acts surjectively on $\Delta$, and there is a distinguished isomorphism of free $k[[t]]$-modules of rank
  $d+1$\begin{gather}\Delta^*\cong k[[t]]\otimes_k(\Delta[t]^*).\label{saabvo50}\end{gather}The map (\ref{bohrlaerm}) is a $k[[t]][\varphi]$-linear isomorphism.

(b) If for any $1\le j\le d$ there is some $0\le i\le d$ with $k_i\ne k_{i+j}$, then $\Delta$ is irreducible as a $k[[t]][\varphi]$-module.

(c) For $0\le i\le d$ let $\eta_i:\Gamma\to k^{\times}$ be the character with
  $\gamma\cdot e_i=\eta_i(\gamma)e_i$ for all $\gamma\in\Gamma$. Suppose that
  for any $1\le j\le d$ which satisfies $k_i=k_{i+j}$ for all $0\le i\le d$
  there is some $0\le i\le d$ with $\eta_{i}\ne\eta_{i+j}$. Then $\Delta$ is
  irreducible as an ${\mathfrak O}$-module.

(d) At least after a finite extension of $k$ we have: $\Delta$ admits a filtration such that each associated
  graded piece is an
  irreducible standard cyclic object in ${\rm Mod}^{\rm ad}({\mathfrak O})$. If $p$ does not divide $d+1$ then $\Delta$ is even the direct sum of
  irreducible standard cyclic objects in ${\rm Mod}^{\rm ad}({\mathfrak O})$.
\end{pro}

{\sc Proof:} (This is very similar to \cite{dfun} Proposition 6.2.) 

(a) For $0\le j\le d$ consider $$w_j=k_j+qk_{j-1}+\ldots+q^jk_0+q^{j+1}k_d+\ldots+q^dk_{j+1}.$$Repeated
substitution of $\varphi t=t^q\varphi$ (recall $\Phi(t)=t^q$ modulo $\pi$)
shows that $t^{w_j}\varphi^{d+1}e_j\in k^{\times}e_j$. As $k_i>0$ for at least
one $i$ we have $w_j>0$, and hence $e_j\in t\Delta$. As $\Delta[t]$ is
generated over $k$ by all $e_j$ it follows that $\Delta[t]\subset t\Delta$. As
$\Delta$ is generated over $k[[t]][\varphi]$ by $\Delta[t]$, the equation
$\varphi t=t^q \varphi$ therefore shows $\Delta\subset t\Delta$, i.e. $t$ acts
surjectively on $\Delta$. We deduce that $\Delta^*$ is a torsion free, and hence free $k[[t]]$-module of rank $d+1$. As $\Delta$ is generated over $k[[t]][\varphi]$ by $\Delta[t]$ the map (\ref{bohrlaerm}) is surjective. But it is also injective, because Lemma \ref{abstrtor} tells us that it induces an isomorphisms between the respective kernels of $t$. We view the bijective map (\ref{bohrlaerm}) as an identification. The proof of Lemma \ref{abstrtor} yielded a canonical $k$-vector space
decomposition $\Delta=C\oplus \Delta[t]$ where the $k$-subvector
space $C$ of $\Delta$ is generated by the image elements of the elements
$t^{\theta}\varphi^r\otimes e\in k[[t]][\varphi]\otimes_{k[[t]]}\Delta[t]$
which do not belong to $1\otimes  \Delta[t]$ (for some $e\in\Delta[t]$, and
some $\theta, r\ge0$). We may thus identify $\Delta[t]^*={\rm Hom}_k(\Delta[t],k)$ with the subspace of
$\Delta^*={\rm Hom}_k(\Delta,k)$ consisting of all
$f\in \Delta^*$ with $f|_C=0$. The composition of this $k$-linear embedding
$\Delta[t]^*\to\Delta^*$ with the projection $\Delta^*\to
(\Delta^*)/t(\Delta^*)$ is a $k$-linear isomorphism. Therefore, and as
$\Delta^*$ is free and finitely generated over $k[[t]]$, the $k[[t]]$-linear
map $k[[t]]\otimes_k(\Delta[t]^*)\to \Delta^*$ extending the $k$-linear
embedding $\Delta[t]^*\to\Delta^*$ is an isomorphism as stated in formula (\ref{saabvo50}). 

(b) Let $Z$ be a non zero $k[[t]][\varphi]$-submodule of $\Delta$. With
$\Delta$ also $Z$ is a torsion $k[[t]]$-module, hence ${\rm ker}(t|_{Z})=Z[t]$
is non zero. For non zero elements $z=\sum_{0\le i\le d}x_ie_i$ of $Z[t]$
(with $x_i\in k$) put \begin{align}{\mathcal D}(z)&=\{0\le i\le
  d\,|\,x_i\ne0\}, &\nu(z)&=|{\mathcal D}(z)|,\notag\\\eta(z)&={\rm
    max}\{k_i\,|\,i\in {\mathcal D}(z)\},&\Lambda(z)&=t^{\eta(z)}\varphi z.\notag\end{align}Then
  $\Lambda(z)$ is again a non zero element of $Z[t]$. We have $${\mathcal
    D}(\Lambda(z))=\{i+1\,|\,\eta(z)=k_i\mbox{ and }i\in {\mathcal
    D}(z)\}$$(we read elements in $\{0\le i\le
  d\}$ modulo $(d+1)$), in particular $\nu(\Lambda(z))\le \nu(z)$. If $\nu(\Lambda(z))= \nu(z)$
  then ${\mathcal
    D}(\Lambda(z))=\{i+1\,|\,i\in {\mathcal D}(z)\}$ and $k_{i}=k_{i+j}$ whenever $i, i+j\in{\mathcal
    D}(z)$. This implies that if we had $\nu(\Lambda^n(z))= \nu(z)>1$ for all $n\ge0$
  then there was some $1\le j\le d$ with $k_i= k_{i+j}$ for all $0\le i\le
  d$. But this
  would contradict our hypothesis. Thus, for sufficiently large $n\ge0$ we
  have $\nu(\Lambda^n(z))=1$, i.e. $\Lambda^n(z)\in k^{\times}e_i$ for some $0\le i\le d$. For such $n$ we then even have $\Lambda^{n+j}(z)\in k^{\times}e_{i+j}$ for all $j\ge0$. It follows that $Z$ contains all $e_i$, hence $Z=\Delta$.

(c) We use the functions $\nu$, $\Lambda$ already employed in the proof of
  (b). Let $0\ne
Z\subset \Delta$ be a nonzero
${\mathfrak O}$-submodule. Choose a non zero $z\in Z[t]$ for which $\nu(z)$ is minimal (for all non
zero $z\in Z[t]$). If $\nu(z)=1$ then we obtain $Z=\Delta$ as in the proof of
(b). Now assume $\nu(z)>1$. For all $n\ge0$ we have $\nu(\Lambda^n(z))\le
\nu(z)$, hence $\nu(\Lambda^n(z))= \nu(z)$ by the choice of $z$. Thus, writing $z=\sum_{0\le i\le
  d}x_{i}e_i$ with $x_i\in k$, we have $x_{i}\ne0$ and $x_{i+j}\ne0$ for some $i,j$, with $j$ violating the hypothesis in (b). By the hypothesis in (c), replacing $i$ by $i+n$ and $z$ by $\Lambda^n(z)$ we may assume that
$\eta_i\ne\eta_{i+j}$. Pick $\gamma\in\Gamma$ with
$\eta_i(\gamma)\ne\eta_{i+j}(\gamma)$, and pick $a\in k^{\times}$ with $a
e_i=\gamma\cdot e_i$. Then $az-\gamma\cdot z$ is a non zero element in $Z[t]$
with $\nu(az-\gamma\cdot z)<\nu(z)$: a contradiction.

(d) Passing to a
finite extension of $k$ if necessary we may assume that there is a $(d+1)$-st
root of $\prod_{i=0}^d\rho_i$ in $k$. Thus, rescaling the $e_i$ if necessary we may
assume $\rho_i=\rho_j$ for all $i, j$. We argue by induction on $d$. If $\Delta$ itself is not irreducible then there is, by (c), some $1\le j\le d$ which satisfies
$k_i=k_{i+j}$ and $\eta_{i}=\eta_{i+j}$ for all $0\le i\le d$. The minimal
such $j$ is a divisor of $d+1$. Consider the $k$-subvector space $V$ of $\Delta[t]$ spanned by the vectors $\epsilon_i=e_{ij}$ for $0\le i< \frac{d+1}{j}$. Then $$(\prod_{i=1}^j\rho_i^{-1})t^{k_j}\varphi\cdots t^{k_1}\varphi$$induces the automorphism $f$ of $V$ with $f(\epsilon_i)=\epsilon_{i+1}$ (where we understand $\epsilon_{\frac{d+1}{j}}=\epsilon_0$). Choose (after passing to a finite extension of $k$ if necessary) an $f$-stable filtration $0=V_0\subset V_1\subset\ldots\subset V_{\frac{d+1}{j}}=V$ such that each $V_i/V_{i-1}$ is one dimensional. Then define for $0\le s\le\frac{d+1}{j}$ the ${\mathfrak O}$-submodule $\Delta_s={\mathfrak O} V_{0}+\cdots+{\mathfrak O} V_{s}$ of $\Delta$. It induces on $\Delta[t]$ the filtration$$\Delta_s[t]= \Delta_{s-1}[t]+V_{s}+t^{k_1}\varphi V_s+\ldots+ t^{k_{j-1}}\varphi\cdots t^{k_1}\varphi V_s.$$By construction, each $\Delta_{i+1}/\Delta_{i}$ is standard cyclic, and the induction hypothesis applies. If $p$ does not divide $\frac{d+1}{j}$ then there is even an $f$-stable direct sum decomposition $V=\oplus_sV_{[s]}$ with one dimensional $V_{[s]}$'s. Then $\Delta=\oplus_s\Delta_{[s]}$ with $\Delta_{[s]}={\mathfrak O}V_{[s]}$ is a direct sum decomposition of $\Delta$, and each $\Delta_{[s]}$ is standard cyclic, and the induction hypothesis applies.\hfill$\Box$\\

\begin{lem}\label{papa82} Let $\Delta\in {\rm Mod}^{\rm ad}({\mathfrak O})$ be standard
cyclic and put ${\bf D}=\Delta^*\otimes_{k[[t]]}k((t))\in {\rm Mod}^{et}(k((t)))$, see Proposition \ref{nopsi}. We have ${\bf D}^{\natural}=\Delta^*={\bf D}^{\sharp}$.
\end{lem}

{\sc Proof:} In Proposition \ref{altfund} we saw that $\Delta^*$ is a free $k[[t]]$-module, hence the natural map $\Delta^*\to{\bf
    D}=\Delta^*\otimes_{k[[t]]}k((t))$ is injective; we view it as an inclusion.

The $\varphi$-operator on $\Delta$ is the adjoint of the $\psi$-operator on ${\bf D}$,
in such a way that $\psi(\Delta^*)=\Delta^*$ since $\varphi$ acts injectively
on $\Delta$. Therefore the definitions of ${\bf D}^{\natural}$ and ${\bf
  D}^{\sharp}$ yield ${\bf D}^{\natural}\subset \Delta^*\subset{\bf
  D}^{\sharp}$. Since ${\bf
  D}^{\natural}$ is a lattice with $\psi({\bf
  D}^{\natural})={\bf
  D}^{\natural}$ we get $t\Delta^{\sharp}\subset{\bf
  D}^{\natural}$, together\begin{gather}{\bf D}^{\natural}\subset \Delta^*\subset{\bf
  D}^{\sharp}\quad\quad\mbox{ and }\quad\quad t{\bf D}^{\sharp}\subset{\bf
    D}^{\natural}.\label{jesusimTempel}\end{gather}

 Let $e_i$ and $k_i$ be as in the definition of $\Delta$ being standard
cyclic.

Formula (\ref{jesusimTempel}) implies $t\Delta^*\subset{\bf
  D}^{\natural}$, hence $t(\Delta^*/{\bf
  D}^{\natural})=0$, hence $\Delta^*/{\bf
  D}^{\natural}$ is dual to a subspace $W$ of $\Delta[t]$ stable under $\varphi$. To prove ${\bf D}^{\natural}=\Delta^*$ it is therefore enough to prove that $\Delta[t]$ does not contain a non-zero subspace $W$ stable under $\varphi$. Assume that such a $W$ does exist. A non-zero element $\beta\in W$ may be written as $\beta=\sum_{i=0}^{d}\alpha_i e_i$ with $\alpha_i\in k$. Let $k={\rm max}\{k_{i+1}\,|\,\alpha_i\ne0\}$. Since by assumption $k_i>0$ for at least one $i$, replacing $\beta$ by $\varphi^r\beta$ for some $r\in{\mathbb N}$ if necessary, we may assume $k>0$. But then $t^k\varphi\beta$ is a non-zero linear combination of the $e_i$, whereas we also have $t\varphi\beta=0$ since $\varphi\beta\in W\subset \Delta[t]$: a contradiction.

Formula (\ref{jesusimTempel}) implies $t{\bf
  D}^{\sharp}\subset\Delta^*$, i.e. $t({\bf
  D}^{\sharp}/\Delta^*)=0$. We endow ${\bf
  D}^{\sharp}$ and all its submodules with the $t$-adic topology. By
Pontrjagin duality (as recalled e.g. in \cite{schven}) we in particular have ${\rm
  Hom}_k^{\rm cont}(\Delta^*,k)=\Delta$. Now $t({\bf
  D}^{\sharp}/\Delta^*)=0$ means that the kernel $W$ of the natural
projection ${\rm
  Hom}_k^{\rm cont}({\bf
  D}^{\sharp},k)\to {\rm
  Hom}_k^{\rm cont}(\Delta^*,k)=\Delta$ is contained in ${\rm
  Hom}_k^{\rm cont}({\bf
  D}^{\sharp},k)[t]$. As $t$ acts injectively on ${\bf
  D}^{\sharp}$, it acts surjectively on ${\rm
  Hom}_k^{\rm cont}({\bf
  D}^{\sharp},k)$. Hence, if $\Delta^*\ne{\bf D}^{\sharp}$ then $W\ne0$ and
there is some $\beta\in {\rm
  Hom}_k^{\rm cont}({\bf
  D}^{\sharp},k)$ with $0\ne t\beta\in W$. Now $t\beta\in W$ means that $\beta$ maps to an
element in $\Delta[t]$. Since on the other hand $tW=0$ (as $W\subset {\rm
  Hom}_k^{\rm cont}({\bf
  D}^{\sharp},k)[t]$) we may write $\beta=\sum_{i=0}^{d}\alpha_i\widetilde{e}_i$ with
$\alpha_i\in k$, where $\widetilde{e}_i\in {\rm
  Hom}_k^{\rm cont}({\bf
  D}^{\sharp},k)$ lifts $e_i$. We
then also have $0\ne t\widetilde{e}_{i_0}\in W$ for some $i_0$. As $\varphi$
is injective on $W$ (which follows
from the surjectivity of $\psi$ on ${\bf D}^{\sharp}$ and hence on $W^*={\bf
  D}^{\sharp}/\Delta^*$) this gives $t^q\varphi\widetilde{e}_{i_0}=\varphi t\widetilde{e}_{i_0}\ne0$ in ${\rm
  Hom}_k^{\rm cont}({\bf
  D}^{\sharp},k)$. Together with $W\subset {\rm
  Hom}_k^{\rm cont}({\bf
  D}^{\sharp},k)[t]$ we get $t^{q-1}\varphi {e}_{i_0}\ne0$ in $\Delta$. Applying the same argument with
$t^{q-1}\varphi\widetilde{e}_{i_0}$ instead of $\widetilde{e}_i$ (again using that
$t^q\varphi\widetilde{e}_{i_0}\ne0$) we see $t^{q-1}\varphi t^{q-1}\varphi
{e}_{i_0}\ne0$. Next we get $t^{q-1}\varphi t^{q-1}\varphi t^{q-1}\varphi
e_{i_0}\ne0$ etc.. But this means $q-1=k_i$ for each $i$, contradicting the
hypothesis. We obtain $\Delta^*={\bf D}^{\sharp}$.\hfill$\Box$\\

{\bf Definition:} Let ${\rm Mod}^{\clubsuit}({\mathfrak O})$ denote
the subcategory of ${\rm Mod}^{\rm ad}({\mathfrak O})$ whose objects admit a filtration such that each associated graded piece becomes a standard
cyclic object in ${\rm Mod}^{\rm ad}({\mathfrak O})$ after a suitable field extension of $k$.\\

{\bf Remark:} Proposition \ref{altfund} (d) implies that each subquotient in ${\rm Mod}^{\rm ad}({\mathfrak O})$ of an object in ${\rm Mod}^{\clubsuit}({\mathfrak O})$ again is an object in ${\rm Mod}^{\clubsuit}({\mathfrak O})$. 

\begin{pro}\label{wienaus} The restriction of the functor (\ref{wienschluss}) to the category
  ${\rm Mod}^{\clubsuit}({\mathfrak O})$ is exact and fully faithful.
\end{pro}

{\sc Proof:} We already know that the functor is exact. Next, we claim \begin{gather}{\bf
  D}^{\natural}=\Delta^*={\bf D}^{\sharp}\quad\mbox{ with }\quad {\bf
  D}=\Delta^*\otimes_{k[[t]]}k((t))\label{vierwovorend}\end{gather}for $\Delta\in {\rm
  Mod}^{\clubsuit}({\mathfrak O})$. Indeed, for standard
cyclic $\Delta$ this is shown in Lemma \ref{papa82}. For $\Delta$ which become standard cyclic after a field extension $k'/k$ it then follows since the definitions of $(.)^{\natural}$ and $(.)^{\sharp}$ in terms of the $k$-linear operator $\psi$ imply ${\bf D}^{\natural}\otimes_kk'=({\bf D}\otimes_kk')^{\natural}$ and ${\bf D}^{\sharp}\otimes_kk'=({\bf D}\otimes_kk')^{\sharp}$. For general $\Delta\in {\rm
  Mod}^{\clubsuit}({\mathfrak O})$ it then follows from Lemma
\ref{goldhochzeit}. We now claim that the reverse functor (on the essential image of the
functor under discussion) is given by sending ${\bf D}$ to the topological dual $({\bf
  D}^{\natural})'$ of ${\bf
  D}^{\natural}$ (where we endow ${\bf D}^{\natural}$ with its $t$-adic topology). Indeed, for ${\bf D}$ in this essential image and for $\Delta\in {\rm
  Mod}^{\clubsuit}({\mathfrak O})$ we have natural isomorphisms$$(({\bf
  D}^{\natural})')^*\otimes_{k[[t]]}k((t))\stackrel{(i)}{\cong} {\bf
  D}^{\natural}\otimes_{k[[t]]}k((t))\cong {\bf
  D},$$$$((\Delta^*\otimes_{k[[t]]}k((t)))^{\natural})'\stackrel{(ii)}{\cong}
(\Delta^*)'\stackrel{(iii)}{\cong} \Delta,$$where $(i)$ and $(iii)$ follow
from Pontrjagin duality, see e.g. Proposition 5.4 in \cite{schven}, and where
$(ii)$ follows from formula (\ref{vierwovorend}).\hfill$\Box$\\

\subsection{Standard cyclic \'{e}tale $(\varphi,\Gamma)$-modules}

\label{standardcycet}

\begin{pro}\label{fudonsup} Let $\Delta\in{\rm Mod}^{\rm ad}({\mathfrak O})$ be a standard
cyclic object, with $d$, $e_i$, $k_i$, $\rho_i$, $\eta_i$ as in the
definition resp. as in Proposition \ref{altfund}. The \'{e}tale $(\varphi,\Gamma)$-module $\Delta^*\otimes_{k[[t]]}k((t))$ over
$k((t))$ admits a $k((t))$-basis $f_0,\ldots, f_d$ such that for all $0\le j\le d$ we have\begin{gather}\varphi(f_{j-1})=\rho_{j-1}^{-1}t^{1+k_j-q}f_j\label{tobfrag}\end{gather} (reading
  $f_{-1}=f_d$), and moreover \begin{gather}\gamma\cdot
    f_j-\eta_j^{-1}(\gamma) f_j\in tk[[t]]f_j\quad\quad\mbox{ for all }\gamma\in\Gamma.\label{monueber}\end{gather}
\end{pro}

{\sc Proof:} We use formula (\ref{psieqqp}).

First we assume $F\ne{\mathbb Q}_p$. Put $N=\oplus_{i=0}^dk.e_i$. As explained in the proof of Proposition \ref{altfund}, we have a bijective map (\ref{bohrlaerm}) which we view as an identification. In particular, Lemma \ref{abstrtor} and its proof apply. In the context of that proof we identify $e_i$ with the class of $1\otimes e_i$ in
$\Delta$. By formula (\ref{bohrnerv1}) we have a $k$-linear isomorphism
$(1\otimes N)\oplus C
\cong \Delta$ with $C$ as in formula (\ref{bohrnerv}). For $0\le j\le d$ we
may therefore define $f_j\in \Delta^*$ by asking $f_j(C)=0$ and $f_j(e_i)=\delta_{ij}$ for
$0\le i\le d$. Proposition \ref{altfund} tells us that $f_0,\ldots, f_d$ is a
$k[[t]]$-basis of $\Delta^*$. For $\theta,r\ge0$ and any $i,j$ we have $f_j(t^{\theta}\varphi^r\otimes e_i)\ne0$ if and only if $r\equiv j-i$ modulo $(d+1){\mathbb Z}$ and $\theta=k_j+qk_{j-1}+\ldots+q^{r-1}k_{j-r+1}$. As before, $\psi\in{\rm End}_k(\Delta^*)$ is defined by $(\psi(f))(x)=f(\varphi(x))$ for $x\in\Delta$, $f\in\Delta^*$. We claim\begin{gather}\psi(t^{m+k_j+1}f_j)=\rho_{j-1}\psi_{k((t))}(t^{m})tf_{j-1}\label{neurech}\end{gather}for
all $j$, all $m\ge-k_j-1$. Indeed, for
$0\le i\le d$ and $\theta, r\ge0$ we have$$(\psi(t^{m+k_j+1}f_j))(t^{\theta}\varphi^r\otimes
e_i)=f_j(t^{m+k_j+1}\varphi t^{\theta}\varphi^r\otimes e_i).$$If
$m+1\notin{\mathbb Z}q$ then this shows $(\psi(t^{m+k_j+1}f_j))(t^{\theta}\varphi^r\otimes
e_i)=0$ by what we pointed out above. But $m+1\notin{\mathbb Z}q$ also implies
$\psi_{k((t))}(t^{m})=0$. In the
case $m+1=qn$ (some $n\in{\mathbb Z}$) we compute
$$(\psi(t^{m+k_j+1}f_j))(t^{\theta}\varphi^r\otimes
e_i)=f_j(t^{k_j+qn}\varphi t^{\theta}\varphi^r\otimes
e_i)=f_j(t^{k_j}\varphi t^{n+\theta}\varphi^r\otimes
e_i)$$$$=\rho_{j-1}f_{j-1}(t^{n+\theta}\varphi^r\otimes
e_i)=(\rho_{j-1}\psi_{k((t))}(t^{m})tf_{j-1})(t^{\theta}\varphi^r\otimes
e_i)$$where we used 
$\psi_{k((t))}(t^{m})=t^{n-1}$. We have proven formula
(\ref{neurech}). 

On the other hand, by tracing the construction in Proposition \ref{nopsi} we see
that $\varphi(tf_{j-1})$ is characterized by satisfying
\begin{gather}\psi(t^{m} \varphi(tf_{j-1}))=\psi_{k((t))}(t^{m})tf_{j-1}\label{sekrgeschwaetz}\end{gather} for all
$m$. Comparing formulae (\ref{neurech}) and (\ref{sekrgeschwaetz}) we find
$\varphi(tf_{j-1})=\rho_{j-1}^{-1}t^{k_j+1}f_j$ which is equivalent with formula (\ref{tobfrag}). Next, for $\gamma\in\Gamma$ we compute$$(\gamma\cdot
f_j)(e_i)=f_j(\gamma^{-1}\cdot
e_i)=f_j(\eta_i(\gamma^{-1})e_i)=(\eta_i(\gamma^{-1})f_j)(e_i)=(\eta_j(\gamma^{-1})f_j)(e_i).$$Here the last equation is trivial if $i=j$, whereas if $i\ne j$ then both sides vanish. This shows $(\gamma\cdot
f_j-\eta_j(\gamma^{-1})f_j)|_N=0$, and hence $\gamma\cdot
f_j-\eta_j(\gamma^{-1})f_j\in t\Delta^*=tk[[t]]\{f_0,\ldots, f_d\}$. On the other hand, by what we pointed out above, $(\gamma\cdot f_j)(t^{\theta}\varphi^r\otimes e_i)=f_j([\gamma]_{\Phi}(t)^{\theta}\varphi^r\otimes e_i)$ vanishes if $r+i-j\notin (d+1){\mathbb Z}$, and this shows $\gamma\cdot f_j\in k[[t]]f_j$. We trivially have $\eta_j(\gamma^{-1})f_j\in  k[[t]]f_j$, and hence altogether $\gamma\cdot
f_j-\eta_j(\gamma^{-1})f_j\in tk[[t]]\{f_0,\ldots, f_d\}\cap k[[t]]f_j=tk[[t]]f_j$, formula (\ref{monueber}).

Now we assume $F={\mathbb Q}_p$. Let us suppose for simplicity that $\pi=q$. For $0\le j\le d$ we
may define $f_j\in \Delta^*$ as follows. For $\theta, r\ge0$ (and any $i, j$) we require $f_j(t^{\theta}\varphi^r\otimes e_i)\ne0$ if and only if $r\equiv j-i$ modulo $(d+1){\mathbb Z}$ and there are $a_1,\ldots,a_{r-1}\in\{0,1\}$ such that $$\theta=k_j+qk_{j-1}+\ldots+q^{r-1}k_{j-r+1}+\sum_{i=1}^{r-1}a_iq^{i-1}(1-q);$$if this is the case we put $$f_j(t^{\theta}\varphi^r\otimes e_i)=\rho_{j-1}\rho_{j-2}\cdots\rho_{j-r}.$$(As usual, the subindices of the $\rho_?$ are read modulo $(d+1){\mathbb Z}$.) Again $f_0,\ldots, f_d$ is a
$k[[t]]$-basis of $\Delta^*$. Again we claim formula (\ref{neurech}). As before we see that both sides vanish if $m\notin{\mathbb Z}q-1\cup {\mathbb Z}q$, and coincide if $m\in{\mathbb Z}q-1$. But the same computation also shows their coincidence if $m=qn$ for some $n\in{\mathbb N}$, as follows:$$(\psi(t^{m+k_j+1}f_j))(t^{\theta}\varphi^r\otimes
e_i)=f_j(t^{k_j+1}\varphi t^{n+\theta}\varphi^r\otimes
e_i)$$$$=\rho_{j-1}f_{j-1}(t^{n+\theta+1}\varphi^r\otimes
e_i)=(\rho_{j-1}\psi_{k((t))}(t^{m})tf_{j-1})(t^{\theta}\varphi^r\otimes
e_i)$$where we used 
$\psi_{k((t))}(t^{m})=t^{n}$. With formula (\ref{neurech}) being established, the remaining arguments are exactly as before.\hfill$\Box$\\

{\bf Definition:} We say that an object ${\bf D}\in{\rm Mod}^{et}(k((t)))$
of dimension $d+1$ is standard cyclic if it admits a $k((t))$-basis $f_0,\ldots, f_d$ such that
there are $\sigma_j\in k^{\times}$, characters $\alpha_j:\Gamma\to k^{\times}$ and $m_j\in \{1-q,\ldots,-1,0\}$ for $0\le j\le d$ satisfying the following conditions:$$(m_0,\ldots, m_d)\notin\{(0,\ldots,0), (1-q,\ldots,1-q)\},$$$$\varphi(f_{j-1})=\sigma_j t^{m_j}f_j\quad\quad\mbox{ for all }j\mbox{ (reading }f_{-1}=f_d),$$$$\gamma\cdot f_j-\alpha_j(\gamma) f_j\in tk[[t]]\{f_0,\ldots,f_d\}\quad\quad\mbox{ for all }\gamma\in\Gamma.$$

\begin{lem}\label{ossa} (a) The constant $\prod_{j=0}^d\sigma_j\in k^{\times}$
  as well as, up to cyclic permutation, the ordered tuple
  $((\alpha_0,m_0),\ldots,(\alpha_d, m_d))$, are uniquely determined by the
  isomorphism class of the $(\varphi,\Gamma)$-module ${\bf D}$.

(b) $\alpha_1,\ldots,\alpha_d$ are uniquely determined by $\alpha_0$ and $m_0,\ldots,m_d$.

\end{lem}

{\sc Proof:} (a) In the following, for elements of ${\rm GL}_{d+1}(k((t)))$ we read the (two) respective indices of their entries always modulo $(d+1){\mathbb Z}$. 

The effect of $\varphi$ on the basis $f_0,\ldots, f_d$ is described by $T=(T_{ij})_{0\le i,j\le d}\in{\rm GL}_{d+1}(k((t)))$ with $T_{i, i+1}=\sigma_it^{m_i}$ for $0\le i\le d$, but $T_{i, j}=0$ for $j\ne i+1$.

Let $\sigma'_j\in k^{\times}$ and $((\alpha'_0,m'_0),\ldots,(\alpha'_d, m'_d))$ be another datum as above, let ${\bf D}'$ be an \'{e}tale $(\varphi,\Gamma)$-module admitting a $k((t))$-basis $f'_0,\ldots, f'_d$ with $\varphi(f'_{j-1})=\sigma'_j t^{m'_j}f'_j$ and $\gamma\cdot f'_j-\alpha'_j(\gamma) f'_j\in tk[[t]]\{f'_0,\ldots,f'_d\}$ for $\gamma\in\Gamma$. Define $T'=(T'_{ij})_{0\le i,j\le d}\in{\rm GL}_{d+1}(k((t)))$ similarly as above.

Suppose that there is an isomorphism of $(\varphi,\Gamma)$-modules ${\bf
  D}'\cong {\bf D}$. With respect to the bases $f_{\bullet}$ and
$f'_{\bullet}$ it is described by some $A(t)=(a_{i,j}(t))_{0\le i,j\le d}\in
{\rm GL}_{d+1}(k((t)))$. In view of $\varphi t= \Phi(t) \varphi$,
the compatibility of the isomorphism with the respective $\varphi$-actions
comes down to the matrix equation $$T\cdot A(t)=A(\Phi(t))\cdot T'.$$For the individual entries this is equivalent
with$$a_{i,j}(t)=\sigma'_j\sigma_i^{-1}t^{m'_j-m_i}a_{i-1,j-1}(\Phi(t))$$for
  all $i,j$. Iteration of this equation yields$$a_{i,j}(t)=(\prod_{\ell=0}^d
  \sigma'_{j-\ell}\sigma_{i-\ell}^{-1}(\Phi^{\ell}(t))^{m'_{j-\ell}-m_{i-\ell}})a_{i,j}(\Phi^{d+1}(t))$$for
  all $i,j$. (Here $\Phi^{\ell}(t)$ resp. $\Phi^{d+1}(t)$ means
  $\Phi(\Phi(\ldots \Phi(t)\ldots))$.) From this we deduce that for fixed $i, j$ either
  $a_{i,j}$ is a non zero constant and $\prod_{\ell=0}^d
  \sigma'_{j-\ell}\sigma_{i-\ell}^{-1}=1$ and $m'_{j-\ell}=m_{i-\ell}$ for all
  $\ell$, or $a_{i,j}=0$. But since $A(t)$ is invertible we do find $i,j$ with
  $a_{i,j}\ne 0$. It already follows that
  $\prod_{j=0}^d\sigma_j=\prod_{j=0}^d\sigma'_j$ and that $(m'_0,\ldots,m'_d)$
  coincides with $(m_0,\ldots, m_d)$ up to cyclic permutation. But since in
  addition we just saw that $A$ is a constant matrix, with $a_{i,j}=0$ if and
  only if $a_{i-1,j-1}=0$, we see that the same index permutation takes
  $\alpha'_j$ to $\alpha_j$.

(b) This follows from the fact that, in view of the defining formulae, ${\bf D}$ is generated by
  $f_0$ as a $\varphi$-module over $k((t))$.\hfill$\Box$\\

\begin{pro}\label{responseday} The functor $\Delta\mapsto
  \Delta^*\otimes_{k[[t]]}k((t))$ induces a bijection between the set of standard cyclic objects in ${\rm Mod}^{\rm ad}({\mathfrak O})$ and the set of standard cyclic objects in ${\rm Mod}^{et}(k((t)))$.
\end{pro}

{\sc Proof:} $\Delta^*\otimes_{k[[t]]}k((t))$ for a standard cyclic object $\Delta\in{\rm Mod}^{\rm ad}({\mathfrak O})$ is a standard cyclic object in ${\rm Mod}^{et}(k((t)))$ by Proposition \ref{fudonsup}. With Lemma \ref{ossa} (a) we see that the assignment
$\Delta\mapsto\Delta^*\otimes_{k[[t]]}k((t))$ is injective on standard cyclic objects in ${\rm Mod}^{\rm ad}({\mathfrak O})$. It is also surjective: Proposition \ref{fudonsup} (together with Lemma \ref{ossa} (b)) explicitly says how to convert the parameter data describing a standard cyclic object in ${\rm Mod}^{et}(k((t)))$ into the parameter data describing a standard cyclic object in ${\rm Mod}^{\rm ad}({\mathfrak O})$. \hfill$\Box$\\ 



{\bf Definition:} A $(d+1)$-dimensional standard cyclic ${\rm Gal}(\overline{F}/F)$-representation is a ${\rm Gal}(\overline{F}/F)$-representation over $k$ which corresponds, under the equivalence of categories in Theorem \ref{sosego}, to an object in ${\rm Mod}^{et}(k((t)))$ of dimension $d+1$ which is standard cyclic.

\section{Hecke algebras and supersingular modules}

\subsection{The pro-$p$-Iwahori Hecke algebra ${\mathcal H}$}

\label{defhsu}

We introduce the pro-$p$ Iwahori Hecke algebra ${\mathcal H}$ of ${\rm GL}_{d+1}(F)$ with coefficients in $k$ in a
slightly unorthodox way, which however is well suited for our later constructions. 

Let $\overline{T}$ be a free ${\mathbb Z}/(q-1)$-module of rank $d+1$. Then
 ${\rm Hom}(\Gamma,\overline{T})$ (with $\Gamma={\mathcal O}_F^{\times}$)
is also free of rank $d+1$ over ${\mathbb Z}/(q-1)$. We write the group law of
$\overline{T}$ multiplicatively, but that of ${\rm
  Hom}(\Gamma,\overline{T})$ we write additively. Let
$e^*,\alpha_1^{\vee},\ldots,\alpha_d^{\vee}$ be a ${\mathbb Z}/(q-1)$-basis of
${\rm Hom}(\Gamma,\overline{T})$. Put
$\alpha_{0}^{\vee}=-\sum_{i=1}^d\alpha_{i}^{\vee}$. We let the symmetric group
${\mathfrak S}_{d+1}$ act on ${\rm Hom}(\Gamma,\overline{T})$ as follows. We
think of ${\mathfrak S}_{d+1}$ as the permutation group of $\{0,1,\ldots,d\}$,
generated by the transposition $s=(01)\in {\mathfrak S}_{d+1}$ and the cycle
$\omega\in {\mathfrak S}_{d+1}$ with $\omega(i)=i+1$ for all $0\le i\le
d-1$. We then put$${\omega}\cdot e^* = e^*+\alpha_{0}^{\vee},\quad\quad{\omega}\cdot
\alpha_{0}^{\vee}=\alpha_{d}^{\vee}\quad\quad\mbox{ and }\quad\quad{\omega}
\cdot \alpha_{i}^{\vee}=\alpha_{i-1}^{\vee}\quad\mbox{ for }1\le i\le d.$$If $d=1$ we put$$s\cdot
e^* = e^*-\alpha_1^{\vee},\quad s\cdot
\alpha_i^{\vee}=-\alpha_i^{\vee}\quad\mbox{ for }i=0,1,$$but if $d\ge2$ we put$$s\cdot
e^* = e^*-\alpha_1^{\vee},\quad s\cdot
\alpha_0^{\vee}=\alpha_0^{\vee}+\alpha_1^{\vee},\quad s\cdot
\alpha_1^{\vee}=-\alpha_1^{\vee},\quad s\cdot
\alpha_2^{\vee}=\alpha_1^{\vee}+\alpha_2^{\vee},$$$$s\cdot
\alpha_i^{\vee}=\alpha_i^{\vee}\quad\mbox{ for }3\le i\le
d.$$ One easily checks that there is a unique action of ${\mathfrak S}_{d+1}$ on $\overline{T}$ such that for $\gamma\in\Gamma$ and $f\in {\rm
  Hom}(\Gamma,\overline{T})$ we have $$\omega\cdot (f(\gamma))=(\omega\cdot f)(\gamma)\quad\quad\mbox{ and }\quad\quad s\cdot (f(\gamma))=(s\cdot f)(\gamma).$$

Define $\alpha_1^{\vee}({\mathbb
    F}_q^{\times})$ to be the image of the composition ${\mathbb
    F}_q^{\times}\to\Gamma\stackrel{\alpha_1^{\vee}}{\to} \overline{T}$ where the first map is the Teichm\"uller homomorphism.\\

{\bf Definition:} (a) The $k$-algebra ${\mathcal H}$ is generated by elements
$T_{\omega}^{\pm 1}$, $T_s$ and $T_t$ for $t\in \overline{T}$, subject
to the following relations (with $t,t'\in\overline{T}$):
\begin{gather}T_sT_{\omega}T_sT_{\omega}^{-1}T_sT_{\omega}=T_{\omega}T_sT_{\omega}^{-1}T_sT_{\omega}T_s\quad\quad\mbox{
    if
  }d>1,\label{braid1}\\T_sT_{\omega}^{-m}T_sT_{\omega}^m=T_{\omega}^{-m}T_sT_{\omega}^{m}T_s\quad\quad\mbox{
    for all }1<m<d,\label{braid2}\\T_s^2=T_s\tau_s=\tau_sT_s\quad\quad\mbox{
  with }\quad\tau_s=\sum_{t\in\alpha_1^{\vee}({\mathbb
    F}_q^{\times})}T_{t},\label{sorge}\\T_{\omega}T_{\omega}^{-1}=1=T_{\omega}^{-1}T_{\omega},\label{hsv16},\\T_{\omega}^{d+1}T_s=T_sT_{\omega}^{d+1},\label{sanafabe}\\ T_{t}T_{t'}=T_{t'
t},\quad\quad T_{1_{\overline{T}}}=1,\label{hsv1}\\T_{t}T_{\omega}=T_{\omega}T_{{\omega}\cdot t},\label{hsvtabfue}\\T_{t}T_s=T_sT_{s\cdot t}.\label{hsvtabfue1}\end{gather}Notice that $T_{\omega}^{d+1}$ is central in ${\mathcal H}$. 

(b) ${\mathcal H}_{\rm aff}$ is the $k$-subalgebra of ${\mathcal
  H}$ generated by all $T_t$ for $t\in\overline{T}$, by
$T_{\omega}^{d+1}$, $T_{\omega}^{-d-1}$ and by
all $T_{\omega}^{m}T_sT_{\omega}^{-m}$ for $m\in{\mathbb Z}$.

(c) ${\mathcal H}^{\flat}$ is the quotient of ${\mathcal H}$ by
the two sided ideal spanned by all elements $T_t-1$ with $t\in\overline{T}$.\\

{\bf Caution:} ${\mathcal H}_{\rm aff}$ differs from the
similarly denoted algebra in \cite{vigneras}. (The difference is that here we
include $(T_{\omega}^{d+1})^{\mathbb Z}$.)\\

{\bf Remark:} Let $\overline{T}$ denote the subgroup of $G={\rm GL}_{d+1}(F)$ consisting of diagonal matrices with entries in the image of the Teichm\"uller homomorphism ${\mathbb F}_q^{\times}\to{\mathcal O}_F^{\times}$. For $\gamma\in\Gamma$ let $\overline{\gamma}$ be its image in ${\mathbb F}_q^{\times}$. In $\overline{T}$ define the elements $e^*(\gamma)={\rm diag}(\overline{\gamma},1_d)$ and $\alpha_i^{\vee}(\gamma)={\rm diag}(1_{i-1},\overline{\gamma},\overline{\gamma}^{-1},1_{d-i})$ for $1\le i\le d$. Define the elements $\omega=(\omega_{ij})_{0\le i,j\le d}$ and $s=(s_{ij})_{0\le i,j\le d}$ of $G$ by $\omega_{d0}=\pi$ and $\omega_{i,i+1}=1$ (for $0\le i\le d-1$) and $\omega_{ij}=0$ for all other pairs $(i,j)$, resp. by $s_{10}=s_{01}=s_{ii}=1$ for $i\ge 2$, and $s_{ij}=0$ for all other pairs $(i,j)$.

Let $I_0$ denote the pro-$p$-Iwahori subgroup of $G$ for which
$g=(g_{ij})_{0\le i,j\le d}\in G$ belongs to $I_0$ if and only if all the
following conditions are satisfied: $g_{ij}\in\pi{\mathcal O}_F$ for $i>j$,
and $g_{ij}\in{\mathcal O}_F$ for $i<j$, and $g_{ii}\in 1+\pi{\mathcal
  O}_F$. 

{\it Claim:} The corresponding pro-$p$-Iwahori Hecke algebra $k[I_0\backslash G/I_0]^{\rm op}\cong{\rm End}_{k[G]}({\rm ind}_{I_0}^G k)^{\rm op}$ is
isomorphic with ${\mathcal H}$, in such a way that the double coset
$I_0gI_0$ for $g\in\overline{T}\cup\{s,\omega\}$ corresponds to the element
$T_g\in {\mathcal H}$.

To prove this claim we use the description of $k[I_0\backslash G/I_0]^{\rm op}$ worked out by Vign\'{e}ras in \cite{vigneras} (or rather: we use the description of $k[I_0\backslash G/I_0]^{\rm op}$ which results from the description of $k[I_0\backslash G/I_0]$ given in \cite{vigneras}). 

Let $T$ denote the maximal torus of diagonal matrices in $G$, let $N(T)$ be its normalizer in $G$. Let $T_1$ (resp. $T_0$) denote the subgroup of $T$ consisting of diagonal matrices with entries in the kernel of ${\mathcal O}_F^{\times}\to {\mathbb F}_q^{\times}$ (resp. in ${\mathcal O}_F^{\times}$); thus $T_0/T_1\cong \overline{T}$. For $0\le i\le d$ define $s_i=\omega^{1-i}s\omega^{i-1}$. The (classes of) $s_0,s_1,\ldots, s_d$ are the Coxeter generators of a Coxeter subgroup $W_{\rm aff}$ of $N(T)/T_0$, and $N(T)/T_0$ is generated by $W_{\rm aff}$ together with the element $\omega$. The length function $\ell:W_{\rm aff}\to{\mathbb Z}_{\ge0}$ can be extended to a function $\ell:N(T)/T_0\to {\mathbb Z}_{\ge0}$ in such a way that $\ell(\omega)=0$. We again denote by $\ell$ the induced function $W^{(1)}=N(T)/T_1\to {\mathbb Z}_{\ge0}$. For $w\in W^{(1)}$ and $w'\in N(T)$ lifting $w$, the double coset $I_0w'I_0$ only depends on $w$; we denote it by $T_w$. For $0\le i\le d$ let $\overline{T}_i$ be the image of one of the two cocharacters ${\mathbb F}_q^{\times}\to \overline{T}$ associated with $s_i$. (Here we identify $\overline{T}$ with the maximal torus of diagonal matrices in ${\rm GL}_{d+1}({\mathbb F}_q)$. If $1\le i\le d$ then $s_i$ is the simple reflection associated with the coroot $\alpha_i^{\vee}$, and $\alpha_i^{\vee}({\mathbb F}_q^{\times})=\overline{T}_i$.) Now, according to \cite{vigneras}, a $k$-basis of $k[I_0\backslash G/I_0]^{\rm op}$ is given by the set of all $T_w$ for $w\in W^{(1)}$, and the multiplication is uniquely determined by the relations\begin{gather}T_wT_{w'}=T_{w'w}\quad\mbox{ for }w,w'\in W^{(1)}\mbox{ with }\ell(w)+\ell(w')=\ell(ww'),\label{gheorgetalk1}\\T^2_{s_i}=T_{s_i}\tau_{i}\quad\mbox{ where }\tau_{i}=\sum_{t\in\overline{T}_i}T_{t}\quad\mbox{ for }0\le i\le d.\label{gheorgetalk2}\end{gather}In the following we repeatedly use that conjugating these relations by powers of $T_{\omega}$ leads to similar relations (since $\ell(\omega)=0$). From formula (\ref{gheorgetalk1}) we first deduce $T_{s_i}=T_{\omega}^{i-1}T_sT_{\omega}^{1-i}$ and then that $T_{\omega}^{\pm 1}$ and $T_s=T_{s_1}$ together with the elements $T_t$ for $t\in\overline{T}$ generate $k[I_0\backslash G/I_0]^{\rm op}$ as a $k$-algebra. Next, from $s_is_{i-1}s_i=s_{i-1}s_is_{i-1}$ in $W_{\rm aff}$ (for $0\le i\le d$; if $i=0$ read $i-1=d$) we get $T_{s_i}T_{s_{i-1}}T_{s_i}=T_{s_{i-1}}T_{s_i}T_{s_{i-1}}$ by applying formula (\ref{gheorgetalk1}) twice, but this comes down to formula (\ref{braid1}) (up to conjugation by a power of $T_{\omega}$). Similarly from $s_is_j=s_{j}s_i$ in $W_{\rm aff}$ for $0\le i<j-1\le d-1$ with $i+d>j$ we get $T_{s_i}T_{s_j}=T_{s_j}T_{s_i}$ by applying formula (\ref{gheorgetalk1}) twice, but this comes down to formula (\ref{braid2}) (up to conjugation by a power of $T_{\omega}$). Formula (\ref{gheorgetalk2}) for any $i$ is a $T_{\omega}$-power conjugate of formula (\ref{sorge}). Finally, formulae (\ref{hsv16}), (\ref{sanafabe}), (\ref{hsv1}), (\ref{hsvtabfue}) and (\ref{hsvtabfue1}) are special instances of formula (\ref{gheorgetalk1}). Conversely, it is not hard to see that these, together with formulae (\ref{braid1}), (\ref{braid2}) and (\ref{sorge}) suffice to generate all relations in $k[I_0\backslash G/I_0]^{\rm op}$. The claim is proven.

We add: If $I$ denotes the Iwahori subgroup of $G$ containing $I_0$, then ${\mathcal H}^{\flat}$
becomes isomorphic with the Iwahori Hecke algebra $k[I\backslash G/I]^{\rm op}$.\\

{\bf Definition:} A character $\chi:{\mathcal H}_{\rm aff}\to k$ is
called {\it supersingular} if the following two conditions are both satisfied:

(a) There is an $m\in{\mathbb Z}$ with $\chi(T_{\omega}^m T_sT_{\omega}^{-m})=0$. 

(b) There is an $m\in{\mathbb Z}$ with either $\chi(T_{\omega}^m T_sT_{\omega}^{-m})=-1$ or
$\chi(T_{\omega}^m \tau_sT_{\omega}^{-m})=0$.\footnote{We have $\chi(T_{\omega}^m \tau_sT_{\omega}^{-m})=0$ if and only if $\chi(T_{\omega}^m T_tT_{\omega}^{-m})\ne 1$ for some $t\in\alpha_1^{\vee}({\mathbb
    F}_q^{\times})$, if and only if $\chi(\alpha_{m+1}^{\vee}(\gamma))\ne 1$ for some $\gamma\in\Gamma$.}\\

{\bf Definition:} (a) An ${\mathcal
  H}$-module $M$ is called {\it standard supersingular} if it is isomorphic with ${\mathcal
  H}\otimes_{{\mathcal H}_{\rm aff},\chi}k.e$, where ${\mathcal H}_{\rm aff}$ acts on the one dimensional $k$-vector space $k.e$ through a supersingular character $\chi$.

Equivalently, $M$ is standard supersingular if and only if $M=\bigoplus_{0\le m\le d}T_{\omega}^m(M_1)$ with an ${\mathcal H}_{\rm
  aff}$-module $M_1$ of $k$-dimension $1$ on which ${\mathcal H}_{\rm aff}$
acts through a supersingular character.\footnote{Then ${\mathcal H}_{\rm aff}$
acts on each $T_{\omega}^m(M_1)$ through a supersingular character.}

(b) An irreducible ${\mathcal
  H}$-module is called {\it supersingular} if it is a subquotient of a standard supersingular ${\mathcal
  H}$-module. 

A finite dimensional ${\mathcal H}$-module is called {\it supersingular} if each of its irreducible subquotients is supersingular. 

More generally, an ${\mathcal H}$-module is called supersingular if it is the inductive limit of its finite dimensional ${\mathcal H}$-submodules and if each finite dimensional ${\mathcal H}$-submodule is supersingular.\footnote{It is easy to see that the irreducible subquotients of a supersingular ${\mathcal H}$-module are the irreducible subquotients of its finite dimensional ${\mathcal H}$-submodules.}\\

{\bf Remark:} For non-zero finite dimensional ${\mathcal
  H}$-modules, the above definition of supersingularity is equivalent with
the one given by Vign\'{e}ras. This follows from the discussion in section 6 of \cite{vigjuss}. There is also a notion of supersingularity for ${\mathcal
  H}{{{}}}$-modules which are not necessarily inductive limits of their finite dimensional submodules. In the present paper however, without further mentioning {\it all ${\mathcal
  H}{{{}}}$-modules will be assumed to be inductive limits of their finite dimensional submodules.}\\

{\bf Remark:} In the literature on modules over Hecke algebras, the term {\it standard} module is occasionally used, but this usage is unrelated to our terminology.

\subsection{The coverings ${\mathcal H}^{\sharp\sharp}$ and ${\mathcal H}^{\sharp}$ of ${\mathcal H}$} 

\label{monachstras}

{\bf Definition:} (a) Let ${\mathcal H}^{\sharp}$ denote the $k$-algebra generated by
elements $T_{\omega}^{\pm 1}$, $T_s$ and $T_{t}$ for $t\in
\overline{T}$, subject to 

${\bullet}$ the relations
(\ref{sorge}), (\ref{hsv16}), (\ref{hsv1}), (\ref{hsvtabfue}),

${\bullet}$ the relations (\ref{hsvtabfue1}) for
$t=\alpha^{\vee}_{i}(\gamma)$ (all $0\le i\le d$, $\gamma\in\Gamma$),

${\bullet}$ the relation \begin{gather}
  T_{\omega}^{d+1}T_s^2=T_s^2T_{\omega}^{d+1},\label{sanafabehsv17}\end{gather}

${\bullet}$ the relations \begin{gather}T_{t}T^2_s=T^2_sT_{t}\quad\quad\mbox{ for all }t\in\overline{T},\label{hsvtabfue11}\end{gather}

${\bullet}$ the relations \begin{gather}T^2_sT_{\omega}T^2_sT_{\omega}^{-1}T^2_sT_{\omega}=T_{\omega}T^2_sT_{\omega}^{-1}T^2_sT_{\omega}T^2_s\quad\quad\mbox{
    if
  }d>1,\label{braid11}\\T^2_sT_{\omega}^{-m}T^2_sT_{\omega}^m=T_{\omega}^{-m}T^2_sT_{\omega}^{m}T^2_s\quad\quad\mbox{
    for all }1<m<d.\label{braid22}\end{gather}

(b) Let ${\mathcal H}^{\sharp\sharp}$ denote the $k$-algebra generated by
the elements $T_{\omega}^{\pm 1}$, $T_s$ and $T_{t}$ for $t\in
\overline{T}$, subject to 

${\bullet}$ the relations
(\ref{sorge}), (\ref{hsv16}), (\ref{hsv1}), (\ref{hsvtabfue}),

${\bullet}$ the relations (\ref{hsvtabfue1}) for
$t=\alpha^{\vee}_{i}(\gamma)$ (all $0\le i\le d$, $\gamma\in\Gamma$),

${\bullet}$ the relations (\ref{hsvtabfue11}).

\begin{lem}\label{2vorstubei} In ${\mathcal H}$ we have the relations (\ref{sanafabehsv17}), (\ref{hsvtabfue11}), (\ref{braid11}) and (\ref{braid22}).
\end{lem} 

{\sc Proof:} It is immediate that the relations (\ref{hsvtabfue1}) and (\ref{sanafabe}), imply the relations
(\ref{hsvtabfue11}) and (\ref{sanafabehsv17}), respectively. For $1<m<d$ and $t\in \alpha_1^{\vee}({\mathbb
  F}_q)$ we have $s\omega^m\cdot t=\omega^m\cdot t$, hence $T_s\sum_{t\in\alpha_1^{\vee}({\mathbb
  F}_q)}T_{\omega^m\cdot t}=\sum_{t\in\alpha_1^{\vee}({\mathbb
  F}_q)}T_{\omega^m\cdot t}T_s$. The same applies with $-m$ instead of $m$, hence$$T_sT_{\omega}^{-m}\tau_sT_{\omega}^{m}=T_{\omega}^{-m}\tau_sT_{\omega}^{m}T_s\quad\quad\mbox{
    and }\quad\quad
  T_sT_{\omega}^{m}\tau_sT_{\omega}^{-m}=T_{\omega}^{m}\tau_sT_{\omega}^{-m}T_s.$$This,
together with $T_s^2=\tau_sT_s=T_s\tau_s$ (formula (\ref{sorge})), justifies (i) and (iii)
in$$T_s^2T_{\omega}^{-m}T_s^2T_{\omega}^{m}\stackrel{(i)}{=}\tau_s(T_{\omega}^{-m}\tau_sT_{\omega}^{m})T_sT_{\omega}^{-m}T_sT_{\omega}^{m}\stackrel{(ii)}{=}\tau_s(T_{\omega}^{-m}\tau_sT_{\omega}^{m})T_{\omega}^{-m}T_sT_{\omega}^{m}T_s\stackrel{(iii)}{=}T_{\omega}^{-m}T_s^2T_{\omega}^{m}T_s^2,$$whereas
(ii) is justified by (\ref{braid2}). We have shown (\ref{braid22}). Finally, to see (\ref{braid11}) comes down, using (\ref{sorge}), (\ref{hsvtabfue}) and (\ref{hsvtabfue1}), to comparing$$T_{\omega}T_s^2T_{\omega}^{-1}T_s^2T_{\omega}T_s^2=(\sum_{t_1,t_2,t_3\in\alpha_1^{\vee}({\mathbb
  F}_q)}T_{\omega^{-1}\cdot
  t_1}T_{\omega^{-1} s\omega\cdot t_2}T_{\omega^{-1} s
  \omega s
  \omega^{-1}\cdot t_3})T_{\omega}T_sT_{\omega}^{-1}T_sT_{\omega}T_s,$$$$T_s^2T_{\omega}^{-1}T_s^2T_{\omega}T_s^2T_{\omega}=(\sum_{t_1,t_2,t_3\in\alpha_1^{\vee}({\mathbb
  F}_q)}T_{
  t_1}T_{s \omega^{-1}\cdot t_2}T_{s\omega^{-1} s\omega\cdot t_3})T_sT_{\omega}T_sT_{\omega}^{-1}T_sT_{\omega}.$$That these are
equal follows from (\ref{braid1}) and equality of the bracketed
terms; for the latter observe $\omega s\omega^{-1} s \omega\cdot t=t$ for any $t\in\alpha_1^{\vee}({\mathbb
  F}_q^{\times})$.\hfill$\Box$\\

In view of Lemma \ref{2vorstubei} we have natural surjections of $k$-algebras $${\mathcal H}^{\sharp\sharp}\longrightarrow {\mathcal H}^{\sharp}\longrightarrow {\mathcal H}\longrightarrow{\mathcal H}^{\flat}.$$

{\bf Remark:} ${\mathcal H}^{\sharp\sharp}$ (and in particular ${\mathcal H}^{\sharp}$ and ${\mathcal H}$) is generated as a $k$-algebra by $T_{\omega}^{\pm 1}$, $T_s$ and the $T_{e^*(\gamma)}$ for $\gamma\in\Gamma$.\\

\begin{lem}\label{involu} There are unique $k$-algebra involutions $\iota$ of ${\mathcal H}$, ${\mathcal H}^{\sharp}$ and ${\mathcal H}^{\sharp\sharp}$ with
$$\iota(T_{\omega})= T_{\omega},\quad\quad  \iota(T_s)=\tau_s-T_s, \quad\quad  \iota(T_t)=T_t\quad \mbox{  for }t\in\overline{T}.$$
\end{lem}

{\sc Proof:} This is a slightly tedious but straightforward computation. (For
${\mathcal H}$ see \cite{vigneras} Corollary 2.)\hfill$\Box$\\

{\bf Remark:} Besides $\iota$ consider the $k$-algebra involution $\beta$ of
${\mathcal H}$, ${\mathcal H}^{\sharp}$ and ${\mathcal H}^{\sharp\sharp}$
given on generators by$$\beta(T_{\omega})= T_{\omega}^{-1},\quad\quad
\beta(T_s)=T_s, \quad\quad  \beta(T_t)=T_{s\cdot t}\quad \mbox{  for
}t\in\overline{T}.$$Moreover, for any automorphism ${\mathfrak o}$ of $\Gamma$
there is an associated automorphism $\alpha_{\mathfrak o}$ of ${\mathcal H}$,
${\mathcal H}^{\sharp}$ and ${\mathcal H}^{\sharp\sharp}$ given on generators by$$\alpha_{\mathfrak o}(T_{\omega})= T_{\omega},\quad\quad  \alpha_{\mathfrak
  o}(T_s)=T_s, \quad\quad  \alpha_{\mathfrak o}(T_{\partial(\gamma)})=
T_{\partial({\mathfrak o}(\gamma))}\quad \mbox{  for }\gamma\in\Gamma,
  \partial\in{\rm Hom}(\Gamma,\overline{T}).$$Do $\iota$, $\beta$ and the $\alpha_{\mathfrak o}$ generate the automorphism group of ${\mathcal H}$ (resp. of ${\mathcal H}^{\sharp}$, resp. of ${\mathcal H}^{\sharp\sharp}$) modulo inner automorphisms ?\\ 

\begin{lem}\label{quadrat} Let $M$ be an ${\mathcal H}^{\sharp\sharp}$-module. We have a direct sum decomposition$$M=M^{T_s=-{\rm
      id}}\bigoplus M^{T_s^2=0}.$$
\end{lem}

{\sc Proof:} One computes $\tau_s^2=(q-1)\tau_s=-\tau_s$ and this shows $T_s=-{\rm id}$ on ${\rm im}(T^2_s)$ as well as $T_s^2=0$ on ${\rm im}(T^2_s-{\rm id})$.\hfill$\Box$\\

Let ${[0,q-2]}^{\Phi}$ be the set of tuples ${\epsilon}=(\epsilon_{i})_{0\le
  i\le d}$ with $\epsilon_i\in\{0,\ldots, q-2\}$ and $\sum_{0\le i\le
  d}{\epsilon}_{i}\equiv 0$ modulo $(q-1)$. We often read the indices as
elements of ${\mathbb Z}/(d+1)$, thus $\epsilon_i=\epsilon_j$ for $i,
j\in{\mathbb Z}$ whenever $i-j\in(d+1){\mathbb Z}$. We let the symmetric group
${\mathfrak S}_{d+1}$ (generated by $s$, $\omega$ as before) act on ${[0,q-2]}^{\Phi}$ as follows:$$({\omega}\cdot\epsilon)_0= \epsilon_{d}\quad\quad\mbox{ and
 }\quad\quad({\omega}\cdot\epsilon)_i= \epsilon_{i-1} \quad \mbox{ for }1\le i\le d.$$If $d=1$ we put$$(s\cdot\epsilon)_i=-\epsilon_i\quad\mbox{ for
 }i=0,1,$$but if $d\ge2$ we put$$(s\cdot\epsilon)_1=-\epsilon_1,\quad
 (s\cdot\epsilon)_0=\epsilon_0+\epsilon_1,\quad
 (s\cdot\epsilon)_2=\epsilon_1+\epsilon_2,\quad (s\cdot\epsilon)_i=\epsilon_i\quad\mbox{ for
 }3\le i\le d.\footnote{Here and below we understand $-\epsilon_i$ to mean the representative in $[0,q-2]$ of the class of $-\epsilon_i$ in ${\mathbb Z}/(q-1)$, and similarly for $\epsilon_0+\epsilon_1$ and $\epsilon_1+\epsilon_2$.}$$

Throughout we assume that all
eigenvalues of the $T_{t}$ for $t\in{\overline{T}}$ acting on an ${\mathcal H}^{\sharp\sharp}$-module belong to $k$.

Let $M$ be an ${\mathcal H}^{\sharp\sharp}$-module. For $a\in[0,q-2]$ and ${\epsilon}=(\epsilon_{i})_{0\le i\le d}
\in{[0,q-2]}^{\Phi}$ and $j\in\{0,1\}$ put \begin{align}M^{{\epsilon}}&=\{x\in
M\,\,|\, T^{-1}_{\alpha_i^{\vee}(\gamma)}(x)=\gamma^{{\epsilon}_{i}}x\mbox{ for all }\gamma\in \Gamma,\mbox{ all }0\le i\le d\},\notag\\M^{\epsilon}_{\underline{a}}&=\{x\in M^{\epsilon}\,|\,T_{e^*(\gamma)}(x)=\gamma^ax\mbox{ for all }\gamma\in\Gamma\},\notag\\M^{\epsilon}_{\underline{a}}[j]&=\{x\in M^{\epsilon}_{\underline{a}}\,|\,T_{s}^2(x)=jx\}.\notag\end{align}

The $T_t$ for $t\in\overline{T}$ are of order divisible by $q-1$, hence are diagonalizable on the $k$-vector space $M$. Since they commute among each other and with $T_s^2$, we may simultaneously diagonalize all these operators (see Lemma \ref{quadrat} for $T_s^2$), hence\begin{gather}M=\bigoplus_{\epsilon, a, j}M^{\epsilon}_{\underline{a}}[j].\label{18okt}\end{gather}

\begin{lem}\label{samserswo} For any $\epsilon\in [0,q-2]^{\Phi}$ and $a\in
  [0,q-2]$ we have$$T_{\omega}(M^{\epsilon}_{\underline{a}})=M^{\omega\cdot\epsilon}_{\underline{a-\epsilon_0}}\quad\quad\mbox{ and }\quad\quad T_{s}(M^{\epsilon})\subset
  M^{s\cdot\epsilon}.$$If $M$ is even an ${\mathcal
  H}$-module then\begin{gather} T_s(M^{\epsilon}_{\underline{a}})\subset M^{s\cdot\epsilon}_{\underline{\epsilon_1+a}}.\label{vorvigilallsaints}\end{gather}
\end{lem}

{\sc Proof:} $T_{\omega}(M^{\epsilon})=M^{\omega\cdot\epsilon}$ and $T_{s}(M^{\epsilon})\subset
  M^{s\cdot\epsilon}$ follows from formulas (\ref{hsvtabfue}) and (\ref{hsvtabfue1}) respectively, for the $t=\alpha^{\vee}_{i}(\gamma)$. For
the following computation recall that $\omega\cdot e^{*}=e^{*}+\alpha_0^{\vee}$:
For $\gamma\in\Gamma$ and $x\in M^{\epsilon}_{\underline{a}}$ we
have$$T_{e^*(\gamma)}T_{\omega}(x)=T_{\omega}T_{(\omega\cdot e^*)(\gamma)}(x)=T_{\omega}T_{e^*(\gamma)}T_{\alpha_0^{\vee}(\gamma)}(x)=\gamma^{a-\epsilon_0}T_{\omega}(x).$$This shows $T_{\omega}(M^{\epsilon}_{\underline{a}})=M^{\omega\cdot\epsilon}_{\underline{a-\epsilon_0}}$. For formula (\ref{vorvigilallsaints}) recall that $s\cdot e^{*}=e^{*}-\alpha_1^{\vee}$ and employ formula (\ref{hsvtabfue1}). \hfill$\Box$\\

Any $x\in M$ can be uniquely written as $$x=\sum_{a\in[0,q-2]}x_{\underline{a}}\quad \mbox{ with }\quad x_{\underline{a}}\in\sum_{\epsilon\in [0,q-2]^{\Phi}} M^{\epsilon}_{\underline{a}}.$$Given $a\in{\mathbb Z}$ and $x\in M$, we write $x_{\underline{a}}=x_{\underline{\widetilde{a}}}$ where $\widetilde{a}\in[0,q-2]$ is determined by $a-\widetilde{a}\in(q-1){\mathbb Z}$.\\

{\bf Definition:} (a) An ${\mathcal
  H}^{\sharp}$-module $M$ is called {\it standard supersingular} if the ${\mathcal
  H}^{\sharp}$-action factors through ${\mathcal H}$, making it a standard supersingular ${\mathcal
  H}$-module. 

 (b) An irreducible ${\mathcal
  H}^{\sharp}$-module is called {\it supersingular} if it is a subquotient of a standard supersingular ${\mathcal
  H}^{\sharp}$-module. An ${\mathcal H}^{\sharp}$-module $M$ is called {\it supersingular} if it is the
inductive limit of finite dimensional ${\mathcal
  H}^{\sharp}$-modules and if each of its irreducible subquotients is supersingular.

 (c) An ${\mathcal
  H}^{\sharp\sharp}$-module $M$ is called {\it supersingular} if it satisfies
the condition analogous to (b).

(d) A supersingular ${\mathcal
    H}^{\sharp}$-module is called {\it $\sharp$-supersingular} if for all
$e\in M^{\epsilon}_{\underline{a}}[0]$ with $\epsilon_1>0$ we have
$$(T_se)_{\underline{c+\epsilon_1+a}}=0\quad\quad \mbox{ for all }\quad q-1-\epsilon_1\le c\le q-2.$$

\begin{lem}\label{vialsa} (a) An ${\mathcal
  H}$-module is supersingular if and only if it is supersingular when viewed as an ${\mathcal
    H}^{\sharp}$-module. A supersingular ${\mathcal
  H}$-module is $\sharp$-supersingular when viewed as an ${\mathcal
    H}^{\sharp}$-module.

  (b) The category of supersingular ${\mathcal
  H}$-modules, the category of supersingular ${\mathcal
  H}^{\sharp}$-modules, the category of supersingular ${\mathcal
  H}^{\sharp\sharp}$-modules and the category of $\sharp$-supersingular ${\mathcal
  H}^{\sharp}$-modules are abelian.
\end{lem}

{\sc Proof:} Statement (a) follows from formula
(\ref{vorvigilallsaints}). Statement (b) is clear from the definitions.\hfill$\Box$\\

\section{Reconstruction of supersingular ${\mathcal H}^{\sharp}$-modules}

\label{may2019}

Given an ${\mathcal
    H}^{\sharp}$-module $M$ together with a submodule $M_0$ such that $M/M_0$ is supersingular, we address the problem of reconstructing the ${\mathcal
    H}^{\sharp}$-module $M$ from the ${\mathcal
    H}^{\sharp}$-modules $M_0$ and $M/M_0$ together with an additional set of data (intended to be sparse). Our proposed solution (Proposition \ref{schutzengel}) critically relies on the braid relations (\ref{braid11}) and (\ref{braid22}).

\begin{lem}\label{brixen} Let $B_0,\ldots,B_n$ be linear operators on a
  $k$-vector space $M$ such that\begin{gather}B_j^2=B_j\quad\mbox{ for all }
  0\le j\le n,\notag\\B_jB_{j'}B_j=B_{j'}B_jB_{j'}\quad\mbox{ for all }0\le
  j',j\le n, \notag\\B_jB_{j'}=B_{j'}B_j\quad\mbox{ for all }0\le j'< j\le
  n\mbox{ with }j-j'\ge 2. \notag\end{gather}Put
$\beta=B_n\cdots B_{1}B_0$ and let $x\in M$ with $\beta^mx=x$ for some $m\ge1$. Then we have $B_jx=x$ for
  each $0\le j\le n$. 
\end{lem}

{\sc Proof:} We first claim\begin{gather}\beta
B_{j+1}=B_{j}\beta\quad\quad \mbox{ for
  all }0\le j< n.\label{vorpiz}\end{gather}Indeed, \begin{align}\beta
  B_{j+1}&=B_n\cdots  B_{j+2}B_{j+1} B_{j}B_{j-1}\cdots B_{1}B_0    B_{j+1}\notag\\{}&=B_n\cdots  B_{j+2}B_{j+1} B_{j} B_{j+1}B_{j-1}   \cdots
  B_{1}B_0\notag\\{}&=B_n\cdots  B_{j+2}  B_{j} B_{j+1}B_{j} B_{j-1}   \cdots
  B_{1}B_0\notag\\{}&=B_{j}\beta.\notag\end{align}Choose $\nu\ge1$ with
  $m\nu\ge n$. For $0\le j\le n$ we then
  compute$$x\stackrel{(i)}{=}\beta^{m\nu}x=\beta^{n-j}\beta^{m\nu-n+j}x\stackrel{(ii)}{=}\beta^{n-j}B_n\beta^{m\nu-n+j}
  x\stackrel{(iii)}{=}B_j\beta^{n-j}\beta^{m\nu-n+j} x{=}B_j\beta^{m\nu} x\stackrel{(iv)}{=}B_j x,$$where $(i)$
  and $(iv)$ follow
  from the hypothesis $\beta^m x=x$, where $(ii)$ follows
  from $B_n\beta=\beta$ and where $(iii)$ follows from repeated application of
  formula (\ref{vorpiz}).\hfill$\Box$\\

\begin{pro}\label{novordschutzengel} Let $M$ be an ${\mathcal
    H}^{\sharp}$-module, let $M_0\subset M$ be an ${\mathcal
    H}^{\sharp}$-submodule such that $M/M_0$ is supersingular. Let
  $\overline{x}\in (M/M_0)^{\epsilon}$ (some $\epsilon\in [0,q-2]^{\Phi}$) be such that
  $\overline{x}\{i\}=T_{\omega}^{i+1}\overline{x}$ is an eigenvector under
  $T_s$, for each $i\in{\mathbb Z}$. For liftings $x\in M$ of $\overline{x}$ put $x\{i\}=T_{\omega}^{i+1}x$.

(a) If the ${\mathcal
    H}^{\sharp}$-action on $M$ factors through ${\mathcal
    H}$ then we may choose $x\in M^{\epsilon}$ such that for each $i$ with
  $T_s(\overline{x}\{i\})= 0$ and $(\omega^{i+1}\cdot\epsilon)_1=0$ we have $T_s(x\{i\})=0$.

(b) If the ${\mathcal
    H}^{\sharp}$-action on $M$ factors through ${\mathcal
    H}$ then we may choose $x\in M^{\epsilon}$ such that for each $i$ with
  $T_s(\overline{x}\{i\})= -\overline{x}\{i\}$ we have $T_s(x\{i\})=-x\{i\}$.

(c) We may choose $x\in M^{\epsilon}$ such that for each $i$ with
  $T_s^2(\overline{x}\{i\})= 0$ we have $T_s^2(x\{i\})=0$.

(d) We may choose $x\in M^{\epsilon}$ such that for each $i$ with
  $T_s^2(\overline{x}\{i\})=\overline{x}\{i\}$ we have $T_s^2(x\{i\})=x\{i\}$.

\end{pro}

{\sc Proof:} (a) Let $i_1<\ldots <i_r$ be the increasing enumeration of the
set of all $0\le i\le d$ with $T_sT_{\omega}^{i+1}(\overline{x})=0$ and
$(\omega^{i+1}\cdot\epsilon)_1=0$. Replacing $M$ by its ${\mathcal
    H}^{\sharp}$-submodule generated by $x$ and $M_0$ we may assume that
$M/M_0$ is a subquotient of a standard supersingular ${\mathcal
    H}$-module, attached to a supersingular character $\chi:{\mathcal
  H}_{\rm aff}\to k$. If we had $T_sT_{\omega}^{i+1}(\overline{x})=0$ and
$(\omega^{i+1}\cdot\epsilon)_1=0$ for all $0\le i\le d$ then this would mean
$\chi(T_{\omega}^{m}T_sT_{\omega}^{-m})=0$ and
$\chi(T_{\omega}^{m}{\tau}_sT_{\omega}^{-m})\ne0$ for all $m\in{\mathbb Z}$,
in contradiction with the supersingularity of $\chi$. Hence there is some
$0\le i\le d$
not occuring among $\{i_1,\ldots,i_r\}$. Thus, after a cyclic index shift, we may assume $i_r<d$. 

Start with an arbitrary lift $x\in M^{\epsilon}$ of $\overline{x}$. 

We claim that for any $j$ with $0\le j\le r$, after modifying $x$ if
necessary, we can achieve $T_s(x\{i_s\})=0$ for all $s$ with $1\le
s\le j$. For $j=r$ this is the desired statement.

Let us illustrate the argument in the case $d=1$ first. (This will
logically not
be needed for the general case. Notice e.g. that the subarguments (2) and (3)
below are required only if $d>1$.) Then we have $r=1$ and $i_1=0$, and the claim for $j=1$
states that there is some $\tilde{x}\in M^{\epsilon}$ lifting $\overline{x}$
with $T_sT_{\omega}(\tilde{x})=0$. But indeed,
$\tilde{x}=x+T_{\omega}^{-1}T_sT_{\omega}x$ works: First, $\tilde{x}$
lifts $\overline{x}$ because of $T_sT_{\omega}\overline{x}=0$. Next, $\tilde{x}$ belongs to
$M^{\epsilon}$ because of $T_{\omega}^{-1}T_sT_{\omega}x\in M^{\epsilon}$ (which
follows from ${x}\in M^{\epsilon}$ and the assumption
$(\omega^{i_{1}+1}\cdot\epsilon)_1=0$). Finally,
$T_sT_{\omega}(\tilde{x})=0$, because $T_{\omega}\tilde{x}\in
M^{\omega\cdot\epsilon}=M^{\omega^{i_{1}+1}\cdot\epsilon}$ and
$(\omega^{i_{1}+1}\cdot\epsilon)_1=0$ imply $(T_s+T_s^2)T_{\omega}\tilde{x}=0$.

Now let us consider the case of a general $d$. Induction on $j$. For $j=0$ there is nothing to do. Now fix $1\le j\le r$ and assume that $x$
satisfies the condition for $j-1$, i.e. assume $T_s(x\{i_s\})=0$ for all $s$ with $1\le
s\le j-1$. For $-1\le i\le d$ and
$0\le m<j$ define inductively $$x\{i\}_0=x\{i\}=T_{\omega}^{i+1}x,$$$$x\{i\}_{m+1}=T_{\omega}^{i-i_{j-m}}T_s(x\{i_{j-m}\}_m).$$ We establish several subclaims.

(1) $x\{i\}_m\in M^{\omega^{i+1}\cdot \epsilon}$. 

For $m=0$ there is nothing to do. Next, if the claim is true
for an arbitrary $m$, then we have in particular $x\{i_{j-m}\}_m\in M^{\omega^{i_{j-m}+1}\cdot \epsilon}$. By assumption
we know $(\omega^{i_{j-m}+1}\cdot\epsilon)_1=0$, which implies
$T_s(M^{\omega^{i_{j-m}+1}\cdot\epsilon})\subset
M^{\omega^{i_{j-m}+1}\cdot\epsilon}$. Thus, we get
$x\{i_{j-m}\}_{m+1}=T_s(x\{i_{j-m}\}_m)\in
M^{\omega^{i_{j-m}+1}\cdot\epsilon}$. From this we get $x\{i\}_{m+1}=T_s(x\{i\}_m)\in
M^{\omega^{i+1}\cdot\epsilon}$ for general $i$ by applying powers of $T_{\omega}$ to $x\{i_{j-m}\}_{m+1}$.

(2) $T_s(x\{i_s\}_m)=0$ for all $1\le s\le j$ and all $0\le m<j-s$.

Induction on $m$. For $m=0$ this is true by induction hypothesis (on $j$). Now
let $0<m <j-s$ and assume that we know the claim for $m-1$ instead of $m$. In
particular we then know $T_s(x\{i_s\}_{m-1})=0$. We deduce\begin{align}T_s(x\{i_s\}_m)&=T_sT_{\omega}^{i_s-i_{j-m+1}}T_sT_{\omega}^{i_{j-m+1}-i_s}T_{\omega}^{i_s-i_{j-m+1}}(x\{i_{j-m+1}\}_{m-1})\notag\\&=T_sT_{\omega}^{i_s-i_{j-m+1}}T_sT_{\omega}^{i_{j-m+1}-i_s}(x\{i_s\}_{m-1})\notag\\&=T_{\omega}^{i_s-i_{j-m+1}}T_sT_{\omega}^{i_{j-m+1}-i_s}T_s(x\{i_s\}_{m-1})\notag\\&=0\notag\end{align}where
we use the braid relation (\ref{braid2}) (which applies since
$|i_s-i_{j-m+1}|>1$ and $i_r<d$). The induction on $m$ is complete. 

(3) $T_s(x\{i_s\}_m)=0$ for all $1\le s\le j$ and all $j-s+1<m\le
j$.

We induct on $m+s-j$. The induction begins with $m+s-j=2$. By (2) we know
$T_s(x\{i_{j-m+1}\}_{m-2})=0$. Thus, if
$i_{j-m+1}+1<i_{j-m+2}$, the same argument as in (2) shows
$T_s(x\{i_{j-m+1}\}_{m-1})=0$ and hence $x\{i\}_m=0$ for all $i$, and
there is nothing more to do. If however  $i_{j-m+1}+1=i_{j-m+2}$ we
compute\begin{align}T_s(x\{i_{j-m+2}\}_m)&=T_sT_{\omega}T_sT_{\omega}^{-1}T_sT_{\omega}(x\{i_{j-m+1}\}_{m-2})\notag\\&=T_{\omega}T_sT_{\omega}^{-1}T_sT_{\omega}T_s(x\{i_{j-m+1}\}_{m-2})\notag\\&=0\notag\end{align}where
we use the braid relation (\ref{braid1}). This settles the case $m+s-j=2$. For
$m+s-j>2$ we now argue exactly as in (2) again: $T_s(x\{i_s\}_m)=0$ implies
$T_s(x\{i_s\}_{m+1})=0$. The induction is complete. 

(4) $T_s(x\{i_{j-m}\}_{m}+x\{i_{j-m}\}_{m+1})=0$ for all $0\le
m<j$. 

Indeed, by (1) and our assumption
$(\omega^{i_{j-m}+1}\cdot\epsilon)_1=0$ we know that $x\{i_{j-m}\}_m$ is fixed under
$T_{\alpha_1^{\vee}(\Gamma)}$ and hence is killed by $T_s^2+T_s$, as follows
from the quadratic relation (\ref{sorge}). As
$x\{i_{j-m}\}_{m+1}=T_s(x\{i_{j-m}\}_m)$ this gives the claim.

(5) $$\tilde{x}=\sum_{0\le m\le j}x\{-1\}_m$$lifts $\overline{x}$. 

Indeed, we have $T_s(x\{i_j\})\in M_0$ by our defining assumption on $i_j$. It follows that $x\{-1\}_{m}\in M_0$ for all $m\ge1$, hence $x-\tilde{x}\in M_0$.

(6) From (1) we deduce $\tilde{x}\{i\}\in M^{\omega^{i+1}\cdot\epsilon}$. Writing$$\tilde{x}\{i_s\}=(\sum_{0\le m<
        j-s}x\{i_s\}_m)+(x\{i_s\}_{j-s}+x\{i_s\}_{j-s+1})+(\sum_{j-s+1<m\le j}x\{i_s\}_m)$$we see that (2), (3) and (4) imply $T_s(\tilde{x}\{i_s\})=0$ for all $s$ with $1\le s\le j$. 

The induction on $j$ is complete: we may substitute $\tilde{x}$ for the old $x$. 

(b) Composing the given ${\mathcal H}$-module
structure on $M$ with the involution $\iota$ of Lemma \ref{involu} we get a new ${\mathcal H}$-module structure on
$M$. Applying statement (a) to this new ${\mathcal H}$-module and then
translating back via $\iota$, we get statement (b). Notice that here, in contrast
to the setting in (a), we
{\it automatically} have $(\omega^{i+1}\cdot\epsilon)_1=0$ for each $i$ with
  $T_s(\overline{x}\{i\})= -\overline{x}\{i\}$.

(c) Statement (c) is proved in the same way as statement (a), with the
following minor modifications: Each occurence of $T_s$ must be replaced by
$T_s^2$, and in the definition of $x\{i\}_{m+1}$ the alternating sign
$(-1)^{m+1}$ must be included,
i.e.\begin{gather}x\{i\}_{m+1}=(-1)^{m+1}T_{\omega}^{i-i_{j-m}}T_s^2(x\{i_{j-m}\}_m)\label{nachbei}\end{gather}In
particular, we then have $x\{i_{j-m}\}_{m+1}=-T^2_s(x\{i_{j-m}\}_m)$. In (2)
and (3), the
appeal to the braid relations (\ref{braid1}), (\ref{braid2}) must be replaced
by an appeal to the braid relations (\ref{braid11}), (\ref{braid22}). In (4), the
appeal to $T_s^2+T_s=0$ on vectors fixed under
$T_{\alpha_1^{\vee}(\Gamma)}$ must be replaced by an appeal to $T_s^4-T_s^2=0$ (it
is here where the alternating sign in the defining formula (\ref{nachbei}) is
needed). Notice that here, in contrast
to the setting in (a), we do not need to impose $(\omega^{i+1}\cdot\epsilon)_1=0$ for each $i$ with
$T_s^2(\overline{x}\{i\})= 0$. (On the one hand, because of $T_s^2(M^{\epsilon})\subset
M^{\epsilon}$ for {\it any} $\epsilon$ the argument analogous to the one
in (a)(1) carries over; on the other hand, because of $T_s^4-T_s^2=0$ on all
of $M$ the argument analogous to the one
in (a)(4) carries over.)

(d) Composing the given ${\mathcal H}^{\sharp}$-module
structure on $M$ with the involution $\iota$ of Lemma \ref{involu} we get a new ${\mathcal H}^{\sharp}$-module structure on
$M$. Applying statement (c) to this new ${\mathcal H}^{\sharp}$-module and then
translating back via $\iota$, we get statement (d).\hfill$\Box$\\

\begin{pro}\label{schutzengel} Let $M$ be an ${\mathcal
    H}^{\sharp}$-module, let $M_0\subset M$ be an ${\mathcal
    H}^{\sharp}$-submodule such that $M/M_0$ is supersingular. The action of
  ${\mathcal H}^{\sharp}$ on $M$ is uniquely determined by the following combined data:

(a) the action of ${\mathcal H}^{\sharp}$ on $M_0$ and on $M/M_0$,

(b) the action of $T_{e^*(\Gamma)}$ and of $T_sT_{\omega}$ on $M$,

(c) the restriction of $T_{\omega}$ to $(T_sT_{\omega})^{-1}(M_0)$, i.e. the map $$\{x\in M\,|\,T_sT_{\omega}(x)\in M_0\}\stackrel{T_{\omega}}{\longrightarrow} M,$$

(d) the subspace $\sum_{\epsilon\in
    [0,q-2]^{\Phi}\atop\epsilon_1=0}M^{\epsilon}$ of $M$.

\end{pro}
 
{\sc Proof:} The $k$-algebra ${\mathcal
    H}^{\sharp}$ is generated by $T_{e^*(\Gamma)}$, by $T_s$ and by
$T^{\pm 1}_{\omega}$. Therefore we only need to see that the action of $T_s$ and
$T_{\omega}$ on $M$ can be reconstructed from the given data (a), (b), (c), (d). Exhausting $M/M_0$ step by step we may assume that $M/M_0$ is an
irreducible supersingular ${\mathcal
  H}^{\sharp}$-module. 

We first show that $T_s$ is uniquely determined. For this we make constant use
of Lemma \ref{quadrat} (and the decomposition (\ref{18okt})). As $T_s|_{M_0}$ is given to us, it is enough to show that for any non-zero $\overline{x}$ in $M/M_0$ with
either $T_s(\overline{x})=-\overline{x}$ or $T_s(\overline{x})=0$ we find
some lifting $x\in M$ such that $T_s(x)$ can be reconstructed. Consider
first the case $T_s(\overline{x})=-\overline{x}$. By the quadratic relation (\ref{sorge}) (see Lemma \ref{quadrat}) we then have $\overline{x}\in\sum_{\epsilon\in
    [0,q-2]^{\Phi}\atop\epsilon_1=0}(M/M_0)^{\epsilon}$, and using the datum
(d) as well as our knowledge of the subspace $T_sM$ (since $T_sM=T_sT_{\omega}M$ this is given to us
in view of datum (b)), we lift $\overline{x}$ to some $x\in T_sM\cap \sum_{\epsilon\in
    [0,q-2]^{\Phi}\atop\epsilon_1=0}M^{\epsilon}$ (use the decomposition (\ref{18okt})). For such $x$ we have
$T_s(x)=-x$. Now consider the case $T_s(\overline{x})=0$. An arbitrary lifting
$x\in M$ of $\overline{x}$ then satisfies $T_s(x)\in M_0$, and $T_s(x)$ is determined by the
given data as
$T_s(x)=(T_sT_{\omega})T_{\omega}^{-1}(x)$ (notice that the datum (c) is
equivalent with the datum $T_s^{-1}(M_0)\stackrel{T_{\omega}^{-1}}{\rightarrow} M$). 

To show that $T_{\omega}$ is uniquely determined, suppose that besides $T_{\omega}\in{\rm Aut}_k(M)$ there is another candidate
$\tilde{T}_{\omega}\in{\rm Aut}_k(M)$ extending the data (a), (b), (c), (d) to another ${\mathcal H}^{\sharp}$-action on $M$. 

We find and choose some non-zero $\overline{x}\in M/M_0$ such that
$T_{\omega}^j(\overline{x})$ is an eigenvector under $T_s$, for each
$j\in{\mathbb Z}$. For any 
$x\in M$ lifting $\overline{x}$ we
have\begin{gather}T_{\omega}=\tilde{T}_{\omega}\quad\mbox{ on
}M_0+k.T_{\omega}^{j-1}(x)\quad\mbox{ if
}T_sT_{\omega}^j(\overline{x})=0\label{nachorkan}\end{gather}as both
$\tilde{T}_{\omega}$ and ${T}_{\omega}$ respect the datum (c).

Let $i_0<\ldots<i_n$ be the increasing enumeration of the set $$\{0\le i\le
d\,|\,T_s^2{T}_{\omega}^i\overline{x}={T}_{\omega}^i\overline{x}\}.$$As
$M/M_0$ is a subquotient of a standard supersingular ${\mathcal H}$-module,
this set is not the full set $\{0\le i\le d\}$. Applying a suitable power of $T_{\omega}$ and reindexing we may assume
that $0$ does not belong to this set, i.e. that $i_0>0$.

Choose a lifting $x\in M$ of $\overline{x}$ such that for each $i\in\{i_0,\ldots,i_n\}+{\mathbb Z}(d+1)$ we have
$T_s^2{T}_{\omega}^ix={T}_{\omega}^ix$. This is possible by Proposition
\ref{novordschutzengel}. Put $z_{0}=x$. For $i\ge 1$
put\begin{gather}z_i=\left\{\begin{array}{l@{\quad:\quad}l}\tilde{T}_{\omega}z_{i-1}
& i\notin\{i_0,\ldots,i_n\}+{\mathbb Z}(d+1) \notag\\T_s^{2}\tilde{T}_{\omega}z_{i-1} &
i\in\{i_0,\ldots,i_n\}+{\mathbb Z}(d+1)\end{array}\right.\notag.\end{gather}We
claim \begin{gather}z_i={T}^{i}_{\omega}x\label{nachstur}\end{gather}for each
$i\ge 0$. We induct on $i$. The case $i=0$ is trivial. For $i\ge1$ with
$i\notin\{i_0,\ldots,i_n\}+{\mathbb Z}(d+1)$ we
compute$$z_{i}=\tilde{T}_{\omega}z_{i-1}\stackrel{(i)}{=}{T}_{\omega}z_{i-1}\stackrel{(ii)}{=}{T}^{i}_{\omega}x$$where
in $(i)$ we use statement (\ref{nachorkan})  and in (ii) we use the induction
hypothesis. For $i\ge1$ with $i\in\{i_0,\ldots,i_n\}+{\mathbb Z}(d+1)$ we
compute$$z_{i}=T^2_s\tilde{T}_{\omega}z_{i-1}\stackrel{(i)}{=}T^2_s{T}_{\omega}z_{i-1}\stackrel{(ii)}{=}{T}^{i}_{\omega}x$$where
in $(i)$ we use the assumption $T_sT_{\omega}=T_s\tilde{T}_{\omega}$, and in
(ii) we use the induction hypothesis ${T}_{\omega}z_{i-1}{=}{T}^{i}_{\omega}x$
and the assumption on $x$. The induction is
complete. Put $$B_{i_j}=\tilde{T}_{\omega}^{-i_j}T_s^2\tilde{T}_{\omega}^{i_j}.$$The relation (\ref{sanafabehsv17}) implies $B_{i_j}=\tilde{T}_{\omega}^{-i_j+(d+1)\nu}T_s^2\tilde{T}_{\omega}^{i_j-(d+1)\nu}$
for each $\nu\in{\mathbb Z}$. Thus$$(B_{i_n}\cdots
B_{i_1}B_{i_0})^mx\stackrel{(i)}{=}\tilde{T}_{\omega}^{-m(d+1)}z_{m(d+1)}\stackrel{(ii)}{=}\tilde{T}_{\omega}^{-m(d+1)}{T}_{\omega}^{m(d+1)}x$$for $m\ge 0$, where
$(i)$ follows from the definition of $z_{m(d+1)}$, whereas $(ii)$ follows from formula
(\ref{nachstur}). Choosing $m$ large enough we may assume
${T}_{\omega}^{m(d+1)}x=x$ and $\tilde{T}_{\omega}^{m(d+1)}x=x$ (as
${T}_{\omega}$ and $\tilde{T}_{\omega}$ are automorphisms of a finite vector space); then$$(B_{i_n}\cdots
B_{i_1}B_{i_0})^mx{=}x.$$The braid relations (\ref{braid11}), (\ref{braid22}) show that the
$B_{i_j}$ satisfy the hypotheses of Lemma \ref{brixen} (in particular, the commutation
$B_{i_0}B_{i_n}=B_{i_n}B_{i_0}$ if $n>1$ follows from $i_0>0$). This
Lemma now tells us $B_{i_j}\cdots
B_{i_1}B_{i_0}x=x$ for each $0\le j\le n$. But by the
definition of the $z_i$ this
means \begin{gather}z_i=\tilde{T}^{i}_{\omega}x\label{nachstur1}\end{gather}for
each $0\le i\le d+1$. When compared with formula (\ref{nachstur}) this yields
${T}_{\omega}=\tilde{T}_{\omega}$ since $M$ is generated as a $k$-vector
space by $M_0$ together with the ${T}^{i}_{\omega}x$ (or: the
$\tilde{T}^{i}_{\omega}x$) for $0\le i\le d$.\hfill$\Box$\\

{\bf Remarks:} The above proof of Proposition \ref{schutzengel} shows the following:

(i) The subspace in (d) could be replaced by the subspace $\{x\in M\,|\,T_s^2(x)=x\}$.

(ii) If the ${\mathcal
    H}^{\sharp}$-action factors through an ${\mathcal
    H}$-action, then the datum (d) can be entirely left out ($T_{\omega}$ can
then be reconstructed without a priori knowledge of $T_s$).

\section{The functor}

Here we define a functor $M\mapsto \Delta(M)$ from supersingular ${\mathcal H}^{\sharp\sharp}$-modules to torsion $k[[t]]$-modules with $\varphi$ and $\Gamma$ actions, as outlined in the introduction. Its entire content is encapsulated in the explicit formula for the elements $h(e)$ introduced below.

Let $M$ be an ${\mathcal H}^{\sharp\sharp}$-module. View $M$ as a $k[[t]]$-module with $t=0$
on $M$. Let $\Gamma$ act on $M$ by $$\gamma\cdot x=T^{-1}_{e^*(\gamma)}(x)$$
for $\gamma\in\Gamma$, making $M$ a $k[[t]][\Gamma]$-module. We have an isomorphism of $k[[t]][\varphi]$-modules$${\mathfrak O}\otimes_{k[[t]][\Gamma]}M\cong k[[t]][\varphi]\otimes_{k[[t]]}M$$and hence an action of ${\mathfrak O}$ on $k[[t]][\varphi]\otimes_{k[[t]]}M$.

For $e\in  M^{\epsilon}_{\underline{a}}[j]$ (any $\epsilon\in{[0,q-2]}^{\Phi}$, any $a\in[0,q-2]$, any $j\in\{0,1\}$) define the element\begin{gather}h(e)=\left\{\begin{array}{l@{\quad:\quad}l}t^{\epsilon_1}\varphi\otimes T^{-1}_{{\omega}}(e)+1\otimes
    e +\sum_{c=0}^{q-2}t^{c}\varphi\otimes
    T^{-1}_{{\omega}}((T_se)_{\underline{c+\epsilon_1+a}}) & j=0 \notag\\t^{q-1}\varphi\otimes T^{-1}_{{\omega}}(e)+1\otimes
    e & j=1\end{array}\right.\notag\end{gather}of $k[[t]][\varphi]\otimes_{k[[t]]}M$. Define $\nabla(M)$ to be the $k[[t]][\varphi]$-submodule of
$k[[t]][\varphi]\otimes_{k[[t]]}M$ generated by the elements $h(e)$ for all
    $e\in  M^{\epsilon}_{\underline{a}}[j]$ (all $\epsilon$, $a$,
    $j$). Define $$\Delta(M)=\frac{k[[t]][\varphi]\otimes_{k[[t]]}M}{\nabla(M)}.$$ 

{\bf Remark:} If $M$ is even an ${\mathcal H}$-module, then in view of formula (\ref{vorvigilallsaints}) the definition of
$h(e)$ simplifies to become\begin{gather}h(e)=\left\{\begin{array}{l@{\quad:\quad}l}t^{\epsilon_1}\varphi\otimes T^{-1}_{{\omega}}(e)+1\otimes
    e +\varphi\otimes
    T^{-1}_{{\omega}}(T_se) & j=0 \notag\\t^{q-1}\varphi\otimes T^{-1}_{{\omega}}(e)+1\otimes
    e & j=1\end{array}\right.\notag.\end{gather}In this case it is not
    necessary to split up $M$ into eigenspaces under the action of
    $T_{e^*(\Gamma)}$, and the {\it notation} of many of the subsequent
    computations simplifies (no underlined subscripts are needed). However, they hardly simplify in mathematical complexity, not even if we restrict to ${\mathcal
      H}^{\flat}$-modules only (in which case always $\epsilon_1=0$ and $T_{e^*(\gamma)}=1$).\\

\begin{lem}\label{strasheim} Let $e\in  M^{\epsilon}_{\underline{a}}[j]$. The integer\begin{gather}k_e=\left\{\begin{array}{l@{\quad:\quad}l}\epsilon_{1}& j=0  \notag\\ q-1
  &  j=1  \notag \end{array}\right.\notag\end{gather}satisfies $k_e\equiv\epsilon_1$ modulo $(q-1)$.
\end{lem}

{\sc Proof:} $j=1$ means $T_s^2(e)=e$, hence the claim follows from the relation (\ref{sorge}).\hfill$\Box$\\

\begin{lem}\label{4stubai} For $e\in  M^{\epsilon}_{\underline{a}}[j]$ we have $\gamma\cdot h(e)=h(T^{-1}_{e^*(\gamma)}(e))$ for all $\gamma\in \Gamma$. In particular, $\nabla(M)$ is stable under the action of $\Gamma$, hence is an ${\mathfrak O}$-submodule of $k[[t]][\varphi]\otimes_{k[[t]]}M$. Hence $\Delta(M)$ is even an ${\mathfrak O}$-module.
\end{lem}

{\sc Proof:} First notice that $T^{-1}_{e^*(\gamma)}(e)\in M^{\epsilon}_{\underline{a}}[j]$. In particular, $h(T^{-1}_{e^*(\gamma)}(e))$ is well defined. For $\gamma\in\Gamma$ we find\begin{gather}\gamma\cdot (1\otimes
    e)=1\otimes \gamma\cdot e=1\otimes  T^{-1}_{e^*(\gamma)}(e).\label{7vorstubai1}\end{gather}Next, we compute\begin{align}\gamma\cdot( t^{k_e}\varphi\otimes T^{-1}_{{\omega}}(e))&\stackrel{(i)}{=}\gamma^{k_e}t^{k_e}\varphi\otimes \gamma\cdot T^{-1}_{{\omega}}(e)\notag\\{}&\stackrel{(ii)}{=}t^{k_e}\varphi\otimes T^{-1}_{{\omega}}T^{-1}_{e^*(\gamma)}(e).\label{7vorstubai2}\end{align}In $(i)$ we used $\gamma
t=[\gamma]_{\Phi}(t) \gamma$ and $[\gamma]_{\Phi}(t)\equiv\gamma t$
modulo $t^qk[[t]]$ (Lemma \ref{dominikanerjub}) and the fact that, since
$\pi=0$ in $k$, we have $t^q\varphi\otimes M=\Phi(t)\varphi\otimes M=\varphi
t\otimes M=0$. To see (ii) observe$$\gamma\cdot
T^{-1}_{{\omega}}(e)=T^{-1}_{e^*(\gamma)}
T^{-1}_{{\omega}}(e)=T^{-1}_{{\omega}}T^{-1}_{({\omega}^{-1}\cdot
  e^*)(\gamma)}(e)$$$$=T^{-1}_{{\omega}}T^{-1}_{(e^*-\alpha_1^{\vee})(\gamma)}(e)=T^{-1}_{{\omega}}T_{\alpha_1^{\vee}(\gamma)}T^{-1}_{e^*(\gamma)}(e)=\gamma^{-k_e}T^{-1}_{{\omega}}T^{-1}_{e^*(\gamma)}(e)$$where
in the last step we use Lemma \ref{strasheim}. Combining formulae (\ref{7vorstubai1}) and (\ref{7vorstubai2}) we are done in the case $j=1$. In the case $j=0$ we in addition need the formula\begin{gather}\gamma\cdot\sum_{c=0}^{q-2}t^{c}\varphi\otimes
    T^{-1}_{{\omega}}((T_se)_{\underline{c+\epsilon_1+a}})=\sum_{c=0}^{q-2}t^{c}\varphi\otimes
    T^{-1}_{{\omega}}((T_sT_{e^*(\gamma)}e)_{\underline{c+\epsilon_1+a}}).\label{7vorstubai4}\end{gather}Let us prove this (for $e\in M^{\epsilon}_{\underline{a}}[0]$). For $f\in{\mathbb Z}$ and $\gamma\in\Gamma$ we compute\begin{gather}T_{(\omega^{-1}\cdot
  e^*)(\gamma)}((T_se)_{\underline{f}})\stackrel{(i)}{=}T_{e^*(\gamma)}T_{\alpha_1^{\vee}(\gamma^{-1})}((T_se)_{\underline{f}})\notag\\\stackrel{(ii)}{=}\gamma^{f-\epsilon_1}(T_se)_{\underline{f}}=\gamma^{f-\epsilon_1-a}(T_s(\gamma^a e))_{\underline{f}}=\gamma^{f-\epsilon_1-a}(T_sT_{e^*(\gamma)}e)_{\underline{f}}.\label{nachvor}\end{gather}In $(i)$ recall that $\omega^{-1}\cdot e^*=e^*-\alpha_1^{\vee}$, in $(ii)$ notice that $(T_se)_{\underline{f}}\in M^{s\cdot\epsilon}$ and $(s\cdot\epsilon)_1=-\epsilon_1$. For $c\in[0,q-2]$ we deduce\begin{align}\gamma\cdot (t^c\varphi\otimes T_{\omega}^{-1}((T_se)_{\underline{c+\epsilon_1+a}}))&=\gamma^c t^c\varphi\otimes\gamma\cdot (T_{\omega}^{-1}((T_se)_{\underline{c+\epsilon_1+a}}))\notag\\{}&=\gamma^ct^c\varphi\otimes T^{-1}_{e^*(\gamma)}T_{\omega}^{-1}((T_se)_{\underline{c+\epsilon_1+a}})\notag\\{}&=\gamma^ct^c\varphi\otimes T_{\omega}^{-1}T^{-1}_{(\omega^{-1}\cdot e^*)(\gamma)}((T_se)_{\underline{c+\epsilon_1+a}})\notag\\{}&=t^c\varphi\otimes T_{\omega}^{-1}((T_sT_{e^*(\gamma)}e)_{\underline{c+\epsilon_1+a}})\notag\end{align}where in the last equality we inserted formula (\ref{nachvor}).\hfill$\Box$\\

\begin{pro}\label{jishnusharp} (a) If $M$ is supersingular and
  finite dimensional, then we have: $\Delta(M)$ is a torsion $k[[t]]$-module, generated by $M$ as a
  $k[[t]][\varphi]$-module, and $\varphi$ acts injectively on it. The dual $\Delta(M)^*={\rm Hom}_k(\Delta(M),k)$ is a free $k[[t]]$-module of rank ${\rm
    dim}_k(M)$. The map $M\to \Delta(M)$ which sends $m\in M$ to the class of
  $1\otimes m$ induces a bijection \begin{gather}M\cong \Delta(M)[t].
    \label{sofeendsharp}\end{gather}

(b) $\Delta(M)$ belongs to ${\rm
    Mod}^{\clubsuit}({\mathfrak O})$.
  
(c) The assignment $M\mapsto \Delta(M)$ is an exact
  functor from the category of supersingular ${\mathcal
    H}^{\sharp\sharp}$-modules to ${\rm
    Mod}^{ad}({\mathfrak O})$. 
\end{pro}

{\sc Proof:} (a) Notice first that it is enough to prove these claims after a finite base extensions of $k$.

Assume first that $M$ is irreducible. It can then be realized as a subquotient
of a standard supersingular ${\mathcal H}$-module $N$ --- in fact, it can even be
realized as a submodule or as a quotient of such an $N$. Observing the decomposition
(\ref{18okt}) for $N$, we see that there are a $k$-basis $e_0,\ldots,e_d$ of $N$ as well as $0\le
k_{e_j}\le q-1$ for $0\le j\le d$, not all of them $=0$ and not all of them
$=q-1$, such that $\nabla(N)$ is generated by the elements
$h(e_j)=t^{k_{e_j}}\varphi\otimes T_{\omega}^{-1}(e_j)+1\otimes e_j$. It
follows that $\Delta(N)$ is standard cyclic. Now it is easy to see that $\Delta(M)$ is a subquotient of $\Delta(N)$. Thus, by Proposition \ref{altfund} (d), it is standard cyclic as well, at least after a finite extension of $k$. Therefore  all our
claims follow from Lemma \ref{abstrtor} and Proposition \ref{altfund} (a).

Now let $M$ be arbitrary (supersingular, finite dimensional). Choose a separated
and exhausting descending
filtration of $M$ by ${\mathcal
  H}^{\sharp\sharp}$-submodules $F^{\mu}M$ with irreducible subquotients $F^{\mu-1}M/F^{\mu}M$. Since on any standard supersingular ${\mathcal
  H}$-module (and hence on any of its subquotients, and hence on any irreducible ${\mathcal
  H}^{\sharp\sharp}$-module) we have $T_s=-T_s^2$ and hence ${\rm ker}(T_s^2)={\rm ker}(T_s)$, the filtration satisfies\begin{gather}{T}_s
(F^{\mu-1}M\cap {\rm ker}(T_s^2))\subset F^{\mu}M\label{6vorstubaisharp}\end{gather}for each $\mu\in{\mathbb Z}$. Putting$${\bf
  F}^{\mu}=k[[t]][\varphi]\otimes_{k[[t]]}F^{{\mu}}M$$we claim\begin{gather}\nabla(F^{\mu}M)=\nabla(M)\cap{\bf F}^{\mu}.\label{savovasharp}\end{gather}Arguing by induction, we may assume
that this is known with ${\mu}-1$ instead of ${\mu}$. Let ${\mathcal E}$ be a family of elements $e\in (F^{{\mu}-1}M)^{\epsilon_e}_{\underline{a_e}}[j_e]$ (for suitable
$\epsilon_e\in [0,q-2]^{\Phi}$ and $a_e\in [0,q-2]$ and $j_e\in\{0,1\}$
depending on $e$) which induces a $k$-basis
of $F^{{\mu}-1}M/F^{{\mu}}M$. We consider an
expression\begin{gather}\sum_{j_1,j_2\in {\mathbb Z}_{\ge0},e\in{\mathcal E}}c_{j_1,j_2,e}t^{j_2}\varphi^{j_1}h(e)\label{safa3sharp}\end{gather}with $c_{j_1,j_2,e}\in k$. Assuming that the expression (\ref{safa3sharp})
belongs to ${\bf F}^{\mu}$ we need to see that it even belongs to
$\nabla(F^{\mu}M)$.

Suppose that this is false. We may then define$$j_1={\rm min}\{j\ge 0\,|\,c_{j,j_2,e}t^{j_2}\varphi^{j}h(e)\notin \nabla(F^{\mu}M)\mbox{ for some }j_2\ge0,\mbox{ some } e\in{\mathcal E}\}.$$

{\it Claim: We find some $j_2$ and some $e$ with $c_{j_1,j_2,e}t^{j_2}\varphi^{j_1}h(e)\in {\bf
  F}^{\mu}-\nabla(F^{\mu}M)$.} 

For $e\in {\mathcal E}$ the expression\begin{gather}1\otimes e+t^{k_e}\varphi\otimes T_{\omega}^{-1}(e)\label{freivorvansharp}\end{gather}is congruent to $h(e)$ modulo ${\bf
  F}^{{\mu}}$, in view of $e\in F^{\mu-1}M$ and formula (\ref{6vorstubaisharp}). Therefore, modulo ${\bf F}^{{\mu}}$ the expression
(\ref{safa3sharp})
reads$$\sum_{j_1,j_2,e}c_{j_1,j_2,e}t^{j_2}\varphi^{j_1}\otimes
e+c_{j_1,j_2,e}t^{j_2}\varphi^{j_1}t^{k_e}\varphi\otimes
T_{\omega}^{-1}(e).$$Notice that $\varphi^{j_1}t^{k_e}\varphi\in
k[[t]]\varphi^{j_1+1}$. The claim now follows in view of\begin{gather}\frac{{\bf F}^{{\mu}-1}}{{\bf
    F}^{{\mu}}}=\bigoplus_{j\ge0}k[[t]]\varphi^j\otimes_{k[[t]]}\frac{F^{{\mu}-1}M}{F^{\mu}M}.\label{savovalasharp}\end{gather}

The claim proven, we may argue by induction on the number
of summands in the expression (\ref{safa3sharp}) which do not belong to
$\nabla(F^{\mu}M)$. We may thus assume from the start that the expression
(\ref{safa3sharp}) consists of a single summand $t^{j_2}\varphi^{j_1}h(e)$, and
that moreover
$e\notin {F}^{\mu}M$ for this $e$. The aim is then to deduce $t^{j_2}\varphi^{j_1}h(e)\in \nabla(F^{\mu}M)$, which contradicts our
 above assumption.

Let us write $\epsilon=\epsilon_e$ and $a=a_e$. The vanishing of
$t^{j_2}\varphi^{j_1}h(e)$ modulo ${\bf F}^{\mu}$ means, by the
decomposition (\ref{savovalasharp}) again,
that$$t^{j_2}\varphi^{j_1}\otimes
e\stackrel{(i)}{=}0\stackrel{(ii)}{=}t^{j_2}\varphi^{j_1}t^{k_e}\varphi\otimes
T_{\omega}^{-1}(e)$$(i.e. absolute vanishing, not just modulo
${\bf F}^{\mu}$). If $T_s^2(e)=e$ then this shows $t^{j_2}\varphi^{j_1} h(e)=0$. Now suppose $T_s^2(e)=0$ (and hence $k_e<q-1$). The definition of $h(e)$ together with the vanishings $(i)$ and $(ii)$ shows$$t^{j_2}\varphi^{j_1}h(e)=t^{j_2}\varphi^{j_1}\sum_{c=0}^{q-2}t^c\varphi \otimes
T_{\omega}^{-1}((T_se)_{\underline{c+\epsilon_1+a}}).$$Since the vanishing $(ii)$ also forces
$t^{j_2}\varphi^{j_1}t^{k_e}\varphi\in k[[t]]\varphi^{j_1+1}t$, there is some $i$ and
  some $j'_2\ge0$ with $$t^{j_2}\varphi^{j_1}=t^{j'_2}\varphi^{j_1}t^i\quad \mbox{ and }\quad i\ge q-k_e.$$If
  $k_e=0$ (and hence $i\ge q$) then again the conclusion is
  $t^{j_2}\varphi^{j_1} h(e)=0$. It remains to discuss the case where
  $0<k_e<q-1$. In this case, $(T_se)_{\underline{c+\epsilon_1+a}}\in
  M^{s\cdot\epsilon}$ and $(s\cdot\epsilon)_1=-\epsilon_1$ implies $q-1-k_e=k_{(T_se)_{\underline{c+\epsilon_1+a}}}$
  for each $c$. We thus see\begin{align}t^{q-k_e+c}\varphi \otimes
  T_{\omega}^{-1}((T_se)_{\underline{c+\epsilon_1+a}})&=t^{1+c}(t^{k_{(T_se)_{\underline{c+\epsilon_1+a}}}}\varphi\otimes
  T_{\omega}^{-1}((T_se)_{\underline{c+\epsilon_1+a}})+1\otimes
  (T_se)_{\underline{c+\epsilon_1+a}})\notag\\{}&=t^{1+c}
  h((T_se)_{\underline{c+\epsilon_1+a}})-\sum_{c'=0}^{q-2}t^{1+c+c'}\varphi\otimes
  T_{\omega}^{-1}((T_s((T_se)_{\underline{c+\epsilon_1+a}}))_{\underline{c'+c+a}})\notag\end{align}by
  the definition of $h((T_se)_{\underline{c+\epsilon_1+a}})$, again since
  $(T_se)_{\underline{c+\epsilon_1+a}}\in M^{s\cdot\epsilon}$ and
  $(s\cdot\epsilon)_1=-\epsilon_1$. For $0\le f\le q-2$ we have$$\sum_{0\le
    c,c'\le q-2\atop c+c'=f}(T_s((T_se)_{\underline{c+\epsilon_1+a}}))_{\underline{f+a}}=\sum_{0\le
    c\le q-2}(T_s((T_se)_{\underline{c+\epsilon_1+a}}))_{\underline{f+a}}=0$$as follows
  from $T^2_s(e)=0$. This shows$$\sum_{c,c'=0}^{q-2}t^{1+c+c'}\varphi\otimes T_{\omega}^{-1}((T_s((T_se)_{\underline{c+\epsilon_1+a}}))_{\underline{c'+c+a}})=0.$$Since $e$ belongs to $F^{\mu-1}M$, formula (\ref{6vorstubaisharp}) shows $h((T_s e)_{\underline{c+\epsilon_1+a}})\in\nabla(F^{\mu}M)$. Together we obtain $t^{q-k_e+c}\varphi \otimes
  T_{\omega}^{-1}((T_se)_{\underline{c+\epsilon_1+a}})\in \nabla(F^{\mu}M)$, hence $t^{i+c}\varphi\otimes
  T_{\omega}^{-1}((T_se)_{\underline{c+\epsilon_1+a}})\in \nabla(F^{\mu}M)$ for $0\le c\le q-2$. This gives$$t^{j_2}\varphi^{j_1}h(e)= \sum_{c=0}^{q-2}t^{j'_2}\varphi^{j_1}t^{i+c}\varphi  \otimes T_{\omega}^{-1}((T_s e)_{\underline{c+\epsilon_1+a}})\in \nabla(F^{\mu}M),$$as desired. 

Formula (\ref{savovasharp}) is proven. It allows us to deduce all our claims
for $M$ from the corresponding claims for the $F^{\mu-1}M/F^{\mu}M$; but for them
they have already been established above.

(b) For each irreducible supersingular ${\mathcal
  H}$-module $M$, extending $k$ if necessary, $\Delta(M)$ admits a filtration such that each associated graded piece is a standard cyclic object in ${\rm Mod}^{\rm ad}({\mathfrak O})$, as pointed out above. Since the functor $\Delta$ is exact (see statement (c)) it therefore takes finite dimensional supersingular ${\mathcal
  H}^{\sharp\sharp}$-modules to objects in ${\rm
  Mod}^{\clubsuit}({\mathfrak O})$.

(c) It is clear that $M\mapsto \Delta(M)$ is a (covariant) right exact functor. To
see left exactness, let $M_1\to M_2$ be injective. Since the kernel of
$\Delta(M_1)\to \Delta(M_2)$ is a torsion $k[[t]]$-module it has, if non zero,
a non zero vector killed by $t$. By formula (\ref{sofeendsharp}) it must belong to (the image of)
$M_1$, contradicting the injectivity of $M_1\to M_2$.

\hfill$\Box$\\


\section{Standard objects and full faithfulness}

\subsection{The bijection between standard supersingular Hecke modules and standard cyclic Galois representations}

\label{vialhe}

Let $M$ be a standard supersingular ${\mathcal
  H}$-module, arising from the supersingular character $\chi:{\mathcal
  H}_{\rm aff}\to k$. There is some $e_0\in M$ such that, putting $e_j=T^{-j}_{{\omega}}e_0$, we have $M=\bigoplus_{j=0}^dk.e_j$ and ${\mathcal
  H}_{\rm aff}$ acts on $k.e_0$ by $\chi$. Denote by $\eta_j:\Gamma\to
k^{\times}$ the character through which
$T_{{e^*}(.)}^{-1}$ acts on $k.e_j$,
i.e. $T^{-1}_{{e^*}(\gamma)}(e_j)=\eta_j(\gamma)e_j$ for $\gamma\in\Gamma$.

\begin{lem} \label{irr} (a) There are $0\le k_{e_j}\le q-1$ for $0\le j\le d$, not all of them $=0$ and not all of them $=q-1$, such that\begin{gather}t^{k_{e_j}}\varphi \otimes T^{-1}_{{\omega}}(e_{j})=-1\otimes e_j\label{gabriel}\end{gather}in $\Delta(M)$ for all $0\le j\le d$.

(b) If for any $1\le j\le d$ there is some $0\le i\le d$ with $k_{e_i}\ne k_{e_{i+j}}$, then $\Delta(M)$ is irreducible as a $k[[t]][\varphi]$-module.

(c) Suppose that for any $1\le j\le d$ which satisfies $k_{e_i}=
  k_{e_{i+j}}$ for all $0\le i\le m$ there is some $0\le i\le d$ with $\eta_{i}\ne\eta_{i+j}$. Then $\Delta(M)$ is irreducible as an ${\mathfrak O}$-module.\end{lem}

{\sc Proof:} For $M$ as above, $\nabla(M)$ is generated by elements of the form
$h(e)=t^{k_e}\varphi\otimes T_{\omega}^{-1}(e)+1\otimes e$. They give rise to formula (\ref{gabriel}), hence statement (a). For statements (b) and (c) apply Proposition \ref{altfund}; in (c) notice that $\gamma\cdot(1\otimes e_j)=\eta_j(\gamma)\otimes e_j$ for $\gamma\in\Gamma$.\hfill$\Box$\\

\begin{lem}\label{mafr} (a) Conjugating $\chi$ by powers of
  $T_{\omega}$ means cyclically permuting the ordered tuple
  $((\eta_0,k_{e_0}),\ldots,(\eta_d,k_{e_d}))$ associated with $\chi$ as
  above. Knowing the conjugacy class of $\chi$ (under powers of
  $T_{\omega}$) is equivalent with knowing the tuple
  $((\eta_0,k_{e_0}),\ldots,(\eta_d,k_{e_d}))$ up to cyclic permutations,
  together with $\chi(T^{d+1}_{\omega})$.

(b) (Vign\'{e}ras) Two standard supersingular ${\mathcal
  H}$-modules are isomorphic if and only if the element
  $T^{d+1}_{\omega}\in {\mathcal
  H}$ acts on them by the same constant in $k^{\times}$ and if they arise from two supersingular characters ${\mathcal
  H}_{\rm aff}\to k$ which are conjugate under some power of $T_{\omega}$. 

(c) (Vign\'{e}ras) A standard supersingular ${\mathcal
  H}$-module $M$ arising from $\chi$ is simple if and only if the orbit of $\chi$ under conjugation by powers of $T_{\omega}$ has cardinality $d+1$.
\end{lem}

{\sc Proof:} Statement (a) is clear. For (b) and (c) see \cite{vigneras} Proposition 3 and Theorem 5.\hfill$\Box$\\

\begin{pro}\label{ennobesuch} The functor $M\mapsto\Delta(M)$ induces a bijection between the set of isomorphim classes of standard supersingular ${\mathcal
  H}$-modules and the set of standard
cyclic objects in ${\rm Mod}^{\rm ad}({\mathfrak O})$ of $k$-dimension $d+1$. If the standard supersingular ${\mathcal
  H}$-module $M$ is simple, then $\Delta(M)\in{\rm Mod}^{\rm ad}({\mathfrak O})$ is simple.
\end{pro} 

{\sc Proof:} This follows from Lemma \ref{irr} and Lemma \ref{mafr}.\hfill$\Box$\\

\begin{satz}\label{molola1} (a) The functor $M\mapsto\Delta(M)^*\otimes_{k[[t]]}k((t))$ induces a bijection between the set of isomorphim classes of standard supersingular ${\mathcal
  H}$-modules and the set of isomorphim classes of standard cyclic \'{e}tale $(\varphi,\Gamma)$-modules of dimension $d+1$.

(b) The functor $M\mapsto\Delta(M)^*\otimes_{k[[t]]}k((t))$ induces a bijection between the set of isomorphim classes of simple supersingular ${\mathcal
  H}$-modules of $k$-dimension $d+1$ and the set of isomorphim classes of simple \'{e}tale $(\varphi,\Gamma)$-modules of dimension $d+1$.\end{satz}

{\sc Proof:} Statement (a) follows from Proposition \ref{responseday} and Proposition \ref{ennobesuch}. Statement (b) follows from statement (a) and the full faithfulness of the functor $M\mapsto\Delta(M)^*\otimes_{k[[t]]}k((t))$ on supersingular ${\mathcal
  H}$-modules, see Theorem \ref{7vorallhe} below. (To see that if $M$ is simple then so is $\Delta(M)^*\otimes_{k[[t]]}k((t))$ one may alternatively use Proposition \ref{ennobesuch} together with Proposition
\ref{irrnopsi}.)\hfill$\Box$\\

{\bf Remark:} We may rewrite
the equation (\ref{gabriel}) as \begin{align}t^{k_{e_j}}\varphi \otimes
  e_{j+1}&=-1\otimes e_j\quad\quad\mbox{ for }0\le j\le
  d-1\notag\\t^{k_{e_d}}\varphi \otimes \chi(T^{-d-1}_{\omega})e_{0}&=-1\otimes e_d\notag\end{align}where we used
$T_{{\omega}}^{-1}(e_d)=T_{{\omega}}^{-d-1}(e_0)=\chi(T^{-d-1}_{\omega})e_0$. Thus $(-1)^{d+1}\chi(T^{-d-1}_{\omega})\in k^{\times}$ is the constant referred to in Lemma \ref{ossa}.

\begin{kor}\label{molola} The functor $M\mapsto\Delta(M)^*\otimes_{k[[t]]}k((t))$, composed with the functor of Theorem \ref{sosego}, induces a bijection between the set of isomorphim classes of standard supersingular ${\mathcal
  H}$-modules of $k$-dimension $d+1$ and the set of isomorphim classes of $(d+1)$-dimensional standard cyclic ${\rm Gal}(\overline{F}/F)$-representations.
\end{kor}

{\sc Proof:} Theorem \ref{molola1}. \hfill$\Box$\\

{\bf Remark:} (a) Combining
Corollary \ref{molola} and Theorem \ref{molola1} one can derive the following "numerical Langlands correspondence": the set of (absolutely) simple $(d+1)$-dimensional ${\mathcal
  H}$-modules with fixed scalar action by $T_{\omega}^{d+1}$ has the same cardinality as the set of (absolutely) irreducible $(d+1)$-dimensional ${\rm Gal}(\overline{F}/F)$-representations with fixed determinant of Frobenius. This numerical Langlands correspondence was proven already in \cite{vigneras} Theorem 5.

(b) There is an alternative and arguably more natural definition of
supersingularity for ${\mathcal
  H}$-modules. Its agreement with the one given in subsection \ref{defhsu},
and hence the "numerical Langlands correspondence" with respect to
this alternative definition of supersingularity, was proven in \cite{ol}. 

\subsection{Reconstruction of an initial segment of $M$ from $\Delta(M)$}

\label{allnov}

Let ${[0,q-1]}^{\Phi}$ be the set of tuples ${\mu}=(\mu_{i})_{0\le
  i\le d}$ with $\mu_i\in\{0,\ldots, q-1\}$ and $\sum_{0\le i\le
  d}{\mu}_{i}\equiv 0$ modulo $(q-1)$. We often read the indices as
elements of ${\mathbb Z}/(d+1)$, thus $\mu_i=\mu_j$ for $i,
j\in{\mathbb Z}$ whenever $i-j\in(d+1){\mathbb Z}$.

Let $\Delta$ be an ${\mathfrak O}$-module. For
$\mu\in[0,q-1]^{\Phi}$ let ${\mathcal F}{\Delta}[t]^{\mu}$ be the $k$-subvector space of
$\Delta[t]=\{x\in \Delta\,|\,tx=0\}$ generated by all $x\in {\Delta}[t]$
satisfying $t^{\mu_{i}}\varphi\ldots t^{\mu_1}\varphi t^{\mu_0}\varphi x\in {\Delta}[t]$ for all $0\le i\le d$, as well as $t^{\mu_{d}}\varphi\ldots t^{\mu_1}\varphi t^{\mu_0}\varphi x\in k^{\times}x$.

Put ${\mathcal F}{\Delta}[t]=\sum_{\mu\in [0,q-1]^{\Phi}}{\mathcal F}{\Delta}[t]^{\mu}$ (sum in ${\Delta}[t]$).

\begin{lem}\label{auchwi} ${\mathcal F}{\Delta}[t]=\bigoplus_{\mu\in [0,q-1]^{\Phi}}{\mathcal F}{\Delta}[t]^{\mu}$, i.e. the sum is direct.
\end{lem}

{\sc Proof:} Consider the lexicographic enumeration $\mu(1),
\mu(2),\mu(3),\ldots $ of $[0,q-1]^{\Phi}$ such that
for each pair $r'>r$ there is some $0\le i_0\le d$ with $\mu_i(r)\ge \mu_i(r')$
for all $i<i_0$, and $\mu_{i_0}(r)> \mu_{i_0}(r')$. Let
$\sum_{r\ge1}x_{r}=0$ with $x_{r}\in
{\mathcal F}{\Delta}[t]^{\mu(r)}$. We prove $x_{r}=0$ for all $r$ by
induction on $r$. So, fix $r$ and assume $x_{r'}=0$ for all $r'<r$, hence $\sum_{r'\ge
  r}x_{r'}=\sum_{r\ge1}x_{r}-\sum_{r'<r}x_{r}=0$. For $r'>r$ we have $t^{\mu_d(r)}\varphi\cdots
t^{\mu_0(r)}\varphi(x_{r'})=0$. Therefore \begin{align}0=t^{\mu_d(r)}\varphi\cdots
t^{\mu_0(r)}\varphi (\sum_{r'\ge r}x_{r'})&=\sum_{r'\ge r}t^{\mu_d(r)}\varphi\cdots
t^{\mu_0(r)}\varphi x_{r'}\notag\\{}&=t^{\mu_d(r)}\varphi\cdots  t^{\mu_0(r)}\varphi x_{r}\in k^{\times}x_{r}\notag\end{align}and hence $x_r=0$.\hfill$\Box$\\

We define $k$-linear endomorphisms $T_{{\omega}}$, $T_s$ and $T_{e^*(\gamma)}$ (for $\gamma\in\Gamma$) of
${\mathcal F}{\Delta}[t]$ as follows. In view of Lemma \ref{auchwi} it is enough to
define their values on $x\in{\mathcal F}{\Delta}[t]^{\mu}$; we
put\begin{gather}T_{{\omega}}(x)=-t^{\mu_0}\varphi x,\quad\quad\quad\quad T_{e^*(\gamma)}(x)=\gamma^{-1}\cdot x,\notag\\T_s(x)=\left\{\begin{array}{l@{\quad:\quad}l}-x  &  \mu_d=q-1\notag\\0& \mu_d<q-1\end{array}\right.\notag.\end{gather}Here $\gamma^{-1}\cdot x$ is understood with respect to the $\Gamma$-action induced by the
  ${\mathfrak O}$-module structure on $\Delta(M)$.\\

{\bf Definition:} For an ${\mathcal
  H}^{\sharp\sharp}$-module $M$ and $\mu\in[0,q-1]^{\Phi}$ let ${\mathcal F}M^{\mu}$
denote the $k$-subvector space of $M$ consisting of $x\in M$ satisfying the following conditions for all $0\le i\le
d$:\begin{gather}T^{-1}_{\alpha_1^{\vee}(\gamma)}(T_{{\omega}}^i(x))=\gamma^{\mu_{i-1}}T_{{\omega}}^i(x)\quad\mbox{ for all }\gamma\in\Gamma,\label{nachsy}\\ T_s(T_{{\omega}}^i(x))=\left\{\begin{array}{l@{\quad:\quad}l}-T_{{\omega}}^i(x) &  \mu_{i-1}=q-1\\0& \mu_{i-1}<q-1\end{array}.\right.\label{dochnachsy}\end{gather}Let ${\mathcal F}M$ denote the subspace of $M$ generated by the ${\mathcal F}M^{\mu}$ for all $\mu\in[0,q-1]^{\Phi}$.

 For $\mu\in[0,q-1]^{\Phi}$ let $\epsilon_{\mu}\in{[0,q-2]}^{\Phi}$ be the unique element with\begin{gather}(\epsilon_{\mu})_{-i}\equiv \mu_{i}\,\, {\rm mod
}(q-1).\label{allsavi}\end{gather}for all $i$.

\begin{lem}\label{ennoad} (a) We have ${\mathcal F}M^{\mu}\subset M^{\epsilon_{\mu}}$.

(b) ${\mathcal F}M$ is an ${\mathcal
  H}^{\sharp\sharp}$-submodule of $M$. 

(c) ${\mathcal F}M$ contains each ${\mathcal
  H}^{\sharp\sharp}$-submodule of $M$ which is a subquotient of a standard supersingular ${\mathcal
  H}^{\sharp\sharp}$-module.

(d) Suppose that $M$ is supersingular. Viewing the isomorphism ${\Delta(M)}[t]\cong M$ (Proposition \ref{jishnusharp}) as an identity, we have ${\mathcal F}M^{\mu}\subset
{\mathcal F}{\Delta(M)}[t]^{\mu}$ for each $\mu\in[0,q-1]^{\Phi}$, and in particular\begin{gather}{\mathcal F}M\subset {\mathcal F}{\Delta(M)}[t]\label{fiedro}.\end{gather} The operators $T_{{\omega}}$, $T_s$ and $T_{e^*(\gamma)}$ acting on ${\mathcal F}{\Delta(M)}[t]$ as defined above restrict to the operators $T_{{\omega}}, T_s, T_{e^*(\gamma)}\in{\mathcal H}^{\sharp\sharp}$ acting on ${\mathcal F}M$.\end{lem} 

{\sc Proof:} (a) Let $\mu\in[0,q-1]^{\Phi}$. For $x\in {\mathcal F}M^{\mu}$, any $\gamma\in\Gamma$ and any $i$ we compute$$T_{\alpha_{1-i}^{\vee}(\gamma)}^{-1}(x)=T_{(\omega^{i}\cdot\alpha_{1}^{\vee})(\gamma)}^{-1}(x)=T_{\omega}^{-i}T_{\alpha_1^{\vee}(\gamma)}^{-1}T_{\omega}^{i}(x)=\gamma^{{\mu}_{i-1}}x=\gamma^{(\epsilon_{\mu})_{1-i}}x,$$i.e. $x\in M^{\epsilon_{\mu}}$.

(b) Let $\mu\in[0,q-1]^{\Phi}$ and define $\mu'\in[0,q-1]^{\Phi}$ by
  $\mu'_i=\mu_{i+1}$ for all $i$. For $x\in {\mathcal F}M^{\mu}$, any
  $\gamma\in\Gamma$ and any $i$ we
  compute$$T^{-1}_{\alpha_1^{\vee}(\gamma)}(T_{{\omega}}^i(T_{{\omega}}(x)))=T^{-1}_{\alpha_1^{\vee}(\gamma)}(T_{{\omega}}^{i+1}(x))=\gamma^{\mu_{i}}T_{{\omega}}^{i+1}(x)=\gamma^{\mu_{i}}
  T_{{\omega}}^i(T_{{\omega}}(x)).$$We also find
  $T_s(T_{{\omega}}^i(T_{{\omega}}(x)))=T_s(T_{{\omega}}^{i+1}(x))=-T_{{\omega}}^{i+1}(x)=-T_{{\omega}}^i(T_{{\omega}}(x))$
  if $\mu_i=q-1$, but
  $T_s(T_{{\omega}}^i(T_{{\omega}}(x)))=T_s(T_{{\omega}}^{i+1}(x))=0$ if
  $\mu_i<q-1$. Together this shows $T_{{\omega}}(x)\in{\mathcal F}M^{\mu'}$,
  i.e. $T_{\omega}({\mathcal F}M^{\mu})\subset{\mathcal F}M^{\mu'}$. It is
  immediate from the definitions that $T_{s}({\mathcal
    F}M^{\mu})\subset{\mathcal F}M^{\mu}$. For $x\in {\mathcal F}M^{\mu}$, any
  $\gamma,\gamma'\in\Gamma$ and any $i$ we
  compute$$T_{\alpha_{1}^{\vee}(\gamma)}^{-1}
  T_{\omega}^{i}(T_{e^*(\gamma')}(x))=T_{\alpha_{1}^{\vee}(\gamma)}^{-1}T_{(\omega^{-i}\cdot
    e^*)(\gamma')}T_{\omega}^{i}(x)=T_{(\omega^{-i}\cdot
    e^*)(\gamma')}\gamma^{\mu_{i-1}}T_{\omega}^{i}(x)=\gamma^{\mu_{i-1}}T_{\omega}^{i}(T_{e^*(\gamma')}(x)).$$If
  $\mu_{i-1}=q-1$ we also compute$$T_sT_{{\omega}}^i(T_{e^*(\gamma')}(x))=T_{(s\cdot\omega^{-i}\cdot
    e^*)(\gamma')}T_sT_{\omega}^{i}(x)=-T_{(s\cdot\omega^{-i}\cdot
    e^*)(\gamma')}T_{\omega}^{i}(x)$$$$=-T_{\omega}^i(T_{(\omega^{i}\cdot
    s\cdot\omega^{-i}\cdot e^*)(\gamma')}(x))=-T_{\omega}^i(T_{
    e^*(\gamma')}(x)).$$Here, in the last equation we use
  $\omega^{i}\cdot
  s\cdot\omega^{-i}\cdot e^*=e^*$ for $2\le i\le d$; for $i=1$ we use $\omega\cdot
  s\cdot\omega^{-1}\cdot e^*-e^*=\alpha_0^{\vee}$ and
  $T_{\alpha_0^{\vee}(\gamma')}(x)=T_{\omega}^{-1}T_{\alpha_1^{\vee}(\gamma')}T_{\omega}(x)=\gamma^{-\mu_0}x=x$
  (as $\mu_{0}=q-1$); for $i=0$ we use $s\cdot e^*-e^*=-\alpha_1^{\vee}$ and
  $T_{-\alpha_1^{\vee}(\gamma')}(x)=\gamma^{\mu_{-1}}x=x$
  (as $\mu_{-1}=q-1$). If however $\mu_{i-1}<q-1$ then $T_sT_{{\omega}}^i(T_{e^*(\gamma')}(x))=T_{(s\cdot\omega^{-i}\cdot
    e^*)(\gamma')}T_sT_{\omega}^{i}(x)=0$. Together this shows $T_{e^*(\gamma')}(x)\in{\mathcal F}M^{\mu}$, i.e. $T_{e^*(\gamma')}({\mathcal F}M^{\mu})\subset{\mathcal F}M^{\mu}$.

  (c) On a standard supersingular ${\mathcal
  H}^{\sharp\sharp}$-module, and hence on its subquotients, the
  actions of $T_{\omega}$, $T_s$ and $T_{\alpha_1^{\vee}(\gamma)}$ satisfy
  formulae (\ref{nachsy}) and (\ref{dochnachsy}), for suitable $\mu$'s.

  (d) Let $\mu\in[0,q-1]^{\Phi}$ and define $\mu'\in[0,q-1]^{\Phi}$ by
  $\mu'_i=\mu_{i+1}$ for all $i$. Let $x\in{\mathcal F}M^{\mu}$. The proof of
  (b) shows $T_{\omega}(x)\in\sum_aM^{\epsilon_{\mu'}}_{\underline{a}}[0]$ if
  $\mu'_d=\mu_0<q-1$,
  resp. $T_{\omega}(x)\in\sum_aM^{\epsilon_{\mu'}}_{\underline{a}}[1]$ if
  $\mu'_d=\mu_0=q-1$. In either case, the definition of $\Delta(M)$ then says
  $T_{{\omega}}(x)=-t^{\mu_0}\varphi x$. This shows ${\mathcal F}M^{\mu}\subset
{\mathcal F}{\Delta(M)}[t]^{\mu}$ and that the action of $T_{{\omega}}$ on
${\mathcal F}M$ is indeed as stated. For the actions of $T_s$ and
$T_{e^*(\gamma)}$ this is clear anyway.\hfill$\Box$\\

{\bf Remark:} The inclusion
(\ref{fiedro}) is in fact an equality.

\subsection{Reconstruction of $\sharp$-supersingular ${\mathcal
  H}^{\sharp}$-modules $M$ from $\Delta(M)$}

\begin{lem}\label{4vorstubai} Let $M$ be an irreducible supersingular ${\mathcal
  H}$-module. Let $\mu\in[0,q-1]^{\Phi}$, let $x\in M$ and $u_{i,c}\in
  M^{\omega^{-1} s\omega^{i+1}\cdot\epsilon_{\mu}}$ for $i\ge0$ and
  $0\le c\le q-2$ (with $\epsilon_{\mu}$ given by formula (\ref{allsavi})). Assume $u_{i,c}=0$ if 

(i) $\mu_i=0$, or 

(ii) $\mu_i=q-1$ and $c>0$, or 

(iii) $\mu_{i}<q-1$ and $c\ge q-1-\mu_{i}$. 

Assume that, if we put $x\{-1\}=x$, then$$x\{i\}=t^{\mu_i}\varphi(x\{i-1\})-\sum_{c=0}^{q-2}t^c\varphi
u_{i,c}$$belongs to $M\cong \Delta(M)[t]$ for each $i\ge0$. Finally, assume that
$x\{D\}=x$ for some $D>0$ with $D+1\in{\mathbb Z}(d+1)$. Then there is some ${x}'\in M$ with $x-x'\in
M^{\epsilon_{\mu}}$ and such that$$x'\{i\}=t^{\mu_i}\varphi(\ldots
(t^{\mu_1}\varphi(t^{\mu_0}\varphi x'))\ldots )$$belongs to $M$ for each $i$,
and $x'\{D\}=x'$. Moreover, if $x$ is an eigenvector
for $T_{e^*(\Gamma)}$, then $x'$ can be chosen to be an eigenvector
for $T_{e^*(\Gamma)}$, with the same eigenvalues.
\end{lem}

{\sc Proof:} It is easy to see that all the irreducible subquotients of a standard supersingular ${\mathcal
  H}$-module are isomorphic. In particular, an irreducible supersingular ${\mathcal
  H}$-module is isomorphic with a submodule of a standard supersingular ${\mathcal
  H}$-module. Therefore we may assume that $M$ itself is a (not necessarily
irreducible) standard supersingular ${\mathcal
  H}$-module. We then have a direct sum decomposition
$M=\oplus_{j=0}^dM^{[j]}$ with ${\rm dim}_k(M^{[j]})=1$ and integers $0\le
k_j\le q-1$ such that \begin{gather}T_{\omega}(M^{[j+1]})=t^{k_j}\varphi
  (M^{[j+1]})=M^{[j]}\label{3vorstubai}\end{gather}(always reading $j$ modulo
$(d+1)$). More precisely, we have $M^{[j]}\subset M^{\epsilon_j}$ for certain $\epsilon_j\in[0,q-2]^{\Phi}$, and chosing the above $k_j$ minimally, we have $k_j\equiv (\omega\cdot\epsilon_{j+1})_1$ modulo $(q-1)$. It follows that\begin{gather}k[t]\varphi M=\bigoplus_{j=0}^dk[t]\varphi M^{[j]}=\bigoplus_{j=0}^d\bigoplus_{c=0}^{k_j}t^c\varphi M^{[j+1]}.\label{3nachstubai}\end{gather}For $m\in
M$ write $m=\sum_jm^{[j]}$ with $m^{[j]}\in M^{[j]}$. By formulae
(\ref{3vorstubai}), (\ref{3nachstubai}), the defining formula for
$x\{i\}$ splits up into the formulae\begin{gather}x\{i\}^{[j]}=t^{\mu_i}\varphi(x\{i-1\}^{[j+1]})-\sum_{c=0}^{q-2}t^c\varphi(u_{i,c}^{[j+1]})\label{3vorstubai1}\end{gather}for
all $j$. We use them to show\begin{gather}t^c\varphi(u_{i,c}^{[j+1]})= 0\quad \mbox{
    if }c-\mu_i\notin(q-1){\mathbb Z}.\label{wienneu}\end{gather}If
$\mu_i\in\{0,q-1\}$ then formula (\ref{wienneu}) follows from our
assumptions on the $u_{i,c}$'s. Now assume $\mu_i\notin\{0,q-1\}$ and $u_{i,c}^{[j+1]}\ne0$ for some $c$. The assumption $u_{i,c}\in
M^{\omega^{-1} s\omega^{i+1}\cdot\epsilon_{\mu}}$ implies $T_{\omega}(u_{i,c}^{[j+1]})\in M^{s\omega^{i+1}\cdot\epsilon_{\mu}}$, and since$$q-1-\mu_i=q-1-\epsilon_{-i}=(s\omega^{i+1}\cdot\epsilon_{\mu})_1\quad\mbox{ if
}\mu_i\notin\{0,q-1\}$$we get $T_{\omega}(u_{i,c}^{[j+1]})=-t^{q-1-\mu_i}\varphi(u_{i,c}^{[j+1]})$, i.e. $k_j=q-1-\mu_i$. Now $\sum_{c=0}^{k_j}t^c\varphi M^{[j+1]}$ is a {\it direct} sum of one dimensional $k$-vector spaces, with $x\{i\}^{[j]}\in t^{k_j}\varphi M^{[j+1]}$, with $t^{\mu_i}\varphi(x\{i-1\}^{[j+1]})\in t^{\mu_i}\varphi M^{[j+1]}$ and with $t^c\varphi(u_{i,c}^{[j+1]})\in t^{c}\varphi M^{[j+1]}$ for all $c$. Since by assumption $u_{i,c}=0$ for $c\ge
q-1-\mu_{i}=k_j$, formula (\ref{3vorstubai1}) shows $t^c\varphi(u_{i,c}^{[j+1]})= 0$ whenever $c\ne\mu_i$.

 Formula (\ref{wienneu}) is proven. Arguing once more with formulae
(\ref{3vorstubai}), (\ref{3nachstubai}) and (\ref{3vorstubai1}) shows
\begin{gather}[t^{\mu_i}\varphi(x\{i-1\}^{[j+1]})=0\quad\mbox{ or
  }\quad\varphi(u_{i,0}^{[j+1]})=0]\quad\quad\quad\quad\mbox{ if
  }\mu_i=q-1.\label{formelcrux}\end{gather}In the following, by $u_{i,q-1}$ we
mean $u_{i,0}$. If $t^{\mu_i}\varphi(u_{i,\mu_i}^{[j+1]})\ne 0$ we may write
$$t^{\mu_i}\varphi(x\{i-1\}^{[j+1]})-t^{\mu_i}\varphi(u_{i,\mu_i}^{[j+1]})=\rho_{i,j}t^{\mu_i}\varphi(u_{i,\mu_i}^{[j+1]})$$
for some $\rho_{i,j}\in k$, since $t^{\mu_i}\varphi(x\{i-1\}^{[j+1]})$ and
$t^{\mu_i}\varphi(u_{i,\mu_i}^{[j+1]})$ belong to the same one-dimensional
$k$-vector space. The upshot of formulae (\ref{wienneu})
and (\ref{formelcrux}) is then that formula (\ref{3vorstubai1}) simplifies to
become either \begin{gather}x\{i\}^{[j]}=t^{\mu_i}\varphi(x\{i-1\}^{[j+1]})\label{pridon1}\end{gather}or\begin{gather}x\{i\}^{[j]}=\rho_{i,j}t^{\mu_i}\varphi(u_{i,\mu_i}^{[j+1]})\label{pridon2}\end{gather}
for some $\rho_{i,j}\in k$. Departing from $x^{[j]}=x\{D\}^{[j]}$ we repeatedly substitute formula (\ref{pridon1}); if this is possible $D+1$ many times we end up with$$x^{[j]}=x\{D\}^{[j]}=t^{\mu_D}\varphi(\ldots
(t^{\mu_1}\varphi(t^{\mu_0}\varphi (x^{[j]})))\ldots ),$$and in this case we put
$n(j)=0$. Otherwise, after $D+1-n(j)$ many substitutions of formula (\ref{pridon1}), for some $1\le n(j)\le D+1$, we end the procedure by substituting formula (\ref{pridon2}) (once) and obtain$$x^{[j]}=x\{D\}^{[j]}=\rho_{j}t^{\mu_D}\varphi(\ldots
(t^{\mu_{n(j)}}\varphi(t^{\mu_{n(j)-1}}\varphi(u_{n(j)-1,\mu_{n(j)-1}}^{[j+1-n(j)]})))\ldots)$$
with $t^{\mu_{n(j)-1}}\varphi u_{n(j)-1,\mu_{n(j)-1}}^{[j+1-n(j)]}\ne0$, for some $\rho_j\in k$. 

We study this second case $n(j)>0$ further. By construction,$$
  w_j\{-1\}=t^{\mu_{n(j)-1}}\varphi(u_{n(j)-1,\mu_{n(j)-1}}^{[j+1-n(j)]})$$is
  non-zero and belongs
  to $M$. On the other hand, $u_{n(j)-1,\mu_{n(j)-1}}\in M^{\omega^{-1} s
    \omega^{n(j)}\cdot\epsilon_{\mu}}$ implies $T_{\omega}(u^{[j+1-n(j)]}_{n(j)-1,\mu_{n(j)-1}})\in M^{s
    \omega^{n(j)}\cdot\epsilon_{\mu}}$ and hence $$t^{(s
    \omega^{n(j)}\cdot\epsilon_{\mu})_1}\varphi(u_{n(j)-1,\mu_{n(j)-1}}^{[j+1-n(j)]})=-T_{\omega}(u^{[j+1-n(j)]}_{n(j)-1,\mu_{n(j)-1}})\in M^{s
    \omega^{n(j)}\cdot\epsilon_{\mu}}.$$Together this means $\mu_{n(j)-1}\equiv(s
    \omega^{n(j)}\cdot\epsilon_{\mu})_1$ modulo $(q-1)$ and $w_j\{-1\}\in
  M^{s
    \omega^{n(j)}\cdot\epsilon_{\mu}}$. But we also have $\mu_{n(j)-1}\equiv(\omega^{n(j)}\cdot\epsilon_{\mu})_1$. Combining we see
    $\mu_{n(j)-1}\equiv-\mu_{n(j)-1}$ modulo $(q-1)$. Hence, we
    either have
    $\mu_{n(j)-1}=0$ or $\mu_{n(j)-1}=\frac{q-1}{2}$ or $\mu_{n(j)-1}=q-1$. In
    view of the
    assumed vanishings of the $u_{i,c}$'s (and of
    $u_{n(j)-1,\mu_{n(j)-1}}^{[j+1-n(j)]}\ne0$) this leaves $\mu_{n(j)-1}=q-1$
    as the only possibility. It follows that$$s
    \omega^{n(j)}\cdot\epsilon_{\mu}=\omega^{n(j)}\cdot\epsilon_{\mu}$$and hence $w_j\{-1\}\in M^{
      \omega^{n(j)}\cdot\epsilon_{\mu}}$. Next, again by construction we know that$$ w_j\{s\}=t^{\mu_{n(j)+s}}\varphi (
 w_j\{s-1\})$$belongs to $M$, for $0\le
  s\le D-n(j)$. By what we learned about $w_j\{-1\}$ this implies $w_j\{s\}=(-1)^{s+1}T_{\omega}^{s+1}w_j\{-1\}\in M^{ \omega^{n(j)+s+1}\cdot\epsilon_{\mu}}$ by an induction on $s$ (and we also see $\mu_{n(j)+s}\in\{k_0,\ldots, k_d\}$ with the $k_{\ell}$ from formula (\ref{3vorstubai})). For $s=D-n(j)$ we obtain $x^{[j]}=x\{D\}^{[j]}\in
  M^{\epsilon_{\mu}}$.

We now put $x'=\sum_{n(j)=0}x^{[j]}$.\hfill$\Box$\\

\begin{lem}\label{fabei} Let $M$ be an irreducible supersingular ${\mathcal
  H}$-module. Let $\mu\in[0,q-1]^{\Phi}$ and $x\in M$ such
  that$$x\{i\}=t^{\mu_i}\varphi(\ldots (t^{\mu_1}\varphi(t^{\mu_0}\varphi x
  ))\ldots)$$belongs to $M\cong \Delta(M)[t]$ for each $i\ge0$, and such
  that $x\{D\}=x$ for some $D>0$ with $D+1\in{\mathbb Z}(d+1)$. Then $x\in
  M^{\epsilon_{\mu}}$ and $x\{i\}=(-T_{\omega})^ix$ for each $i$.
\end{lem}

{\sc Proof:} This follows from the formulae (\ref{3vorstubai}) and
(\ref{3nachstubai}) in the proof of Lemma \ref{4vorstubai}. The argument is very similar to the one given in the proof of Lemma \ref{auchwi}.\hfill$\Box$\\

\begin{satz}\label{shmichaelmas} Let $M$ be a $\sharp$-supersingular ${\mathcal
    H}^{\sharp}$-module. Via the isomorphism $M\cong \Delta(M)[t]$, the action of ${\mathcal H}^{\sharp}$ on $M$ can be recovered from the action of ${\mathfrak O}$ on $\Delta(M)$.
\end{satz}

{\sc Proof:} We may assume ${\rm dim}_k(M)<\infty$. Define inductively the filtration $(F^iM)_{i\ge0}$ of $M$ by ${\mathcal
  H}^{\sharp}$-submodules as follows: $F^0M=0$, and $F^{i+1}M$ is the preimage of
${\mathcal F}(M/F^{i}M)$ under the projection $M\to M/F^{i}M$. The ${\mathcal
  H}^{\sharp}$-action on the graded pieces can be recovered in view of Lemma
\ref{ennoad}. Exhausting $M$ step by step it is therefore enough to consider
the following setting: The action of ${\mathcal H}^{\sharp}$ has already been recovered
on an ${\mathcal H}^{\sharp}$-submodule $M_0$ of $M$ and on the quotient $M/M_0$, and
the latter is irreducible. 

We reconstruct the action of $T_{e^*(\Gamma)}$ on $M$ by means
of $$T_{e^*(\gamma)}(x)=\gamma^{-1}\cdot x\quad \mbox{ for }\gamma\in
\Gamma$$as is tautological from our definitions. Next we are going to reconstruct the
decomposition\begin{gather}M=\bigoplus_{\epsilon\in[0,q-2]^{\Phi},\atop a\in
  [0,q-2]}M^{\epsilon}_{\underline{a}}.\label{shletzthuerd}\end{gather}Let
$D>0$ be such that $D+1\in{\mathbb Z}(d+1)$ and $f^{D+1}={\rm id}$ for each
$k$-vector space automorphism $f$ of $M$. (Such a $D$ does exist. Indeed, $M$
is finite, hence ${\rm Aut}_k(M)$ is finite, hence there is some $n\in{\mathbb
  N}$ with $f^n={\rm id}$ for each $f\in {\rm Aut}_k(M)$. Now take $D=(d+1)n-1$.) For
$\epsilon\in[0,q-2]^{\Phi}$ and $a\in [0,q-2]$ define $M^{[{\epsilon}]}_{\underline{a}}$ to be the
$k$-subspace of $M$ generated by all $x\in M$ with $\gamma\cdot x=\gamma^ax$
(all $\gamma\in\Gamma$) and
satisfying the following condition: There is some $\mu\in [0,q-1]^{\Phi}$
(depending on $x$) with $\epsilon_{\mu}=\epsilon$, and there are
$u_{i,c}\in M_0^{\omega^{-1}s\omega^{i+1}\cdot\epsilon}$ for
$i\ge0$ and $0\le c\le q-2$ with the following properties: Firstly,
$u_{i,c}=0$ if

(i) $\mu_i=0$, or 

(ii) $\mu_i=q-1$ and $c>0$, or 

(iii) $\mu_{i}<q-1$ and $c\ge q-1-\mu_{i}$. 

Secondly, putting $x\{-1\}=x$ and\begin{gather}x\{i\}=t^{\mu_i}\varphi(x\{i-1\})-\sum_{c}t^c\varphi
  u_{i,c},\label{shsabzwgo}\end{gather}we have $x\{i\}\in M\cong \Delta(M)[t]$ for any
    $i$, as well as $x\{D\}=x$. 

It will be enough to prove $M^{\epsilon}_{\underline{a}}=M^{[{\epsilon}]}_{\underline{a}}$. We first show\begin{gather}M^{\epsilon}_{\underline{a}}\subset
    M^{[{\epsilon}]}_{\underline{a}}.\label{shcla1}\end{gather}

We start with $\overline{x}\in {\mathcal F}(M/M_0)^{\mu}\cap (M/M_0)^{\epsilon}_{\underline{a}}$ for some $\mu$ with
$\epsilon_{\mu}=\epsilon$. By
Proposition \ref{novordschutzengel} we may lift it to some $x\in M^{\epsilon}$ such that for each $i$ with
  $T_sT_{\omega}^{i+1}\overline{x}=0$ we have
$T_s^2T_{\omega}^{i+1}x=0$. As $T_{\omega}$ maps simultaneous eigenspaces for the $T_{t}$ (with $t\in\overline{T}$) again to such simultaneous eigenspaces, and as $T_s^2$ commutes with the $T_t$, we may assume $x\in M^{\epsilon}_{\underline{a}}$. Putting
$$x\{i\}=(-T_{\omega})^{i+1}x$$for $-1\le i\le D$, repeated application of Lemma \ref{samserswo} shows $x\{i\}\in
M^{\omega^{i+1}\cdot\epsilon}_{\underline{a_{\epsilon,i}}}$ with$$a_{\epsilon,-1}=a,\quad 
a_{\epsilon,0}=a-\epsilon_0\quad\mbox{ and }\quad
a_{\epsilon,i}=a-\epsilon_0-\epsilon_{d-i+1}-\ldots-\epsilon_d$$for $i\le d$,
and then $a_{\epsilon,i}=a_{\epsilon,i'}$ for $i-i'\in{\mathbb Z}(d+1)$.

If $0\le\mu_i<q-1$ put
$$u_{i,c}=T_{\omega}^{-1}((T_s(x\{i\}))_{\underline{c+\mu_i+a_{\epsilon,i}}}).$$As $\overline{x}\in {\mathcal F}(M/M_0)^{\mu}$ and $\mu_i<q-1$ we have $u_{i,c}\in M_0$, and as $x\{i\}\in
M^{\omega^{i+1}\cdot\epsilon}$ we have $u_{i,c}\in
M^{\omega^{-1}s\omega^{i+1}\cdot\epsilon}$, together $u_{i,c}\in
M_0^{\omega^{-1}s\omega^{i+1}\cdot\epsilon}$. From $\mu_i<q-1$ we furthermore
deduce $k_{x\{i\}}=(\omega^{i+1}\cdot\epsilon)_{1}= \mu_i$, and since
$T_s^2x\{i\}=0$ we then see\begin{gather}t^{\mu_i}\varphi(x\{i-1\})-x\{i\}-\sum_{c}t^c\varphi
u_{i,c}=h(-x\{i\})=0.\label{shgruendo}\end{gather}Since furthermore
$(T_s(x\{i\}))_{\underline{c+\mu_i+a_{\epsilon,i}}}=0$ and hence $u_{i,c}=0$ for
$q-1-\mu_{i}\le c\le q-2$ by $\sharp$-supersingularity (if $0<\mu_i<q-1$ then $\mu_{i}=\epsilon_{-i}$),
all the conditions on the $u_{i,c}$'s in the definition of
$x\in M^{[{\epsilon}]}_{\underline{a}}$ are satisfied.

If $\mu_i=q-1$ we have $T_s^2(T_s^2x\{i\})=T_s^2x\{i\}$ and hence $k_{T_s^2x\{i\}}=q-1$ (independently of the value of $\mu_i$ we have $(\omega^{i+1}\cdot\epsilon)_{1}\equiv \mu_i$ modulo
$(q-1)$), hence\begin{gather}t^{q-1}\varphi T_{\omega}^{-1}(T_s^2x\{i\})+T_s^2x\{i\}=h(T_s^2x\{i\})=0.\label{papsthaeresie1}\end{gather}Similarly we see $k_{(x\{i\}-T_s^2x\{i\})}=0$ and hence\begin{gather}\varphi
T_{\omega}^{-1}(x\{i\}
-T_s^2x\{i\})+x\{i\}-T_s^2x\{i\}=h(x\{i\}-T_s^2x\{i\})=0.\label{papsthaeresie2}\end{gather}We
compute\begin{align} t^{q-1}\varphi(x\{i-1\})&=-t^{q-1}\varphi
T_{\omega}^{-1}(x\{i\})\notag\\{}&=-t^{q-1}\varphi
T_{\omega}^{-1}T_s^2(x\{i\})\notag\\{}&=T_s^2(x\{i\})\notag\\{}&=\varphi
T_{\omega}^{-1}(x\{i\}-T_s^2x\{i\})+x\{i\}\notag\end{align}where the second
equality is the result of applying $t^{q-1}$ to formula
(\ref{papsthaeresie2}), where the third equality is formula
(\ref{papsthaeresie1}) and where the fourth equality is formula
(\ref{papsthaeresie2}). Thus, putting
$u_{i,0}=T_{\omega}^{-1}(x\{i\}-T_s^2x\{i\})$ and $u_{i,c}=0$ for $c>0$, we
again get formula
(\ref{shgruendo}). Moreover, $u_{i,0}$ belongs to
$M_0$ as $\overline{x}\in {\mathcal F}(M/M_0)^{\mu}$ and $\mu_i=q-1$; but it
also belongs to $M^{\omega^{-1}s\omega^{i+1}\cdot\epsilon}$
since $\mu_i=q-1$ implies
$\omega^{-1}s\omega^{i+1}\cdot\epsilon=\omega^{i}\cdot\epsilon$. By
construction, $x\{d\}=(-T_{\omega})^{d+1}(x)$, hence $x\{D\}=(-T_{\omega})^{D+1}x=x$.

It follows that $x\in M^{[{\epsilon}]}_{\underline{a}}$. We have shown that
any element in ${\mathcal F}(M/M_0)^{\mu}\cap (M/M_0)^{\epsilon}_{\underline{a}}$, for $\mu$ with
$\epsilon_{\mu}=\epsilon$, lifts to an element in
$M^{{\epsilon}}_{\underline{a}}\cap M^{[{\epsilon}]}_{\underline{a}}$. Since we have $(M/M_0)^{\epsilon}=\sum_{\mu\in
  [0,q-1]^{\Phi}\atop\epsilon_{\mu}=\epsilon}{\mathcal F}(M/M_0)^{\mu}$
(see Lemma \ref{ennoad}), and since this is respected by the action of
$T_{e^*(\Gamma)}$, we thus have reduced our problem to showing $(M_0)^{\epsilon}_{\underline{a}}\subset
    M^{[{\epsilon}]}_{\underline{a}}$. But for this we may appeal to an induction on ${\rm
  dim}_k(M)$ (which we may assume to be finite). 

We have shown formula (\ref{shcla1}). Now we show\begin{gather}M^{[{\epsilon}]}_{\underline{a}}\subset
    M^{\epsilon}_{\underline{a}}.\label{shcla2}\end{gather}Let $x\in M^{[{\epsilon}]}_{\underline{a}}$,
$\mu\in [0,q-1]^{\Phi}$ (with $\epsilon_{\mu}=\epsilon$) and $u_{i,c}$ be as in the definition of
$M^{[{\epsilon}]}_{\underline{a}}$. Define $x\{i\}$ for $-1\le i\le D$ as in that definition. By
Lemma \ref{fabei} and the proof of the inclusion (\ref{shcla1}) we find
$\tilde{x}\in M^{\epsilon}_{\underline{a}}$ and $\tilde{u}_{i,c}\in M_0^{\omega^{-1}s\omega^{i+1}\cdot\epsilon}$ for $0\le
i\le D$ such that, after replacing $x$ by $x-\tilde{x}$ and $u_{i,c}$ by
$u_{i,c}-\tilde{u}_{i,c}$, we may assume $x\in M_0$.

{\it Claim: If $x\in M_0$ and if $M_0$ is irreducible, then there is some
  ${x}'\in (M_0)_{\underline{a}}$ with $x-{x}'\in
  (M_0)_{\underline{a}}^{\epsilon}$ and such that\begin{gather}x'\{i\}=t^{\mu_i}\varphi(\ldots
(t^{\mu_1}\varphi(t^{\mu_0}\varphi x'))\ldots )\notag\end{gather}belongs to $M_0$ for all $i$, and $x'\{D\}=x'$.}

This follows from Lemma \ref{4vorstubai}.

If $M_0$ is not irreducible, choose an ${\mathcal
  H}$-submodule $M_{00}$ in $M_0$ such that $M_0/M_{00}$ is irreducible. By the above claim and again invoking the proof of the inclusion (\ref{shcla1}), after
modifying $x$ by another element of $M^{\epsilon}_{\underline{a}}$ (now even
of $(M_0)^{\epsilon}_{\underline{a}}$) and suitably modifying the $u_{i,c}$, we may assume $u_{i,c}\in
M_{00}$. Thus,
it is now enough to solve the problem for the new $x\in (M_0)_{\underline{a}}$ (and the new $u_{i,c}\in
M_{00}$). We continue in this way. Since we may assume that ${\dim}_k(M)$
is finite, an induction on the dimension of $M$ allows us to conclude. 

We have reconstructed the decomposition (\ref{shletzthuerd}) of $M$. 

Now we reconstruct $T_sT_{\omega}$ acting on $M$. As we already know the decomposition (\ref{shletzthuerd}), it is enough to reconstruct
$T_sT_{\omega}(e)$ for $e\in M^{\epsilon'}_{\underline{a'}}$, all $\epsilon'$, $a'$. Given
such $e$, let $\overline{e}$ be its class in $M/M_0$. By Lemma \ref{samserswo}
there are then $\epsilon$, $a$ such that $T_{\omega}\overline{e}\in (M/M_0)^{\epsilon}_{\underline{a}}$. 

First assume $\epsilon_1=0$. We then reconstruct $T_sT_{\omega}(e)$ as
$T_sT_{\omega}(e)=t^{q-1}\varphi(e)$. Indeed, to see this we may assume (by Lemma \ref{quadrat}) that
$T_{\omega}(e)$ is an eigenvector for $T_s^2$. If
$T_s^2T_{\omega}(e)=T_{\omega}(e)$ and hence $T_sT_{\omega}(e)=-T_{\omega}(e)$, the claim follows from the definition of
$h(T_{\omega}(e))$. If $T_s^2T_{\omega}(e)=0$ then in fact $T_sT_{\omega}(e)=0$
(since also $\epsilon_1=0$), and the definition of
$h(T_{\omega}(e))$ shows $t^{q-1}\varphi(e)=0$.

Now assume $\epsilon_1>0$. This implies $T_s^2T_{\omega}(e)=0$ and
$k_{T_{\omega}(e)}=\epsilon_1$, and by $\sharp$-supersingularity we
get$$t^{k_{T_{\omega}(e)}+1}\varphi e=-\sum_{0\le c<q-1-k_{T_{\omega}(e)}}t^{c+1}\varphi
T_{\omega}^{-1}((T_sT_{\omega}e)_{\underline{c+\epsilon_1+a}}).$$Here
$(T_sT_{\omega}e)_{\underline{c+\epsilon_1+a}}\in M_0^{s\cdot\epsilon}$ and
$q-1-k_{T_{\omega}(e)}=(s\cdot\epsilon)_1$. The
map $$\bigoplus_{0\le c<q-1-k_{T_{\omega}(e)}}M_0^{s\cdot\epsilon}\longrightarrow M_0,\quad (y_c)_c\mapsto\sum_{0\le c<q-1-k_{T_{\omega}(e)}}t^{c+1}\varphi
T_{\omega}^{-1}(y_c)$$is injective. This is first seen in the case where $M_0$
is irreducible; it then follows by an obvious devissage argument. We therefore
see that
the $(T_sT_{\omega}e)_{\underline{c+\epsilon_1+a}}$ for $0\le c<q-1-k_{T_{\omega}(e)}$ can be read off from $t^{k_{T_{\omega}(e)}+1}\varphi e$,
hence also $T_sT_{\omega}e$ can be read off from $t^{k_{T_{\omega}(e)}+1}\varphi e$ (by $\sharp$-supersingularity).

The restriction of $T_{\omega}$ to $\{x\in M\,|\,T_sT_{\omega}(x)\in M_0\}$ is
reconstructed as follows. Given
$\overline{x}\in (M/M_0)^{\omega^{-1}\cdot\epsilon}_{\underline{a-\epsilon_1}}$ (for some
$\epsilon$, some $a$) with $T_sT_{\omega}\overline{x}=0$, we use the decomposition (\ref{18okt}) to lift
$\overline{x}$ to some $x\in M^{\omega^{-1}\cdot\epsilon}_{\underline{a-\epsilon_1}}$. Since $(\omega^{-1}\cdot\epsilon)_0=\epsilon_1$, Lemma
\ref{samserswo} says $T_{\omega}{x}\in M^{\epsilon}_{\underline{a}}$. It then follows from the
definitions that$$T_{\omega}x=-t^{\epsilon_1}\varphi
x-\sum_{c\ge0}t^c\varphi T_{\omega}^{-1}((T_sT_{\omega}x)_{\underline{c+\epsilon_1+a}}).$$We have now collected all the data required in Proposition \ref{schutzengel} for reconstructing $M$ as an ${\mathcal H}^{\sharp}$-module.\hfill$\Box$\\

\subsection{Full faithfulness on $\sharp$-supersingular ${\mathcal
  H}^{\sharp}$-modules}

\label{fullfafu}

Let ${\rm Rep}({\rm Gal}(\overline{F}/F))$
denote the category of representations of ${\rm
  Gal}(\overline{F}/F)$ on $k$-vector spaces which are projective limits of finite
dimensional continuous ${\rm Gal}(\overline{F}/F)$-representations.

Let ${\rm Mod}_{ss}({\mathcal H}^{\sharp})$ denote the category of
$\sharp$-supersingular ${\mathcal H}^{\sharp}$-modules. Let ${\rm
  Mod}_{ss}({\mathcal H})$, resp. ${\rm Mod}_{ss}({\mathcal
  H}^{\sharp\sharp})$, denote the category of supersingular ${\mathcal
  H}$-modules and supersingular ${\mathcal H}^{\sharp\sharp}$-modules, respectively.

Let $M\in{\rm Mod}_{ss}({\mathcal
  H}^{\sharp\sharp})$ with ${\rm dim}_k(M)<\infty$. By Proposition \ref{jishnusharp} we have $\Delta(M)\in {\rm Mod}^{\rm ad}({\mathfrak O})$, thus $\Delta(M)^*\otimes_{k[[t]]}k((t))\in{\rm Mod}^{et}(k((t)))$ (see Proposition \ref{nopsi}). Let $V(M)$ be the object in ${\rm Rep}({\rm Gal}(\overline{F}/F))$ assigned to $\Delta(M)^*\otimes_{k[[t]]}k((t))$ by Theorem \ref{sosego}. Exhausting an object in ${\rm Mod}_{ss}({\mathcal
  H}^{\sharp\sharp})$ by its finite dimensional subobjects we see that this
construction extends to all of ${\rm Mod}_{ss}({\mathcal
  H}^{\sharp\sharp})$.

\begin{satz}\label{7vorallhe} The assignment\begin{gather}{\rm Mod}_{ss}({\mathcal
  H}^{\sharp\sharp})\longrightarrow {\rm Rep}({\rm
      Gal}(\overline{F}/F)),\quad\quad M\mapsto
    V(M)\label{torschluss}\end{gather}is an exact contravariant functor, with
  ${\rm dim}_k(M)={\rm dim}_k(V(M))$ for any $M$. Also,\begin{gather}{\rm Mod}_{ss}({\mathcal H}^{\sharp})\longrightarrow {\rm Rep}({\rm Gal}(\overline{F}/F)),\quad\quad M\mapsto V(M),\notag\\{\rm Mod}_{ss}({\mathcal H})\longrightarrow {\rm Rep}({\rm Gal}(\overline{F}/F)),\quad\quad M\mapsto V(M)\label{torschluss1}\end{gather}are exact and fully faithful contravariant functors. 
\end{satz}

{\sc Proof:} Exactness follows from exactness of $M\mapsto \Delta(M)$
(Proposition \ref{jishnusharp}), exactness of $\Delta\mapsto
\Delta^*\otimes_{k[[t]]}k((t))$ (Proposition \ref{nopsi}) and exactness of the
equivalence functor in Theorem \ref{sosego}. From Proposition \ref{jishnusharp} we get
${\rm dim}_k(M)={\rm dim}_{k((t))}(\Delta(M)^*\otimes_{k[[t]]}k((t)))$, from
Theorem \ref{sosego} we get ${\rm dim}_{k((t))}(\Delta(M)^*\otimes_{k[[t]]}k((t)))={\rm dim}_k(V(M))$. 

To prove faithfulness on ${\rm Mod}_{ss}({\mathcal H}^{\sharp})$, suppose that we are given finite dimensional objects $M$, $M'$ in ${\rm
  Mod}_{ss}({\mathcal H}^{\sharp})$ and a morphism $\mu:V(M')\to V(M)$ in ${\rm Rep}({\rm Gal}(\overline{F}/F))$. By Theorem \ref{sosego}, the latter corresponds to a unique morphism of \'{e}tale $(\varphi,\Gamma)$-modules $$\mu:\Delta(M')^*\otimes_{k[[t]]}k((t)) \to\Delta(M)^*\otimes_{k[[t]]}k((t)).$$By Proposition \ref{wienaus} (which applies since Proposition \ref{jishnusharp} tells us $\Delta(M), \Delta(M')\in {\rm Mod}^{\clubsuit}({\mathfrak O})$) it is induced by a unique morphism of ${\mathfrak O}$-modules $\mu:\Delta(M)\to\Delta(M')$. Clearly $\mu$ takes ${\Delta(M)}[t]$ to ${\Delta(M')}[t]$, i.e. it takes
$M$ to $M'$. Applying Theorem \ref{shmichaelmas} to both $M$ and $M'$ we see that $\mu:M\to M'$ is ${\mathcal
  H}^{\sharp}$-equivariant. If $M, M'\in{\rm Mod}_{ss}({\mathcal H}^{\sharp})$ are not necessarily finite dimensional, the same conclusion is obtained by exhausting $M$, $M'$ by its finite dimensional submodules. We deduce the stated full faithfulness on ${\rm Mod}_{ss}({\mathcal H}^{\sharp})$. It implies full faithfulness on ${\rm Mod}_{ss}({\mathcal H})$ (see Lemma \ref{vialsa}).\hfill$\Box$\\

{\bf Example:} The analogs of Proposition
\ref{schutzengel} and Theorem \ref{7vorallhe} (on the functor in formula (\ref{torschluss1})) fail for supersingular ${\mathcal
  H}^{\sharp\sharp}$-modules. To see this, take $d=2$, and endow the $6$-dimensional
$k$-vector space $M$ with basis $e_0, e_1, e_2, f_0, f_1, f_2$ with the
structure of an ${\mathcal
  H}^{\sharp\sharp}$-module as follows. $T_{t}$ for each $t\in
\overline{T}$ acts trivially. Put $T_s(f_0)=T_s(e_1)=T_s(e_2)=0$ and
$T_s(e_0)=-e_0$, $T_s(f_1)=-f_1$, $T_s(f_2)=-f_2$. Fix $\alpha\in k$ and put
$T_{\omega}(e_0)=e_1$, $T_{\omega}(e_1)=e_2$, $T_{\omega}(e_2)=e_0$,
$T_{\omega}(f_0)=f_1+\alpha e_1$, $T_{\omega}(f_1)=f_2-\alpha e_2$,
$T_{\omega}(f_2)=f_0$. This is even an ${\mathcal
  H}^{\sharp}$-module if and only if $\alpha=0$, if and only if it is decomposable (as an ${\mathcal
  H}^{\sharp\sharp}$-module). The corresponding ${\mathfrak O}$-module
$\Delta(M)$ is
defined by the relations $\varphi e_0=-e_1$, $\varphi e_1=-e_2$,
$t^{q-1}\varphi e_2=-e_0$, $\varphi f_2=-f_0$, $t^{q-1}\varphi(f_0-\alpha
e_0)-f_1$, $t^{q-1}\varphi(f_1+\alpha
e_1)-f_2$. But this ${\mathfrak O}$-module is in fact independent of
$\alpha$, since $t^{q-1}\varphi  e_1=t^{q-1}\varphi e_0=0$. Thus, an ${\mathcal
  H}^{\sharp\sharp}$-analog of
Theorem \ref{7vorallhe} fails. To see that an ${\mathcal
  H}^{\sharp\sharp}$-analog of Proposition
\ref{schutzengel} fails take $M_0$ to be the $k$-subvector space of $M$
spanned by $e_0, e_1, e_2$; it is stable under ${\mathcal
  H}^{\sharp\sharp}$. The action of ${\mathcal
  H}^{\sharp\sharp}$ on $M_0$ and on $M/M_0$ does not depend on $\alpha$. The
actions of $T_{\omega}^{d+1}=T_{\omega}^3$, of $T_{e^*(\Gamma)}$ and of
$T_sT_{\omega}$ do not depend on $\alpha$. We have
$(T_sT_{\omega})^{-1}(M_0)=M_0+kf_2$ and hence the restriction of $T_{\omega}$
to $(T_sT_{\omega})^{-1}(M_0)$ does not depend on $\alpha$. We have
$M=\sum_{\epsilon}M^{\epsilon}$ with $M^{\epsilon}=0$ whenever $\epsilon_1\ne 0$. Thus, an ${\mathcal
  H}^{\sharp\sharp}$-analog of Proposition
\ref{schutzengel} would predict that also the action of $T_{\omega}$ (even of ${\mathcal
  H}^{\sharp\sharp}$) is independent
of $\alpha$, which however is apparantly not the case.

\subsection{The essential image}

\label{essim}

{\bf Definition:} Let ${\rm Hom}(\Gamma,k^{\times})^{\Phi}$ denote the group of $(d+1)$-tuples $\alpha=(\alpha_0,\ldots,\alpha_d)$ of characters $\alpha_j:\Gamma\to k^{\times}$. Let ${\mathfrak S}_{d+1}$ act on ${\rm Hom}(\Gamma,k^{\times})^{\Phi}$ by the formulae$$(\omega \cdot\alpha)_0=\alpha_d\quad\mbox{ and }\quad(\omega \cdot\alpha)_i=\alpha_{i-1}\mbox{ for }1\le i\le d,$$$$(s \cdot\alpha)_0=\alpha_1,\quad (s \cdot\alpha)_1=\alpha_0 \quad\mbox{ and }\quad(s \cdot\alpha)_i=\alpha_{i}\mbox{ for }2\le i\le d.$$Recall the action of ${\mathfrak S}_{d+1}$ on $[0,q-2]^{\Phi}$. Combining both (diagonally), we obtain an action of ${\mathfrak S}_{d+1}$ on ${\rm Hom}(\Gamma,k^{\times})^{\Phi}\times [0,q-2]^{\Phi}$. 

In Lemma \ref{ossa} we attached to each standard cyclic \'{e}tale $(\varphi,\Gamma)$-module ${\bf D}$ of dimension $d+1$ an ordered tuple
  $((\alpha_0,m_0),\ldots,(\alpha_d, m_d))$ (with integers $m_j\in[1-q,0]$ and characters $\alpha_j:\Gamma\to k^{\times}$), unique up to a cyclic permutation. Sending each $m_j$ to the representative in $[0,q-2]$ of its class in ${\mathbb Z}/(q-1)$, the tuple $(m_0,\ldots,m_d)$ gives rise to an element in $[0,q-2]^{\Phi}$. On the other hand, the tuple $(\alpha_0,\ldots,\alpha_d)$ constitutes an element in ${\rm Hom}(\Gamma,k^{\times})^{\Phi}$. Taken together we thus attach to ${\bf D}$ an element in ${\rm Hom}(\Gamma,k^{\times})^{\Phi}\times [0,q-2]^{\Phi}$, unique up to cyclic permutation. Equivalently, we attach to ${\bf D}$ an orbit in ${\rm Hom}(\Gamma,k^{\times})^{\Phi}\times [0,q-2]^{\Phi}$ under the action of the subgroup of ${\mathfrak S}_{d+1}$ generated by $\omega$. 

Now let ${\bf D}'_1$, ${\bf D}'_2$ be irreducible \'{e}tale $(\varphi,\Gamma)$-modules over
$k((t))$. We say that ${\bf D}'_1$, ${\bf D}'_2$ are strongly ${\mathfrak S}_{d+1}$-linked if they are subquotients of $(d+1)$-dimensional standard cyclic \'{e}tale $(\varphi,\Gamma)$-modules ${\bf D}_1$, ${\bf D}_2$ respectively, and if ${\bf D}_1$, ${\bf D}_2$ give rise to the same ${\mathfrak S}_{d+1}$-orbit in ${\rm Hom}(\Gamma,k^{\times})^{\Phi}\times [0,q-2]^{\Phi}$. We say that 
${\bf D}'_1$, ${\bf D}'_2$ are ${\mathfrak S}_{d+1}$-linked if they are subquotients of $(d+1)$-dimensional standard cyclic \'{e}tale $(\varphi,\Gamma)$-modules ${\bf D}_1$, ${\bf D}_2$ respectively, and if ${\bf D}_1$, ${\bf D}_2$ give rise to the same ${\mathfrak S}_{d+1}$-orbit in $[0,q-2]^{\Phi}$ (or equivalently, if the assigned tuples (up to cyclic permutation) in $[0,q-2]^{\Phi}$ coincide as {\it unordered} tuples (with multiplicities)).\\

{\bf Remark:} (a) Let ${\bf D}$ denote the \'{e}tale $(\varphi,\Gamma)$-module over
$k((t))$ corresponding to $V(M)$, for a finite dimensional supersingular ${\mathcal
    H}^{\sharp\sharp}$-module $M$. Our constructions show $M={\rm Hom}_k^{\rm cont}({\bf
  D}^{\natural},k)[t]$ (where ${\bf
  D}^{\natural}$ is given the $t$-adic topology). Moreover:

{\bf (i)} Consider the natural map of $k[[t]][\varphi]$-modules$$\kappa_{{\bf D}}:k[[t]][\varphi]\otimes_{k[[t]]}M\longrightarrow {\rm Hom}_k^{\rm cont}({\bf
  D}^{\natural},k).$$As a
$k[[t]][\varphi]$-module, ${\rm ker}(\kappa_{{\bf D}})$ is generated by ${\rm ker}(\kappa_{{\bf
    D}})\cap (k\otimes M+k[[t]]\varphi\otimes M)$.

{\bf (ii)} Each irreducible subquotient of ${\bf D}$ is a subquotient of a
$(d+1)$-dimensional standard cyclic \'{e}tale $(\varphi,\Gamma)$-module; more precisely: 

{\bf (ii)(1)} If ${\bf D}$ (or equivalently, $M$) is indecomposable, then any two irreducible subquotients of ${\bf D}$ are ${\mathfrak S}_{d+1}$-linked. 

{\bf (ii)(2)} If $M$ is even a supersingular ${\mathcal
    H}$-module, and if ${\bf D}$ (or equivalently, $M$) is indecomposable,
then any two irreducible subquotients of ${\bf D}$ are strongly ${\mathfrak
  S}_{d+1}$-linked.

{\bf (ii)(3)} If $M$ is even a supersingular ${\mathcal H}^{\flat}$-module, then each irreducible subquotient of ${\bf D}$ is a subquotient of a $(d+1)$-dimensional standard cyclic \'{e}tale $(\varphi,\Gamma)$-module with parameters $m_j\in\{1-q,0\}$ and $\alpha_j=1$ for all $j$.

{\bf (iii)} For any $(\varphi,\Gamma)$-submodule ${\bf D}_0$ of ${\bf D}$ the
$\psi$-operator on ${\bf D}_0\cap {\bf D}^{\natural}$ is surjective.

(b) Does property {\bf (i)} mean (at
least if property {\bf (iii)} is assumed) that ${\bf D}$ is the reduction of a crystalline $p$-adic ${\rm
  Gal}(\overline{F}/F)$-representation with Hodge-Tate weights in $[-1,0]$ ?

(c) Property {\bf (iii)} means that the functor ${\bf D}_0\mapsto {\bf
  D}_0^{\natural}$ is exact on the category of subquotients ${\bf D}_0$ of ${\bf D}$.

(d) It should not be too hard to show that properties {\bf (i)}, {\bf
  (ii)(1)} and {\bf (iii)} together in fact {\it characterize} the essential image of the
functor (\ref{torschluss}).

(e) On the other hand, properties {\bf (i)}, {\bf
  (ii)(2)} and {\bf (iii)} together do {\it not} characterize the essential image of the
functor (\ref{torschluss1}). To see this for $d=1$
consider the following \'{e}tale
$(\varphi,\Gamma)$-module ${\bf D}$ (which satisfies {\bf (i)}, {\bf (ii)(2)}, {\bf
  (iii)}). It is given by a $k$-basis $e_0$, $e_1$, $f_0$, $f_1$, $g_0$,
$g_1$ of $({\bf
  D}^{\natural})^*[t]$ and the following relations:$$\varphi e_1=e_0,\quad \varphi
f_1=f_0,\quad \varphi g_1=g_0,\quad
t^{q-1}\varphi e_0=e_1,\quad t^{q-1}\varphi f_0=f_1+e_1,\quad t^{q-1}\varphi
g_0=g_1+f_0.$$Another object not in the essential image is defined by the set of relations$$\varphi e_1=e_0,\quad \varphi
f_1=f_0,\quad \varphi g_1=g_0,\quad 
t^{q-1}\varphi e_0=e_1,\quad t^{q-1}\varphi f_0=f_1+e_0,\quad t^{q-1}\varphi
g_0=g_1+f_1.$$
 
\section{From $G$-representations to ${\mathcal H}$-modules}

\subsection{Supersingular cohomology}

\label{susicoho}

Put $G={\rm GL}_{d+1}(F)$, let $I_0$ be a pro-$p$-Iwahori subgroup in $G$, and
fix an isomorphism between ${\mathcal H}$ and the pro-$p$-Iwahori Hecke
algebra $k[I_0\backslash G/I_0]$ corresponding to $I_0\subset G$. For a smooth
$G$-representation $Y$ (over $k$) the subspace $Y^{I_0}$ of $I_0$-invariants
then receives a natural action by ${\mathcal H}$. Let us denote by
$H^0_{ss}(I_0,Y)$ the maximal supersingular ${\mathcal H}$-submodule of
$Y^{I_0}$. It is clear that this defines a left exact functor$${\rm
  Mod}(G)\longrightarrow {\rm Mod}_{ss}({\mathcal H}),\quad\quad Y\mapsto
H^0_{ss}(I_0,Y)$$where ${\rm Mod}(G)$ denotes the category of smooth
$G$-representations. The category ${\rm Mod}(G)$ is a Grothendieck category
(\cite{peterdga} Lemma 1) and has enough injective objects (\cite{vigbuch}
I.5.9). Let $D^+(G)$ denote the derived category of complexes of
smooth $G$-representations vanishing in negative degrees, let
$D_{ss}^+({\mathcal H})$ denote the derived category of complexes of
supersingular ${\mathcal H}$-modules vanishing in negative degrees. The above
functor gives rise to a right derived
functor\begin{gather}R_{ss}(I_0,.):D^+(G)\longrightarrow D_{ss}^+({\mathcal
  H}).\label{allesee}\end{gather}Let $D^+({\rm Gal}(\overline{F}/F)))$ denote
the derived category of complexes in ${\rm Rep}({\rm Gal}(\overline{F}/F))$
vanishing in negative degrees. Since the functor $V$ is exact, it induces a
functor$$V:D_{ss}^+({\mathcal H})\longrightarrow D^+({\rm
  Gal}(\overline{F}/F))).$$We may compose them with $R_{ss}(I_0,.)$ to obtain a functor$$V\circ R_{ss}(I_0,.):D^+(G)\longrightarrow D^+({\rm Gal}(\overline{F}/F))).$$

{\bf Remark:} The functor $H^0_{ss}(I_0,.)$ is the composite of the left exact
functor ${\rm Mod}(G)\to {\rm Mod}({\mathcal H})$, $Y\mapsto Y^{I_0}$
(taking $I_0$-invariants) and the left exact functor ${\rm Mod}({\mathcal H})\to {\rm
  Mod}_{ss}({\mathcal H})$, $M\mapsto M_{ss}$ which takes an ${\mathcal
  H}$-module to its maximal supersingular ${\mathcal H}$-submodule. Also ${\rm
  Mod}({\mathcal H})$ is a Grothendieck category with enough injective
objects. Writing $R(I_0,.)$ and $R_{ss}(.)$ for the respective right derived
functors, we have a morphism $R_{ss}(I_0,.)\to R_{ss}(.)\circ R(I_0,.)$.\\

{\bf Remark:} Of course, we expect the functor $V\circ R_{ss}(I_0,.)$ to be meaningful only when restricted to (complexes of) supersingular $G$-representations. The reason is the following theorem of Ollivier and Vign\'{e}ras \cite{ollvig17}: A smooth admissible irreducible $G$-representation $Y$ over an algebraic closure $\overline{k}$ of $k$ is supersingular if and only if $Y^{I_0}$ is a supersingular ${\mathcal H}\otimes_k{\overline{k}}$-module, if and only if $Y^{I_0}$ admits a supersingular subquotient.

It is known that, beyond the case where $G={\rm GL}_2({\mathbb Q}_p)$, a smooth admissible irreducible supersingular $G$-representation $Y$ over $k$ is not uniquely determined by the ${\mathcal H}$-module $Y^{I_0}$. Is it perhaps uniquely determined by the derived object $R_{ss}(I_0,Y)\in D^+_{ss}({\mathcal H})$ ? It would then also be uniquely determined by the derived object ${V}(R_{ss}(I_0,Y))  \in D^+({\rm Gal}(\overline{F}/F)))$.\\

{\bf Remark:} For the universal module $Y={\rm ind}_{I_0}^G k$ we have $H^0_{ss}(I_0,Y)=0$ since ${\mathcal
  H}=({\rm ind}_{I_0}^G k)^{I_0}$ does not contain non-zero finite dimensional ${\mathcal
  H}$-submodules (let alone supersingular ones).

\subsection{An exact functor from $G$-representations to ${\mathcal H}$-modules}

We fix a $(d+1)$-st root of unity $\xi\in k^{\times}$ with
$\sum_{j=0}^{d}\xi^j=0$. 

For an ${\mathcal H}$-module $M$ and $j\in{\mathbb Z}$ let $M^{{\xi}^j}$ be the ${\mathcal H}$-module which coincides with $M$ as a module over the $k$-subalgebra $k[T_s, T_{t}]_{t\in\overline{T}}$, but with $T_{\omega}|_{M^{{\xi}^j}}={\xi}^j T_{\omega}|_{M}$.

Let $\delta:M_0\to M_1$ be a morphism of ${\mathcal
  H}$-modules. For $(x_0,x_1)\in M_0\oplus M_1$ put\begin{align}T_{\omega}((x_0,x_1))&=(T_{\omega}(x_0),T_{\omega}(\delta(x_0))+\xi T_{\omega}(x_1)),\notag\\T_s((x_0,x_1))&=(T_s(x_0),T_s(x_1)),\notag\\T_{t}((x_0,x_1))&=(T_{t}(x_0),T_{t}(x_1))\quad\mbox{ for }t\in\overline{T}.\notag\end{align}

\begin{lem}\label{14vorstubei} These formulae define an ${\mathcal
    H}$-module structure on $M_0\oplus M_1$; we denote this new
  ${\mathcal H}$-module by $M_0\oplus^{\delta} M_1$. We have an exact
  sequence of ${\mathcal H}$-modules \begin{gather}0\longrightarrow M_1^{\xi}\longrightarrow
  M_0\oplus^{\delta} M_1\longrightarrow
  M_0\to0.\label{lotharaugust}\end{gather}The morphism $\delta:M_0\to M_1$ can be recovered from the exact sequence (\ref{lotharaugust}). 

If there is some $\lambda\in k^{\times}$ with $T_{\omega}^{d+1}=\lambda$ on $M_0$ and on $M_1$, then also $T_{\omega}^{d+1}=\lambda$ on $M_0\oplus^{\delta} M_1$,
\end{lem}

{\sc Proof:} By induction on $i$ one
shows $$T_{\omega}^i((x_0,x_1))=(T_{\omega}^i(x_0),\xi^iT_{\omega}^i(x_1)
+\sum_{j=0}^{i-1}\xi^jT_{\omega}^i(\delta(x_0)))$$for $i>0$, and hence
$T_{\omega}^{d+1}((x_0,x_1))=(T_{\omega}^{d+1}(x_0),T_{\omega}^{d+1}(x_1))$. From
here, all the required relations are straighforwardly verified, showing that
indeed we have defined an ${\mathcal H}$-module.

Obviously, from the exact sequence (\ref{lotharaugust}) both $M_0$ and $M_1$
can be recovered. That also $\delta$ can be recovered follows from the
following more general consideration. Suppose that we are given $\delta:M_0\to M_1$ and
$\epsilon:N_0\to N_1$ and a morphism of ${\mathcal H}$-modules
$f:M_0\oplus^{\delta} M_1\to N_0\oplus^{\epsilon} N_1$ with
$f(M_1^{\xi})\subset N_1^{\xi}$. Then there are ${\mathcal H}$-module
homomorphisms $f_0:M_0\to N_0$, $f_1:M_1^{\xi}\to N_1^{\xi}$ and
$\tilde{f}:M_0\to N_1^{\xi}$ with
$f((x_0,x_1))=(f_0(x_0),f_1(x_1)+\tilde{f}(x_0))$. For $x_0\in M_0$ we
compute $$f(T_{\omega}(x_0,0))=f(T_{\omega}(x_0),T_{\omega}(\delta(x_0)))=(T_{\omega}(f_0(x_0)),T_{\omega}(f_1(\delta(x_0)))+\xi
T_{\omega}(\tilde{f}(x_0))),$$$$T_{\omega}(f(x_0,0))=T_{\omega}(f_0(x_0),\tilde{f}(x_0))=(T_{\omega}(f_0(x_0)),T_{\omega}(\epsilon(f_0(x_0)))+\xi T_{\omega}(\tilde{f}(x_0))).$$As
$f(T_{\omega}(x_0,0))=T_{\omega}(f(x_0,0))$ we deduce
$T_{\omega}(\epsilon(f_0(x_0)))=T_{\omega}(f_1(\delta(x_0)))$, and since
$T_{\omega}$ is an isomorphism even $\epsilon(f_0(x_0))=f_1(\delta(x_0))$.\hfill$\Box$\\

Let $$(M_{\bullet},\delta_{\bullet})=[\ldots
  \stackrel{\delta_{-2}}{\longrightarrow} M_{-1} \stackrel{\delta_{-1}}{\longrightarrow} M_0\stackrel{\delta_0}{\longrightarrow} M_1\stackrel{\delta_1}{\longrightarrow}M_2\stackrel{\delta_2}{\longrightarrow} \ldots]$$be a complex of ${\mathcal H}$-modules. 

\begin{lem}\label{claudeles} (a) There is a unique ${\mathcal
    H}$-module $\oplus^{\delta_{\bullet}}_{j\in\mathbb Z}M_j$ with the following properties:

${\bullet}$ As a $k$-vector space, $\oplus^{\delta_{\bullet}}_{j\in\mathbb Z}M_j=\oplus_{j\in\mathbb Z}M_j$.

${\bullet}$ For any $j$ we have $\tau(M_j)\subset M_j+M_{j+1}$ for each $\tau\in{\mathcal H}$; in particular, the subspace $M_{\ge j}=\oplus_{j'\ge j}M_{j'}$ is an ${\mathcal H}$-submodule. 

${\bullet}$ The ${\mathcal H}$-module $M_{\ge j}/M_{\ge j+2}$ is isomorphic with $M_j^{\xi^j}\oplus^{\delta_j} M_{j+1}^{\xi^{j}}$ as defined in Lemma \ref{14vorstubei}.

(b) If there is some $\lambda\in k^{\times}$ with $T_{\omega}^{d+1}=\lambda$
  on each $M_j$, then $T_{\omega}^{d+1}=\lambda$ on
  $\oplus^{\delta_{\bullet}}_{j\in\mathbb Z}M_j$.

(c) The assignment $(M_{\bullet},\delta_{\bullet})\mapsto
  (\oplus^{\delta_{\bullet}}_{j\in\mathbb Z}M_j,(M_{\ge j})_{j\in{\mathbb
      Z}})$ is an exact and faithful functor
  from the category of complexes of ${\mathcal H}$-modules to the
  category of filtered ${\mathcal H}$-modules. The isomorphism class of
  the complex
  $(M_{\bullet},\delta_{\bullet})$ can be recovered from the isomorphism class
  of the filtered ${\mathcal H}$-module $(\oplus^{\delta_{\bullet}}_{j\in\mathbb Z}M_j,(M_{\ge j})_{j\in{\mathbb
      Z}})$. 
\end{lem}

{\sc Proof:} This is clear from Lemma \ref{14vorstubei}.\hfill$\Box$\\

{\bf Definition:} (a) For a smooth $G$-representation $Y$ over $k$ and $i\ge0$ let us denote by $H_{ss}^i(I_0,Y)$ the $i$-th cohomology group of $R_{ss}(I_0,Y)$, see formula (\ref{allesee}). 

(b) We say that a smooth $G$-representation $Y$ over $k$ is {\it exact} if for each $i\ge0$ the functor $Y'\mapsto H_{ss}^i(I_0,Y')$ is exact on the category of $G$-subquotients $Y'$ of $Y$. 

(c) An exhaustive and
separated decreasing filtration $({Y}^{j})_{j\in{\mathbb Z}}$ of a smooth
$G$-representation $Y$ over $k$ is {\it exact} if ${Y}^{j}/{Y}^{j+1}$ is exact
for each $j$. \\

{\bf Example:} A semisimple smooth $G$-representation is exact. \\

Let ${\mathfrak R}_G$ denote the following category: objects are smooth
$G$-representations with an exact filtration, morphisms are
$G$-equivariant maps respecting the filtrations (i.e. $f:Y\to W$ with
$f(Y^i)\subset W^i$ for all $i$). We denote objects $(Y,({Y}^{i})_{i\in{\mathbb Z}})$ in ${\mathfrak R}_G$
simply by $Y^{\bullet}$. 

Let ${\mathfrak E}({\mathcal H})$ denote the category of $E_1$-spectral sequences in the category of ${\mathcal H}$-modules.

For $Y^{\bullet}\in {\mathfrak R}_G$ we have the spectral
sequence$$E(Y^{\bullet})=[E_1^{m,n}(Y^{\bullet})=H^{m+n}_{ss}(I_0,{Y}^{m}/{Y}^{m+1})\Rightarrow
H_{ss}^{m+n}(I_0,Y)].$$A morphism $f:Y^{\bullet}\to W^{\bullet}$
in ${\mathfrak R}_G$ induces morphisms 
$H^{m}_{ss}(I_0,{Y}^{i}/{Y}^{i+1})\to H^{m}_{ss}(I_0,{W}^{i}/{W}^{i+1})$
for any $m$ and $i$, and these induce a morphism of spectral
sequences $E(Y^{\bullet})\to E(W^{\bullet})$. We thus obtain a functor $${\mathfrak R}_G\to {\mathfrak E}({\mathcal H}),\quad\quad Y^{\bullet}\mapsto E(Y^{\bullet}).$$

For $r\ge1$ let ${\mathcal Y}_r$ be the set of equivalence classes of pairs of
integers $(m,n)$, where $(m,n)$ is declared to be equivalent with $(m',n')$ if
and only if there is some $j\in{\mathbb Z}$ with $(m,n)=(m'+jr,n'-j(r-1))$. For
$y\in {\mathcal Y}_r$ let $E_r^y(Y^{\bullet})$ be the complex of ${\mathcal H}$-modules whose terms are the $E_r^{m,n}(Y^{\bullet})$ with $(m,n)\in y$, and whose
differentials $d_r:E_r^{m,n}(Y^{\bullet})\to E_r^{m+r,n-r+1}(Y^{\bullet})$ are given by the spectral sequence. We apply the functor of
Lemma \ref{claudeles} to $E_r^y(Y^{\bullet})$ to obtain a (filtered) supersingular ${\mathcal H}$-module ${\bf E}_r^y(Y^{\bullet})$.

For a morphism $f:Y^{\bullet}\to W^{\bullet}$
in ${\mathfrak R}_G$ we have induced ${\mathcal H}$-linear maps $f_r:\oplus_{y\in {\mathcal Y}_r}{\bf
  E}_r^y(Y^{\bullet})\to \oplus_{y\in {\mathcal Y}_r}{\bf E}_r^y(
W^{\bullet})$. Notice however that, in general, for a given $y\in{\mathcal
  Y}_r$ there is no $y'\in {\mathcal Y}_r$ such that $f_r({\bf
  E}_r^y(Y^{\bullet}))\subset{\bf
  E}_r^{y'}(W^{\bullet})$, even if $r=1$.

\begin{lem}\label{13forstubai} Let $Y^{\bullet}\to W^{\bullet}\to
  X^{\bullet}$ be a complex in
  ${\mathfrak R}_G$ such that for each $i$ the induced sequence $0\to Y^i/Y^{i+1}\to W^i/W^{i+1}\to
  X^{i}/X^{i+1}\to 0$ is exact. We then have an exact sequence of supersingular ${\mathcal H}$-modules$$0\longrightarrow\bigoplus_{y\in {\mathcal Y}_1}{\bf E}_1^y(
Y^{\bullet})\longrightarrow\bigoplus_{y\in {\mathcal Y}_1}{\bf E}_1^y(
W^{\bullet})\longrightarrow\bigoplus_{y\in {\mathcal Y}_1}{\bf E}_1^y(
X^{\bullet})\longrightarrow0.$$
\end{lem}

{\sc Proof:} This follow from the constructions.\hfill$\Box$\\

{\bf Remark:} The analog of Lemma \ref{13forstubai} is false for the maps $f_r$ for $r>1$.\\

{\bf Remark:} For a smooth $G$-representation $Y$ endowed with an exact filtration, we may apply the functor $V$ of subsection \ref{fullfafu} to the supersingular ${\mathcal H}$-module ${\bf E}_r^y(
Y^{\bullet})$ (any $r$). In this way, we assign a ${\rm Gal}(\overline{F}/F)$-representation to $Y$. We propose this construction as a non-derived alternative to that of subsection \ref{susicoho}. Of course, again it will be meaningful only on supersingular $G$-representations.

We expect that for $G={\rm GL}_2({\mathbb Q}_p)$, this construction, with $r=1$, essentially recovers the restriction of Colmez's functor to all \footnote{i.e. not only to those generated by their $I_0$-invariants} supersingular $G$-representations.

\begin{flushleft} \textsc{Humboldt-Universit\"at zu Berlin\\Institut f\"ur Mathematik\\Rudower Chaussee 25\\12489 Berlin, Germany}\\ \textit{E-mail address}:
gkloenne@math.hu-berlin.de \end{flushleft} 
\begin{thebibliography}{abcdefgh} 

\bibitem{col}{\it P. Colmez}, $(\varphi,\Gamma)$-modules et
  repr\'{e}sentations du mirabolique de ${\rm GL}_2({\mathbb Q}_p)$,
  Ast\'{e}risque {\bf 330} (2010), 61--153.

\bibitem{em}{\it M. Emerton}, On a class of coherent rings, with applications
  to the smooth representation theory of ${\rm GL}_2({\mathbb Q}_p)$ in
  characteristic $p$, preprint 2008
  
\bibitem{fon}{\it J.-M. Fontaine}, Repr\'{e}sentations $p$-adiques des corps locaux. I, The Grothendieck
Festschrift, Vol. II, Progr. Math. {\bf 87}, Birkh¨auser Boston, Boston, MA, 1990, p. 249–309.

\bibitem{multivar}{\it E. Grosse-Kl\"onne}, A note on multivariable $(\varphi,\Gamma)$-modules, Research in Number Theory (2019) 5: 6
  
\bibitem{dfun}{\it E. Grosse-Kl\"onne}, From pro-$p$-Iwahori Hecke modules to
  $(\Phi,\Gamma)$-modules I, Duke Mathematical J. {\bf  165}, no. 8 (2016), 1529--1595 

\bibitem{kire}{\it M. Kisin, W. Ren}, Galois representations and Lubin-Tate groups, Doc. Math.
{\bf 14} (2009), p. 441–461.

\bibitem{ol}{\it R. Ollivier}, Parabolic Induction and Hecke modules in
  characteristic $p$ for $p$-adic ${\rm GL}_n$, Algebra and Number Theory
  {\bf 4}(6) (2010), 701--742.

\bibitem{ollvig17}{\it R. Ollivier, M.-F. Vigneras}, Parabolic induction in characteristic $p$, Selecta Math. (N.S.) {\bf 24}(5) (2018), 3973 --– 4039.

\bibitem{peterdga}{\it Smooth representations and Hecke modules in
  characteristic $p$}, Pacific J. Math. {\bf  279} (2015), 447 --- 464
  
\bibitem{peterlec}{\it P. Schneider}, Galois representations and $(\varphi,\Gamma)$-modules, course at M\"unster in 2015, Cambridge Studies in Advanced Mathematics (2017)

\bibitem{schven}{\it P. Schneider, O. Venjakob}, Coates-Wiles homomorphisms
  and Iwasawa cohomology for Lubin-Tate extensions, in: Elliptic Curves, Modular Forms and Iwasawa Theory. Springer Proc. Math. Stat. {\bf 188}, 401--468 (2016) 

 \bibitem{vigbuch}{\it M.F. Vign\'{e}ras}, Repr\'{e}sentations
  $\ell$-modulaires d'un
  groupe  r\'{e}ductif $p$--adique avec $\ell\ne p$, Progress in Math. {\bf 137}, Birkh\"auser: Boston 1996
  
\bibitem{vigneras}{\it M.F. Vign\'{e}ras}, Pro-$p$-Iwahori Hecke ring and supersingular $\overline{\bold F}\sb p$-representations,  Math. Ann.  {\bf 331}  (2005),  no. 3, 523--556 and Math. Ann. {\bf  333}  (2005),  no. 3, 699--701. 

\bibitem{vigjuss}{\it M.F. Vign\'{e}ras}, The pro-$p$-Iwahori Hecke algebra of afddc reductive $p$-adic group III, Journal of the Institute of Mathematics of Jussieu {\bf 16}(3), 571–-608 (2017)

\end{thebibliography}
\end{document}